\documentclass[11pt]{article}
\usepackage{amssymb}
\usepackage{amsmath}
\usepackage{bm}
\newtheorem{lemma}{Lemma}[section]
\newtheorem{example}{Example}[section]
\newtheorem{theorem}{Theorem}[section]
\newtheorem{note}{Note}[section]

\begin{document}

\begin{center}
{\LARGE\bf Unbiased estimates for products}\\[1ex]
{\LARGE\bf of moments and cumulants}\\[1ex]
{\LARGE\bf for finite and infinite populations}\\[1ex]
by\\[1ex]
Christopher S. Withers, Industrial Research Limited, Lower Hutt, New Zealand\\
email: {\sf kit.withers@gmail.com}\\
Saralees Nadarajah, University of Manchester, Manchester M13 9PL, UK\\
email: {\sf mbbsssn2@manchester.ac.uk}
\end{center}
\vspace{1.5cm}
{\bf Abstract:}~~Let $F=F_N$ be the distribution of a finite real population of size $N$.
Let $\widehat{F}=F_N$ be the empirical distribution of a sample of size $n$ drawn
from the population without replacement.
We prove the following
remarkable {\it inversion principle} for obtaining unbiased estimates.
Let $ T \left(F_N\right)$ be any product of the moments or cumulants of $F_N$.
Let $T_{n, N} \left( F_N \right) = E T \left( F_n \right)$.
Then $E T_{N, n} \left( F_n \right) = T \left( F_N \right)$.
We also obtain an explicit expression for $T_{n, N} \left(F_N\right)$ for all $ T \left( F_N \right)$ of order up to 6.

We also prove the following related result.
If $F_n$  and $F_N$ are the sample and population distributions, the only functionals
for which $E T \left( F_n \right) = \lambda_{n, N} T \left( F_N \right)$ are noncentral moments, and
generalized second and third order central moments.
For these three cases the eigenvalues are $\lambda_{n, N}=1$, $\left( 1 - n^{-1} \right) \left( 1 - N^{-1} \right)^{-1}$,
and $\left( 1 - n^{-1} \right) \left( 1 - 2n^{-1} \right) \left( 1 - N^{-1} \right)^{-1} \left( 1 - 2N^{-1} \right)^{-1}$ respectively.



\section{Introduction}

Given a random sample of size $n$ without replacement say $X_1,\ldots,X_n$ from
a finite real population $x_1,\ldots,x_N$ with distribution $F(x)=N^{-1} \displaystyle \sum_{i=1}^N 1 \left( x_i\leq x \right)$,
mean $\mu=m_1 = N^{-1} \displaystyle \sum_{i=1}^N x_i$,
$r$th moment $m_r = \displaystyle N^{-1}\sum_{i=1}^Nx_i^r$,
$r$th central moment $\mu_r=\mu_r(F) = \displaystyle N^{-1}\sum_{i=1}^N \left(x_i-\mu\right)^r$
and $r$th cumulant  $\kappa_r$,
we obtain unbiased estimates (UEs) for products of them
\begin{eqnarray*}
&&
\displaystyle
m \left( 1^{r_1}2^{r_2} \cdots \right) = m_1^{r_1}m_2^{r_2} \cdots,
\nonumber
\\
&&
\displaystyle
\mu \left( 1^{r_1}2^{r_2} \cdots \right) = \mu^{r_1}\mu_2^{r_2}\mu_3^{r_3} \cdots,
\nonumber
\\
&&
\displaystyle
\kappa \left( 1^{r_1}2^{r_2} \cdots \right) = \kappa_1^{r_1}\kappa_2^{r_2} \cdots
\end{eqnarray*}
and for products of the joint central moments
\begin{eqnarray*}
\displaystyle
\mu \left( 1^{r_1}2^{r_2} \cdots \right) =
E \left( \widehat{\mu}-\mu \right)^{r_1} \left( \widehat{\mu}_2- E  \  \widehat{\mu}_2 \right)^{r_2} \cdots
\end{eqnarray*}
and for products of the corresponding joint cumulants $\kappa \left( 1^{r_1}2^{r_2} \cdots \right)$,
for {\it weight} $ \leq 6$, where the weight is
\begin{eqnarray*}
\displaystyle
r \left( 1^{r_1}2^{r_2} \cdots \right) =
1 \cdot r_1 + 2 \cdot r_2 + \cdots
\mbox{ that is } r({\bm \pi}) = \pi_1 + \pi_2 + \cdots
\end{eqnarray*}
for  ${\bm \pi} = \left( \pi_1, \pi_2, \ldots \right)$  a partition of $ r$.
We assume all partitions ${\bm \pi}$ are put into ascending order $\left( 1^{r_1}2^{r_2} \cdots \right)$.
For example, we write $\left( 1^22 \right)$ or $(112)$ rather than $\left( 21^2 \right)$ or $(121)$.

Our UEs are given in terms of
\begin{eqnarray*}
\displaystyle
\widehat{m} \left( 1^{r_1}2^{r_2} \cdots \right) = \widehat{m}_1^{r_1}\widehat{m}_2^{r_2} \cdots,
\mbox{ and }
\displaystyle
\widehat{\mu} \left( 1^{r_1}2^{r_2} \cdots \right) = \widehat{\mu}^{r_1}\widehat{\mu}_2^{r_2} \cdots,
\end{eqnarray*}
respectively,
where
\begin{eqnarray*}
\displaystyle
\widehat{\mu}=\widehat{m}_1=\overline{X}=n^{-1}\sum_{i=1}^nX_i,
\
\displaystyle
\widehat{m}_r=n^{-1}\sum_{i=1}^nX_i^r,
\
\displaystyle
\widehat{\mu}_r=\mu_r \left(\widehat{F}\right) = n^{-1}\sum_{i=1}^n \left(X_i-\overline{X}\right)^r
\end{eqnarray*}
and $\widehat{F}(x)=n^{-1} \displaystyle \sum_{i=1}^n 1 \left(X_i\leq x\right)$, the sample distribution.

For ${\bm \pi}$ a partition of $r$, Section  2 gives $E\widehat{m}({\bm \pi})$
and an UE of $m({\bm \pi})$ for $r \leq 6$.
Section 3 gives $E\widehat{\mu}({\bm \pi})$ and an UE of $\mu({\bm \pi})$ for $r \leq 6$.

We discover a remarkable {\it inversion principle}.
A result of this is that
these UEs do not require having to invert any matrices or solve any sets of
linear equations.
Set
\begin{eqnarray*}
&&
\displaystyle
{\bf m}_{(r)}= \left\{ m ({\bm \pi}): \ {\bm \pi} \mbox{ a partition of } r \right \},
\nonumber
\\
&&
\displaystyle
{\bm \mu}_{(r)} = \left\{ \mu ({\bm \pi}): {\bm \pi} \mbox{ a partition of } r \right\},
\nonumber
\\
&&
\displaystyle
{\bf K}_{(r)} = \left\{ K ({\bm \pi}): \ {\bm \pi} \mbox{ a partition of } r \right \}.
\end{eqnarray*}
In Section 2, we derive a matrix ${\bf B}_r = {\bf B}_r(N,n)$  from Skellam (1949) such that
$E \widehat{\bf m}_{(r)} = {\bf B}_r {\bf m}_{(r)}$, so
${\bf B}_r(N,n)^{-1} \widehat{\bf m}_{(r)}$ is an UE of ${\bf m}_{(r)}$.

In Section 3, we derive a matrix ${\bf C}_r = {\bf C}_r(N,n)$ From Sukhatme (1944) such that
$E \widehat{\bm \mu}_{(r)} = {\bf C}_r {\bm \mu}_{(r)}$, so
${\bf C}_r(N,n)^{-1}\widehat{\bm \mu}_{(r)}$ is an UE of ${\bm \mu}_{(r)}$.
So, expressing cumulants in terms of moments as ${\bf K}_{(r)} = {\bf G}_r {\bf m}_{(r)}$,
where ${\bf G}_r$ is a matrix of constants, it follows that
\begin{eqnarray}
\displaystyle
E \ \widehat{\bf K}_{(r)} = {\bf D}_r(N,n) {\bf K}_{(r)} \mbox{ where }
{\bf D}_r(N,n) = {\bf G}_r {\bf C}_r (N,n) {\bf G}_r^{-1}.
\label{1.8}
\end{eqnarray}
The inversion principle proved in Section 7 states the following amazing result.
\begin{eqnarray}
\displaystyle
{\bf B}_r(N,n)^{-1} = {\bf B}_r(n,N),
\
\displaystyle
{\bf C}_r(N,n)^{-1} = {\bf C}_r(n,N),
\
\displaystyle
{\bf D}_r(N,n)^{-1} = {\bf D}_r(n,N).
\label{1.9}
\end{eqnarray}
This implies that for $T (F)$ a product of moments or cumulants,
if $T_{n, N} (F) = E T \left( \widehat{F} \right)$
then $E T_{N, m} \left( \widehat{F} \right) = T (F)$.
In a later paper we shall extend this result to more general functionals.

Section  4 gives equivalent multivariate results to those of Section  3.
For partitions $\pi_1, \pi_2, \ldots$, set
\begin{eqnarray*}
\displaystyle
\mu \left(\pi_1,\pi_2, \ldots\right) = \mu \left(\pi_1\right) \mu \left(\pi_2\right) \cdots
\end{eqnarray*}
and \\
\begin{eqnarray*}
\displaystyle
\kappa \left( \pi_1, \pi_2, \ldots \right) = \kappa \left(\pi_1\right) \kappa \left(\pi_2\right) \cdots
\end{eqnarray*}

Section  5 gives $E\widehat{\mu} \left( \pi_1, \pi_2, \ldots \right)$ and an UE of $\mu \left( \pi_1, \pi_2, \ldots \right)$
up to total order 6, in particular UEs for $\mu \left(1^r\right)$.
MAPLE was used to
simplify the UE of $\mu \left(1^r\right)$ given by Dwyer and Tracy (1980) for $r \leq 5$
and to confirm they agreed with our results.
Earlier, Nath (1968,1969) gave $\mu \left(1^r\right)$ for $r=3,4$ and UEs for them.

Section  6 gives $E\widehat{\kappa} \left( \pi_1, \pi_2, \ldots \right)$ and an UE of
$\kappa \left( \pi_1, \pi_2, \ldots \right)$ up to total order 6.

Section 7 also proves a multivariate inversion principle: in this case the dimension
$n_r \times n_r$ in (1.7) jumps to $N_r \times N_r$, where
\begin{eqnarray*}
\displaystyle
n_r \mbox{ is the number of partitions of } r,
\
\displaystyle
N_r=\sum_r^{\bm \pi} P({\bm \pi})
\end{eqnarray*}
and $P({\bm \pi})$ is the partition function
\begin{eqnarray}
\displaystyle
P \left( 1^{r_1}2^{r_2} \cdots \right) = r!/ \prod_{i=1} \left( i!^{r_i}r_i! \right)
\label{1.11}
\end{eqnarray}
for $r = \displaystyle \sum_{i=1}ir_i$.

Section 8 proves the following related result.
The only functionals $T(F)$
for which
\begin{eqnarray*}
\displaystyle
\lambda_{n, N} = ET \left(F_n\right) / T(F) \mbox{ does not depend on } F
\end{eqnarray*}
are of the form
\begin{eqnarray*}
\displaystyle
T(F)=ES(X),
\
\displaystyle
\mu^t_{1, 2}(F),
\
\displaystyle
\mu^u_{1, 2, 3}(F),
\end{eqnarray*}
where
\begin{eqnarray*}
\displaystyle
\mu_{1, 2}^t(F)=E \left\{ t(X,X)-t(X,Y) \right\},
\
\displaystyle
\mu_{1, 2, 3}^u(F)=E \left\{ u(X,X,X)-u(X,Y,Y)-u(Y,X,Y)-u(Y,Y,X)+2u(X,Y,Z) \right\},
\end{eqnarray*}
where $X,Y,Z$ are independent with distribution $F$.
For these three cases $\lambda_{n, N}=1$, $\left( 1 - n^{-1} \right) / \left( 1 - N^{-1} \right)$ and
$\prod_{i=1}^2 \left\{ \left( 1 - in^{-1} \right) / \left( 1 - iN^{-1} \right) \right\}$.

Pierce (1940) gives transformations that allow one to obtain noncentral
moments up to order six for finite populations from those for infinite
populations.
He also covers the case where the sample values $X_1,\ldots,X_n$
may have different moments.

We shall use ${\bm \pi}$ for a partition of $r$, ${\bm \pi}_{-}$ for a partition of $r$
excluding one's say ${\bm \pi}_{-} = \left( 2^{r_2}3^{r_3} \cdots \right)$ with $2r_2+3r_3+ \cdots=r$,
and ${\bm \pi}_{+}$ for a partition of $r$ including at least one $1$.
Let $n_r$ be the number of partitions of $r$, and $N_r$ the sum of the partitions of $r$.
We also partition vectors and matrices using subscript $+$ and $--$.
For example,
\begin{eqnarray*}
\displaystyle
{\bm \mu}_{-(r)} = \left\{ \mu \left( {\bm \pi}_{-} \right) \right\},
\
\displaystyle
{\bm \mu}_{+(r)}= \left\{ \mu \left( {\bm \pi}_{+} \right) \right\},
\end{eqnarray*}
and
\begin{eqnarray*}
\displaystyle
{\bf B}_r = {\bf B}_r(N,n) = \left ( \begin{array}{ll}
\displaystyle
{\bf B}_{--r} & {\bf 0}
\\
\displaystyle
{\bf B}_{+-r} & {\bf B}_{++r}
\end{array} \right )
\end{eqnarray*}
since ${\bf B}_{-+r} = {\bf 0}$.

Since also ${\bf C}_{-+r} = {\bf D}_{-+r} = {\bf 0}$,
the inversion principle (\ref{1.8}) implies
\begin{eqnarray}
\displaystyle
{\bf B}_{--r}(N,n)^{-1} = {\bf B}_{--r}(n, N),
\
{\bf C}_{--r}(N,n)^{-1} = {\bf C}_{--r}(n, N),
\
{\bf D}_{--r}(N,n)^{-1} = {\bf D}_{--r}(n, N),
\label{1.13}
\end{eqnarray}
so that
\begin{eqnarray*}
\displaystyle
{\bf B}_{--r}(n,N) \widehat{\bf m}_{-(r)} \mbox{ is an UE of } {\bf m}_{-(r)},
\end{eqnarray*}
\begin{eqnarray*}
\displaystyle
{\bf C}_{--r}(n,N) \widehat{\bm \mu}_{-(r)} \mbox{ is an UE of } {\bm \mu}_{-(r)},
\end{eqnarray*}
and
\begin{eqnarray*}
\displaystyle
{\bf  D}_{--r}(n,N) \widehat{\bf K}_{-(r)} \mbox{ is an UE of } {\bf K}_{-(r)}.
\end{eqnarray*}
The number of parts in ${\bm \pi}$ is denoted by $q({\bm \pi})$:
\begin{eqnarray*}
\displaystyle
q \left( 1^{r_1}2^{r_2} \cdots \right) = r_1+r_2+ \cdots
\end{eqnarray*}
Set  $(n)_i=n(n-1) \cdots(n-i+1)=n!/(n-i)!$.

\section{Products of non central moments}

Here, we derive the result
\begin{eqnarray}
\displaystyle
E \widehat{\bf m} (r) = {\bf B}_{r} {\bf m} (r),
\label{2.1}
\end{eqnarray}
where ${\bf m} (r) =  \left\{ m ({\bm \pi}) \right\}$ and
$m ({\bm \pi})= m_{\pi_{1}} m_{\pi_{2}} \ldots$ for ${\bm \pi} = \left( \pi_1, \pi_2, \ldots \right)$.
Skellam (1949) showed that
\begin{eqnarray*}
\displaystyle
E S_{a_{1}} \cdots S_{a_{s}} = \sum_{\bm \pi}^s \lambda ({\bm \pi}) \sum^{P ({\bm \pi})}
s_{R_{\pi_{1}}} \ s_{R_{\pi_{2}}} \cdots
\end{eqnarray*}
for
\begin{eqnarray*}
\displaystyle
S_{r} = \sum_{i=1}^{n} X_{i}^{r} = n \widehat{m}_{r},
\
\displaystyle
s_{r} = \sum_{i=1}^{N} x_{i}^{r} = N m_{r},
\nonumber
\end{eqnarray*}
where $\displaystyle \sum_{\bm \pi}^s$ sums over partitions ${\bm \pi}$ of $s$,
the partition function $P ({\bm \pi})$ is given by (\ref{1.11}),
$R_{\pi_{1}} = a_{1} + \cdots + a_{\pi_{1}}$,
$R_{\pi_{2}} = a_{\pi_{1}+1} + \cdots + a_{\pi_{1} + \pi_{2}}$, etc
and $\displaystyle \sum^{P ({\bm \pi})} $ sums over all $P ({\bm \pi})$
such terms and $\lambda ({\bm \pi}) = \lambda_{{\bm \pi}}$ is {\it the Carver function}.
For example,
\begin{eqnarray*}
\displaystyle
E \ S_{a_1}S_{a_2}S_{a_3} = \lambda (3)s_{a_1+a_2+a_3}+ \lambda (12) \sum^3 s_{a_1+a_2}s_{a_3}+
\lambda \left(1^3\right)s_{a_1}s_{a_2}s_{a_3}
\end{eqnarray*}
and
\begin{eqnarray}
\displaystyle
E S_{a_{1}} \cdots S_{a_{5}}
&=&
\displaystyle
\lambda (5) s_{a_{1} + \ldots + a_{5}} + \lambda(14) \sum^{3}
s_{a_1 + \ldots + a_{4}} s_{a_{5}}
\nonumber
\\
&&
\displaystyle
+ \lambda (23) \sum^{10} s_{a_{1} + a_{2} + a_{3}}s_{a_4+a_5}
+ \lambda \left(1^{2}3\right) \sum^{10} s_{a_{1} + a_{2} + a_{3}}  s_{a_{4}} s_{a_{5}}
\nonumber
\\
&&
\displaystyle
+ \lambda \left( 12^{2} \right)
\sum^{15} s_{a_{1}+a_{2}} s_{a_{3} + a_{4}}
s_{a_{5}} + \lambda \left( 1^{3}2 \right)
\sum^{10} s_{a_{1}+a_{2}} s_{a_{3}}  s_{a_{4}} s_{a_{5}}
\nonumber
\\
&&
\displaystyle
+ \lambda \left( 1^{5} \right) s_{a_{1}} \cdots s_{a_{5}},
\label{2.4}
\end{eqnarray}
where
\begin{eqnarray*}
&&
\displaystyle
\lambda (3) = e_1 - 3 e_2 + 2 e_3,
\
\displaystyle
\lambda (12)= e_2 - e_3,
\
\displaystyle
\lambda \left(1^3\right) = e_3,
\nonumber
\\
&&
\displaystyle
\lambda (5) = e_{1} - 15 e_{2} + 5 0 e_{3} - 6 0 e_{4} + 24 e_{5},
\nonumber
\\
&&
\displaystyle
\lambda (14) = e_{2} - 7 e_{3} + 12 e_{4} - 6 e_{5},
\nonumber
\\
&&
\displaystyle
\lambda (23) = e_{2} - 4 e_{3} + 5 e_{4} - 2 e_{5},
\nonumber
\\
&&
\displaystyle
\lambda \left( 1^{2}3 \right) = e_{3} - 3 e_{4} + 2 e_{5},
\nonumber
\\
&&
\displaystyle
\lambda \left( 12^{2}  \right) = e_{3} - 2 e_{4} + e_{5},
\nonumber
\\
&&
\displaystyle
\lambda \left( 1^{3}2 \right) = e_{4} - e_{5},
\nonumber
\\
&&
\displaystyle
\lambda \left( 1^{5} \right) = e_{5},
\nonumber
\end{eqnarray*}
where $e_j = (n)_j / (N)_j$.
He wrote (\ref{2.1}) out in full for $s \leq 4$.
We extend this to $s \leq 6$ in Appendix E.
Set
\begin{eqnarray*}
\displaystyle
{\bf S}_{(r)} = \left\{ S_{\pi_{1}} S_{\pi_{2}} \cdots : \left( \pi_{1}, \pi_2, \ldots \right) \mbox{ a partition of } r \right\},
\end{eqnarray*}
and similarly for ${\bf s}_{(r)}$.
It follows from (\ref{2.1}) that
\begin{eqnarray}
\displaystyle
E {\bf S}_{(r)} = {\bf A}_{r} {\bf s}_{(r)},
\label{2.5}
\end{eqnarray}
where
\begin{eqnarray*}
\displaystyle
{\bf A}_r = \left\{ A_{{\bm \pi}, {\bm \pi}^{'}} : {\bm \pi}, {\bm \pi}^{'} \mbox{ partitions of } r \right\}
\end{eqnarray*}
and ${\bf A}_{1}, \ldots, {\bf A}_{6}$ are as follows:
\begin{eqnarray*}
&&
\displaystyle
{\bf A}_1 = \lambda_1,
\
\displaystyle
{\bf A}_2 = \left( \begin{array}{cc}
\displaystyle
\lambda_1 &
\\
\displaystyle
\lambda_2 & \lambda_{11}
\end{array} \right),
\
\displaystyle
{\bf A}_3 = \left( \begin{array}{ccc}
\displaystyle
\lambda_1 & &
\\
\displaystyle
\lambda_2 & \lambda_{11}
\\
\displaystyle
\lambda_3 & 3 \lambda_{12} & \lambda_{1^3}
\end{array} \right),
\\
&&
\displaystyle
{\bf A}_4 = \left[ \begin{array}{lllll}
\displaystyle
\lambda_1 & & & &
\\
\displaystyle
\lambda_2 & \lambda_{11} & & &
\\
\displaystyle
\lambda_2 & 0 & \lambda_{11} & &
\\
\displaystyle
\lambda_3 & 2\lambda_{12} & \lambda_{12} & \lambda_{1^3} &
\\
\displaystyle
\lambda_4 & 4\lambda_{13} & 3 \lambda_{22} & 6 \lambda_{1^22} & \lambda_{1^4}
\end{array} \right],
\\
&&
\displaystyle
{\bf A}_5 = \left[ \begin{array}{lllllll}
\displaystyle
\lambda_1 & & & & & &
\\
\displaystyle
\lambda_2 & \lambda_{11} & & & & &
\\
\displaystyle
\lambda_2 & 0 & \lambda_{11} & & & &
\\
\displaystyle
\lambda_3 & 2 \lambda_{12} & \lambda_{12} & \lambda_{1^3} & & &
\\
\displaystyle
\lambda_3 & \lambda_{12} & 2 \lambda_{12} & 0 & \lambda_{1^3} & &
\\
\displaystyle
\lambda_4 & 3 \lambda_{13} & \lambda_{13} + 3 \lambda_{22} & 3 \lambda_{1^22} & 3\lambda_{1^22} & \lambda_{1^4} &
\\
\displaystyle
\lambda_5 & 5 \lambda_{14} & 10 \lambda_{23} & 10 \lambda_{113} & 15 \lambda_{22} & 10 \lambda_{12} & 3 \lambda_{1^5}
\end{array} \right],
\\
&&
\displaystyle
{\bf A}_6 = \tiny
\left[ \begin{array}{lllllllllll}
\displaystyle
\lambda_1 &&&&&&&&&&
\\
\displaystyle
\lambda_2 & \lambda_{11} &&&&&&&&&
\\
\displaystyle
\lambda_2 & 0 & \lambda_{11} &&&&&&&&
\\
\displaystyle
\lambda_2 & 0 & 0 & \lambda_{11} &&&&&&&
\\
\displaystyle
\lambda_3 & 2 \lambda_{12} & \lambda_{12} & 0 & \lambda_{11} &&&&&&
\\
\displaystyle
\lambda_3 & \lambda_{12} & \lambda_{12} & \lambda_{12} & 0 & \lambda_{1^3} &&&&&
\\
\displaystyle
\lambda_3 & 0 & 3 \lambda_{12} & 0 & 0 & 0 & \lambda_{1^3} &&&&
\\
\displaystyle
\lambda_4 & 3 \lambda_{13} & \lambda_{13} + 3 \lambda_{22} & 0 & 3 \lambda_{1^22} & 3
\lambda_{1^22} & 0 & \lambda_{1^4} &&&
\\
\displaystyle
\lambda_4 & 2 \lambda_{13} & 2 \lambda_{13} + \lambda_{22} & 2 \lambda_{22} & \lambda_{1^22} & 4 \lambda_{1^22} & \lambda_{1^22} & 0 & \lambda_{1^4} & &
\\
\displaystyle
\lambda_5 & 4 \lambda_{14} & \lambda_{14}+6\lambda_{23} & 4 \lambda_{23} & 6
\lambda_{113} & 4 \lambda_{113} +12 \lambda_{122} & 3 \lambda_{122} & 4 \lambda_{1^32} & 6 \lambda_{1^32} & \lambda_{1^5}  &
\\
\displaystyle
\lambda_6 & 6 \lambda_{15} & 15 \lambda_{24} & 10 \lambda_{33} & 15 \lambda_{1^24} &
60 \lambda_{123} & 15 \lambda_{2^3} & 10 \lambda_{1^33} & 45 \lambda_{1^22^2} & 15 \lambda_{1^42} & \lambda_{1^6}
\end{array}
\right].
\end{eqnarray*}
Since $ S_r = n \widehat{m}_r $ and $ s_r = N m_r$, (\ref{2.5}) implies (\ref{2.1}) with
\begin{eqnarray*}
\displaystyle
{\bf B}_r = {\bf D}_{r, n}^{-1} {\bf A}_r {\bf D}_{r, N},
\end{eqnarray*}
where ${\bf D}_{r, N} = {\rm diag} \left\{ N^{q ({\bm \pi})} : r ({\bm \pi}) = r \right\}$.
For example, ${\bf D}_{4, N} = {\rm diag} \left( N, N^2, N^2, N^3, N^4 \right)$.
Writing
\begin{eqnarray*}
&&
\displaystyle
{\bf A}_r (N, n) = {\bf A}_r = \left\{ A_{{\bm \pi}, {\bm \pi}^{'}} (N, n) : {\bm \pi}, {\bm \pi}^{'}
\mbox{ partition of } r \right\},
\nonumber
\\
&&
\displaystyle
{\bf B}_r (N, n) = {\bf B}_r = \left\{ B_{{\bm \pi}, {\bm \pi}'} (N, n) : {\bm \pi}, {\bm \pi}^{'}
\mbox{ partitions of } r \right\},
\nonumber
\end{eqnarray*}
it follows from (\ref{2.1}) that
\begin{eqnarray*}
\displaystyle
E \widehat{m}_{\pi_{1}} \widehat{m}_{\pi_{2}} \cdots = \sum_{{\bm \pi}^{'}}^{r}
B_{{\bm \pi}, {\bm \pi}^{'}} (N, n) m_{\pi_{1}^{'}} m_{\pi_{2}^{'}} \cdots =
b_{\bm \pi} (N, n, {\bf m})
\end{eqnarray*}
say.
By the inversion principle
\begin{eqnarray*}
\displaystyle
{\bf B}_{r} (N, n)^{-1} = {\bf B}_{r} (n, N),
\end{eqnarray*}
so $b_{\bm \pi} \left( n, N, \widehat{\bm m} \right)$ is an UE of $m_{\pi_{1}} m_{\pi_{2}} \cdots$.
For example,
\begin{eqnarray*}
\displaystyle
E \widehat{m}_{1}^{3}  = n^{-3} \left( \lambda_{3} N m_{3} + 3 \lambda_{21} N^{2}
m_{2} m_{1} + \lambda_{1^{3}} N^{3} m_1^{3} \right),
\nonumber
\end{eqnarray*}
so $N^{-3} \left( \overline{\lambda}_{3} n \widehat{m}_{3} + 3 \overline{\lambda}_{21} n^{2}
\widehat{m}_{2} \widehat{m}_{1} + \overline{\lambda}_{1^{3}}
n^{3} \widehat{m}_{1}^{3} \right)$ is an UE of $m_{1}^{3}$,
where $\overline{\lambda}_{\bm \pi} $ is $\lambda_{\bm \pi}$ with $n$ and $N$ reversed.

\begin{example}
Suppose that the population consists of $Np$ ones and $ N (1-p)$ zeros.
Then $m_{r} = p$ for $ r > 0$, so
\begin{eqnarray*}
\displaystyle
E \left( n \widehat{p} \right)^{r} = \sum_{\bm \pi}^{r} \lambda_{\bm \pi} P ({\bm \pi}) (N p)^{q ({\bm \pi})}.
\end{eqnarray*}
So,
\begin{eqnarray*}
\displaystyle
E \widehat{p}^{r} = n^{-r} \sum_{i=1}^{r} a_{r, i} (N,n) (Np)^{i} = b_{r} (N, n, p)
\end{eqnarray*}
say, where
\begin{eqnarray*}
\displaystyle
a_{r, i} (N,n) = \sum_{\bm \pi}^{r} \left\{ \lambda_{\bm \pi} P ({\bm \pi}) : q ({\bm \pi}) = i \right\},
\end{eqnarray*}
and an UE of $p^{r}$ is $b_{r} \left( n, N, \widehat{p} \right)$.
This implies another inversion principle:
${\bf a}_{r} = {\bf a}_{r} (N, n) = \left\{ a_{i, j} (N,n) \right\}$, $r \times r$ with $a_{i, j} = 0$
for $i < j$ has inverse ${\bf a}_{r} (n, N)$.
The first six are given by
\begin{eqnarray*}
\displaystyle
{\bf a}_r = \left( \begin{array}{ll}
\displaystyle
{\bf a}_{r-1} & 0
\\
\displaystyle
{\bf c}_{r}^{'} & \lambda_{1^r}
\end{array} \right),
\end{eqnarray*}
where ${\bf c}_r^{'} = \left( a_{r, 1}, \ldots, a_{r, r-1} \right)$ are
\begin{eqnarray*}
&&
\displaystyle
c_{2} = \lambda_{2},
\
\displaystyle
{\bf c}_{3}^{'} = \left( \lambda_{3}, 3 \lambda_{12} \right),
\
\displaystyle
{\bf c}_{4}^{'} = \left( \lambda_{4}, 4 \lambda_{13} + 3 \lambda_{22}, 6 \lambda_{1^22} \right),
\nonumber
\\
&&
\displaystyle
{\bf c}_{5}^{'} = \left( \lambda_{5}, 5 \lambda_{14} + 10 \lambda_{23},
10 \lambda_{113} + 15 \lambda_{122}, 10 \lambda_{1^{3}2} \right),
\nonumber
\\
&&
\displaystyle
{\bf c}_{6}^{'} = \left( \lambda_{6}, 6 \lambda_{15} + 15 \lambda_{24} + 10 \lambda_{33},
15 \lambda_{1^24} + 60 \lambda_{123} + 15 \lambda_{2^3},
20 \lambda_{1^33} + 45 \lambda_{1^22^2}, 15 \lambda_{1^4}2 \right).
\nonumber
\end{eqnarray*}
\end{example}

The block matrix ${\bf A} = \left(  \begin{array}{cc}
{\bf A}_{1, 1} & {\bf 0}
\\
{\bf A}_{2, 1} & {\bf A}_{2, 2} \end{array} \right) $
has ${\bf A}^{-1} = \left( \begin{array}{ll}
{\bf A}_{1, 1}^{-1} & {\bf 0}
\\
{\bf A}^{2, 1} & {\bf A}_{2, 2}^{-1} \end{array} \right)$,
where ${\bf A}^{2, 1} = -{\bf A}_{2, 2}^{-1} {\bf A}_{2, 1} {\bf A}_{1, 1}^{-1}$.
So, if ${\bf A} = {\bf A} (N, n)$ then
${\bf A} (N, n)^{-1}  = {\bf A} (n, N)$ if and only if
${\bf A}_{i, i} (N, n)^{-1} = {\bf A}_{i, i} (n, N)$, $i = 1, 2$ and
$-{\bf A}_{2, 2} (n, N) {\bf A}_{2, 1} (N, n) {\bf A}_{1, 1} (n, N) = {\bf A}_{2, 1} (n,N)$.

\begin{example}
Suppose $F$ is Poisson $(\lambda)$, that is for $i \geq 0$
a proportion $ p_{i} = e^{-\lambda} \lambda^{i} / i !$ of the population equal $i$.
This implies that $ N = \infty$, but the algebra can be carried out as
if $ N < \infty$, then $N=\infty$ substituted at the end.
For $ i \geq 1$ the Stirling numbers of the second kind $S_{i, 1}, \ldots, S_{i, i}$ are defined by
\begin{eqnarray*}
\displaystyle
j^{i} = \sum_{k=1}^{i} S_{i, k} (j)_{k}.
\end{eqnarray*}
Then $m_{i} = E X^{i} = \displaystyle \sum_{k=1}^{i} S_{i, k} \lambda^{k}$.
So, $\left( m_{1}, \ldots, m_{r} \right) = \left( \lambda, \ldots, \lambda^{r} \right) {\bm \alpha}_{r}^{'}$,
where ${\bm \alpha}_{r} = \left( S_{i, j} \right)$ is $r \times r$ and $S_{i, j} = 0$ for $j > i$.
So, ${\bm \alpha}_{r}^{'-1} \left( \widehat{m}_{1}, \ldots, \widehat{m}_{r} \right)$ is an UE of
$\left( \lambda, \ldots, \lambda^{r} \right)$.

But other UEs of $\lambda^{r}$ can be obtained from products of $\left\{ \widehat{m}_{i} \right\}$.
For example,
\begin{eqnarray*}
\displaystyle
E \left( n \widehat{m}_{1} \right)^{r} = \sum_{i=1}^{r} d_{r, i} N^{i},
\end{eqnarray*}
where
\begin{eqnarray*}
\displaystyle
d_{r, i} = \sum_{\bm \pi}^{r} \left\{ \lambda_{\bm \pi} P ({\bm \pi}) m_{\pi_{1}} m_{\pi_{2}} \cdots : q ({\bm \pi}) = i \right\},
\end{eqnarray*}
Thus,
\begin{eqnarray*}
\displaystyle
E \left( n \widehat{m}_{1} \right)^2
&=&
\displaystyle
\left( \lambda^{2} + \lambda \right) \lambda_{2}N + \lambda^{2} \lambda_{11} N^{2}
\nonumber
\\
&=&
\displaystyle
\left( \lambda^{2} + \lambda \right) n + \lambda^{2} (n)_{2} \mbox{ at } N = \infty .
\nonumber
\end{eqnarray*}
So, at $N= \infty$, an UE of $\theta = \lambda^{2}$ is $ \widehat{\theta}_{1} = \widehat{m}_{1}^{2} - n^{-1} \widehat{m}_{1}$.

But by above $ \widehat{\theta}_{2} = \widehat{m}_{2} - \widehat{m}_{1}$ is also an UE.

One can now apply (\ref{2.4}) to calculate var $\widehat{\theta}_{1}$,
var $\widehat{\theta}_{2} $ to see when $\widehat{\theta}_{1}$ is more efficient than $\widehat{\theta}_{2}$.
\end{example}

\section{Products of powers of the mean and central moments and the Inversion Principle}

Here, we derive $E \widehat{\mu} ({\bm \pi})$ and the UE of $\mu({\bm \pi}) = \mu_{\pi_1} \mu_{\pi_2} \cdots $,
where $\mu_1=\mu$.
We also elaborate on the inversion principle (\ref{1.9}).

By Sukhatme (1944),
\begin{eqnarray}
\displaystyle
E\widehat{\mu} \left( {\bm \pi}_{-} \right) = \sum_{{\bm \pi}^{'}_{-}}^r C_{{\bm \pi}_{-}, {\bm \pi}'_{-}}
\mu \left( {\bm \pi}'_{-} \right),
\label{3.1}
\end{eqnarray}
where $C_{{\bm \pi}_{-}, {\bm \pi}'_{-}}$ is given for $ r \leq 6 $ in terms of $e_j=(n)_j/(N)_j$.

Set
\begin{eqnarray*}
&&
\displaystyle
{\bm \mu}_{-(r)} = \left\{ \mu \left( {\bm \pi}_{-} \right) : r \left( {\bm \pi}_{-} \right) =r \right\},
\
\displaystyle
{\bm \mu}_{+(r)} = \left\{ \mu \left( {\bm \pi}_{+} \right) : r \left( {\bm \pi}_{+} \right) =r \right\},
\nonumber
\\
&&
\displaystyle
{\bm \mu}_{(r)}' = \left( {\bm \mu}_{-(r)}', {\bm \mu}_{+(r)}' \right) = \left\{ \mu ({\bm \pi}) : r({\bm \pi})=r \right\}',
\nonumber
\end{eqnarray*}
where $r({\bm \pi}) = \pi_1+\pi_2+ \cdots$, the order of ${\bm \pi}$.
In particular,
\begin{eqnarray*}
&&
\displaystyle
\mu_{-(2)} = \mu_2,
\
\displaystyle
\mu_{-(3)}=\mu_3,
\
\displaystyle
{\bm \mu}_{-(4)} = \left(\mu_4,\mu_2^2\right)',
\
\displaystyle
{\bm \mu}_{-(5)} = \left( \mu_5, \mu_2\mu_3 \right)',
\nonumber
\\
&&
\displaystyle
{\bm \mu}_{-(6)} = \left( \mu_6, \mu_2\mu_4, \mu_3^2, \mu_2^3 \right)',
\nonumber
\\
&&
\displaystyle
\mu_{+(2)} = \mu^2,
\
\displaystyle
{\bm \mu}_{+(3)}' = \left(\mu\mu_2, \mu^3\right),
\nonumber
\\
&&
\displaystyle
{\bm \mu}_{+(4)}' = \left( \mu\mu_3, \mu^2\mu_2, \mu^4 \right),
\nonumber
\\
&&
\displaystyle
{\bm \mu}_{+(5)}' = \left( \mu\mu_4, \mu\mu_2^2, \mu^2\mu_3, \mu^3\mu_2, \mu^5 \right),
\nonumber
\\
&&
\displaystyle
{\bm \mu}_{+(6)}' = \left( \mu\mu_5, \mu\mu_2\mu_3, \mu^2\mu_4, \mu^2\mu_2^2, \mu^3\mu_3, \mu^4\mu_2, \mu^6 \right).
\nonumber
\end{eqnarray*}
Let $\widehat{\bm \mu}_{-(r)}, \widehat{\bm \mu}_{+(r)}, \widehat{\bm \mu}_(r)$ denote their sample
versions, that is with $F$ replaced by $\widehat{F}$.

Then the result (\ref{3.1}) can be written as
\begin{eqnarray*}
\displaystyle
E\widehat{\bm \mu}_{-(r)} = {\bf C}_{--r} {\bm \mu}_{-(r)},
\end{eqnarray*}
where ${\bf C}_{--r} = \left( C_{{\bm \pi}_{-}, {\bm \pi}'_{-}}: {\bm \pi}_{-}, {\bm \pi}'_{-}
\mbox{ partitions of $r$ excluding  ones } \right)$.
So, $\widetilde{\bm \mu}_{-(r)} = {\bf C}_{--r}^{-1} \widehat{\bm \mu}_{-(r)}$ is an UE of ${\bm \mu}_{-(r)}$.
The coefficients of $C_{{\bm \pi}_{-}, {\bm \pi}'_{-}}$ are given in Appendix A for $r \leq 6$.

Sukhatme (1944) gives $\left\{ C_{{\bm \pi}, {\bm \pi}'} \right\}$ such that
\begin{eqnarray*}
\displaystyle
E \left( \widehat{\mu} - \mu \right)^{r_1} \widehat{\mu}_2^{r_2} \widehat{\mu}_3^{r_2} \cdots =
\sum_{{\bm \pi}_{-}^{'}}^r C_{{\bm \pi}, {\bm \pi}^{'}_{-}} \mu \left( {\bm \pi}^{'}_{-} \right)
\end{eqnarray*}
for ${\bm \pi} = \left( 1^{r_1} 2^{r_2} \cdots \right)$.
To obtain $ E \widehat{\mu}^{r_1} \widehat{\mu}_2^{r_2} \widehat{\mu}_3^{r_3} \cdots$,
one just expands $\widehat{\mu}^{r_1}=(\delta+\mu)^{r_1}$, where $\delta=\widehat{\mu}-\mu$.
So, one obtains
\begin{eqnarray*}
\displaystyle
E\widehat{\mu} ({\bm \pi}) = \sum_{{\bm \pi}'}^r C_{{\bm \pi}, {\bm \pi}'} \mu \left( {\bm \pi}' \right),
\end{eqnarray*}
that is
\begin{eqnarray*}
\displaystyle
E \widehat{\bm \mu}_{(r)} = {\bf C}_r {\bm \mu}_{(r)},
\end{eqnarray*}
where ${\bf C}_r = \left( C_{{\bm \pi}, {\bm \pi}'} : {\bm \pi}, {\bm \pi}' \mbox{ partitions of } r \right)$ and
\begin{eqnarray}
\displaystyle
C_{1^i {\bm \pi}_{-}, 1^j {\bm \pi}'_{-}} = {i \choose j} C_{1^{i-j} {\bm \pi}_{-}, {\bm \pi}'_{-}}
\label{3.7}
\end{eqnarray}
for ${\bm \pi}_{-}, {\bm \pi}'_{-}$ partitions excluding ones with $i + r ({\bm \pi}) = j + r \left( {\bm \pi}' \right)$,
and the left hand side of (\ref{3.7}) is zero if $j>i$.
For completeness $\left\{ C_{{\bm \pi}, {\bm \pi}'} \right\}$ are given in Appendix A for $r \leq 6$.
An UE of ${\bm \mu}_{(r)}$ is $\widetilde{\bm \mu}_{(r)} = {\bf C}_r^{-1} \widehat{\bm \mu}_{(r)}$.
So, by (\ref{3.9}) below
\begin{eqnarray*}
\displaystyle
\widetilde{\bm \mu}_{+(r)} = C_r^{2, 1} \widehat{\bm \mu}_{-(r)} +
C_r^{2, 2} \widehat{\bm \mu}_{+(r)}
\end{eqnarray*}
is an UE of ${\bm \mu}_{+(r)}$.

Note that for ${\bm \pi}_1, {\bm \pi}_2, \ldots$ partitions of $r_1, r_2, \ldots$
\begin{eqnarray*}
\displaystyle
E \widehat{\mu} \left( {\bm \pi}_1\right)
E \widehat{\mu} \left( {\bm \pi}_2 \right) \cdots =
\sum_{{\bm \pi}'_1}^{r_1}\sum_{{\bm \pi}'_2}^{r_2} \cdots
C_{{\bm \pi}_1, {\bm \pi}'_1} C_{{\bm \pi}_2, {\bm \pi}'_2}
\cdots
\mu \left( {\bm \pi}'_1 + {\bm \pi}'_2+ \cdots \right),
\end{eqnarray*}
where
\begin{eqnarray*}
\displaystyle
\mu \left( 1^{r_{1, 1}} 2^{r_{1, 2}} +\cdots+1^{r_{2, 1}}2^{r_{2, 2}} +\cdots+ \cdots \right) = \mu \left( 1^{r_1}2^{r_2} \cdots \right)
\end{eqnarray*}
for $r_1=r_{1, 1}+r_{2, 1}+ \cdots$, $r_2=r_{1, 2}+r_{2, 2}+ \cdots$, etc.
So, $\displaystyle \prod_i E\widehat{\mu} \left( {\bm \pi}_i \right)$
has the form $\displaystyle \sum_{\bm \pi}^r D_{\bm \pi} \mu ({\bm \pi})$,
where $r=r_1+r_2+ \cdots$ and has UE $\displaystyle \sum_{{\bm \pi}'}^r D^{*}_{{\bm \pi}'} \widehat{\mu} ({\bm \pi})$,
where $D_{{\bm \pi}'}^{*} = \displaystyle \sum_{\bm \pi}^r D_{\bm \pi} C^{{\bm \pi}, {\bm \pi}'}$.
For example, in this way we can write down UEs for
$\mu^{r_1} \left( E\widehat{\mu}_2 \right)^{r_2} \left( E\widehat{\mu}_3 \right)^{r_3} \cdots$
for $r_1+2r_2+3r_3+ \cdots \leq 6$.

Since ${\bf C}_r$ has the form
$\left( \begin{array}{cc}
{\bf C}_{1, 1} & {\bf 0}
\\
{\bf C}_{2, 1} & {\bf C}_{2, 2}
\end{array}\right)$
with ${\bf C}_{1, 1} = {\bf C}_{--r}$, its inverse is
\begin{eqnarray}
\displaystyle
{\bf C}_r^{-1} =
\left( \begin{array}{cc}
\displaystyle
{\bf C}_r^{1, 1} & {\bf 0}
\\
\displaystyle
{\bf C}_r^{2, 1} & {\bf C}_r^{2, 2}
\end{array} \right)
\label{3.9}
\end{eqnarray}
with ${\bf C}_r^{1, 1} = {\bf C}_{1, 1}^{-1}$,
${\bf C}_r^{2, 1} = -{\bf C}_{2, 2}^{-1} {\bf C}_{2, 1} {\bf C}_{1, 1}^{-1}$
and ${\bf C}_r^{2, 2} = {\bf C}_{2, 2}^{-1}$.

However, this dimension-reduction method to obtain ${\bf C}_r^{-1}$
is not necessary due to the inversion principle,
${\bf C}_r(N, n)^{-1} = {\bf C}_r(n, N)$.
This implies ${\bf C}_{--r}(N, n)^{-1} = {\bf C}_{--r}(n, N)$ as noted in (\ref{1.13}).
This result was discovered using MAPLE and has
been verified for the ${\bf C}_{--r}$, ${\bf C}_r$ of Appendix A,
that is for $ \left\{ {\bf C}_{--r}: r\leq 6 \right\}$, $\left\{ {\bf C}_r : r\leq 5 \right \}$.
Its proof is given in Section  7.
Setting $\left( C^{{\bm \pi}, {\bm \pi}'} \right) = {\bf C}_r^{-1}$
and ${\bf C}^{{\bm \pi}, {\bm \pi}'}(N, n) = {\bf C}^{{\bm \pi}, {\bm \pi}'}$,
we can write these as $C^{{\bm \pi}, {\bm \pi}'} (N, n) = C_{{\bm \pi}, {\bm \pi}'}(n, N)$
and $C^{{\bm \pi}_{-}, {\bm \pi}'_{-}}(N, n) = C_{{\bm \pi}_{-}, {\bm \pi}'_{-}} (n, N)$.

We now show how to obtain UEs for products of cumulants.
The $r$th cumulant $\kappa_r$ may be written as
\begin{eqnarray*}
\displaystyle
\kappa_r = {\bf a}_r' {\bm \mu}_{-(r)}
\end{eqnarray*}
for $r \geq 2$, where $a_2 = a_3 = 1$, ${\bf a}_4 = (1,-3)$, ${\bf a}_5 = (1, -10)$, ${\bf a}_6 = (1,-15,-10,30)$.
So, $\widetilde{\kappa}_r = {\bf a}_r' \widetilde{\bm \mu}_{-(r)}$ is an UE for $\kappa_r$.
For the case $N=\infty$, Fisher (1929) gave these {\it $k$-statistics}, ${\bm \mu}_{-(r)}$
for $r\leq6$ and their joint cumulants.
Sukhatme obtained $C_{{\bm \pi}, {\bm \pi}'}$ using Fisher's method.

Wishart (1952) gave tables for the poly-kays, $k_{r_1, r_2, \ldots}$, the UE of
$\kappa_{r_1}\kappa_{r_2} \cdots$, in terms of the symmetric functions.
This is reproduced in Appendix Table 11 of Stuart and Ord (1987) who use
$l_{r_1, r_2, \ldots}$ for $\kappa_{r_1, r_2, \ldots}$.
The symmetric functions are given by
their Appendix Table 10 in terms of the power sums $s_r = \displaystyle \sum_{i}x_i^r.$
However, no explicit formulae for the polykays in terms of the sample moments or
cumulants appear to be available.
We rectify this in Appendix D, not by using Tables 11 and 10, but by writing
\begin{eqnarray*}
\displaystyle
\kappa_{\pi_1} \kappa_{\pi_2} \cdots =
\left\{
\begin{array}{cl}
\displaystyle
{\bf a} ({\bm \pi})' {\bm \mu}_{-(r)}, & \mbox{if } {\bm \pi} = \left( \pi_1, \pi_2, \ldots\right) \mbox{ does not contain $1$,}
\\
\displaystyle
{\bf b} ({\bm \pi})' {\bm \mu}_{(r)}, & \mbox{if } {\bm \pi} \mbox{ contains $1$,}
\end{array}
\right.
\end{eqnarray*}
where $r=\pi_1+\pi_2+ \cdots$.
So the UE of $\kappa_{\pi_1}\kappa_{\pi_2} \cdots$ is
${\bf a} \left( {\bm \pi} \right)' \widetilde{\bm \mu}_{-(r)}$ or
${\bf b} ({\bm \pi})' \widetilde{\bm \mu}_{-(r)}$.

In Section 12.22, Stuart and Ord (1987) give a number of references on this subject,
and say that Dwyer and Tracy (1980) give UEs for products of central moments.
This is not so: they give the multivariate equivalent of the UE
for $E\left( \overline{X}-\mu \right)^r$ for $r\leq5$.
They also give $E\widehat{\mu}_r$ and $E\left( \overline{X}-\mu \right)^r$ for $r\leq6$ by
a different method, or rather the equivalent multivariate results.

An alternative method of getting UEs for $E\left( \overline{X}-\mu \right)^r$, and indeed
for $E \left( \overline{X}-\mu \right)^i \widehat{\mu} \left( {\bm \pi}_{-} \right)$ is to note that
\begin{eqnarray*}
\displaystyle
E \left( \overline{X} - \mu \right)^i \widehat{\mu} \left( {\bm \pi}_{-} \right) =
\sum_{{\bm \pi}'_{-}}^{r+i}C_{1^i {\bm \pi}_{-}, {\bm \pi}'_{-}} \mu \left( {\bm \pi}'_{-} \right)
\end{eqnarray*}
and so by (\ref{1.8}) has UE
\begin{eqnarray*}
\displaystyle
\sum_{{\bm \pi}''_{-}}^{r+i} D_{1^i {\bm \pi}_{-}, {\bm \pi}''_{-}} \widehat{\mu} \left( {\bm \pi}''_{-} \right),
\end{eqnarray*}
where
\begin{eqnarray*}
\displaystyle
D_{1^i {\bm \pi}_{-}, {\bm \pi}''_{-}} = \sum_{{\bm \pi}'_{-}}^{r+i} C_{1^i {\bm \pi}_{-}, {\bm \pi}'_{-}}
C^{{\bm \pi}'_{-}, {\bm \pi}''_{-}},
\end{eqnarray*}
and $C^{{\bm \pi}'_{-}, {\bm \pi}''_{-}} = C_{{\bm \pi}'_{-}, {\bm \pi}''_{-}} (n, N)$.

Set $\mu_{a}^{*} = n^{-1} \displaystyle \sum_{i=1}^{n} \left( X_{i} - \mu \right)^{a}$.
Given a partition ${\bm \pi}$ of $r$ and a set of numbers
$\left\{ a_1, \ldots, a_r \right\}$, the set can be divided into groups
of sizes $\pi_{1}, \pi_{2}, \ldots$ in $P ({\bm \pi})$ ways.
Let $ R_{j} = \displaystyle \sum{a_{i}}$ summed over the $j$th group.
By (\ref{3.1}),
\begin{eqnarray}
\displaystyle
E \prod_{i=1}^{r} \mu_{a_{i}}^{*} = n^{-r} \sum_{\bm \pi}^{r} \lambda ({\bm \pi})
N^{q ({\bm \pi})} \sum^{P ({\bm \pi})} \mu_{R_{1}} \mu_{R_{2}} \cdots .
\label{2.15}
\end{eqnarray}

Skellam (1949) checked Sukhatme's first seventeen formulas using (\ref{3.1}) and also checked all his $\lambda ({\bm \pi})$.
We now use (\ref{3.1}) to obtain the formula for $E \left( \widehat{\mu} - \mu \right) \widehat{\mu}_{5}$
overlooked by Sukhatme.
Note that
\begin{eqnarray*}
\displaystyle
\left( \widehat{\mu} - \mu \right) \widehat{\mu}_{5} =
\mu_{1}^{*} \left\{ \mu_{5}^{*} - 5 \mu_{1}^{*} \mu_{4}^{*} + 10 \mu_{1}^{*2} m_{3}^{*} -
10 \mu_{1}^{*3} m_{2}^{*} + 4 \mu_{1}^{*5} \right\}.
\end{eqnarray*}
By (\ref{2.15}),
\begin{eqnarray*}
&&
\displaystyle
n^{2} E \mu_{1}^{*} \mu_{5}^{*} = \lambda (2) N \mu_{6},
\nonumber
\\
&&
\displaystyle
n^{3} E \mu_{1}^{*2} \mu_{4}^{*} = \lambda (3) N \mu_{6} + \lambda (21) N^{2} \mu_{2} \mu_{4},
\nonumber
\end{eqnarray*}
and so on.
Collecting terms, $E \left( \widehat{\mu} - \mu \right) \widehat{\mu}_{5} =
\displaystyle \sum_{{\bm \pi}_{-}}^6  \  C_{15, {\bm \pi}_{-}}  \  m \left( {\bm \pi}_{-} \right)$,
where $\left\{ C_{15, {\bm \pi}_{-}} \right\}$ are given in Appendix A.

By (\ref{3.7}), $\left\{ C_{15, {\bm \pi}_{+}} \right\}$
are given by  $C_{15, 1^{j} {\bm \pi}} = C_{5, {\bm \pi}}$ if $j = 1$
and $C_{15, 1^{j} {\bm \pi}} = 0$ if $j > 1$.
Note that $\lambda ({\bm \pi})$ has the form
\begin{eqnarray*}
\displaystyle
\lambda ({\bm \pi}) = \sum_{j=q({\bm \pi})}^{r ({\bm \pi})}
c_je_j
\end{eqnarray*}
with
\begin{eqnarray*}
\displaystyle
c_{q({\bm \pi})} = 1,
\
\displaystyle
c_{r({\bm \pi})} = \prod_{i=1}(-1)^{\pi_i-1}
\left( \pi_i - 1 \right)! = (-1)^{r({\bm \pi}) - q({\bm \pi})} \prod_i \left( \pi_i-1 \right)!.
\end{eqnarray*}

Now suppose $X_1,\ldots,X_n$ is a random sample from an {\it infinite} population.
Then Sukhatme noted that the same results (1.6), (1.17) hold with
$N^ie_j=0$ for $i<j$ and $N^j e_j=(n)_j$.
Thus only the first term in the expansion
\begin{eqnarray*}
\displaystyle
\lambda({\bm \pi}) = \sum_{j=q({\bm \pi})}^{r({\bm \pi})} c_j e_j
\end{eqnarray*}
makes a contribution: $N^{q({\bm \pi})} \lambda({\bm \pi}) = c_j(n)_j$  at $j = q({\bm \pi})$.

Since the infinite population is that of most interest, we give these simplified
formulas for $C_{{\bm \pi}, {\bm \pi}'}$ in Appendix B.
Because of the inversion principle, for obtaining UEs for $N=\infty$ we need
$C_{{\bm \pi}, {\bm \pi}'}$ at $n=\infty$.
These are given in Appendix C.
Putting $n=\infty$ is the same as putting
$n^{-i}e_j=0$ for $i>j$ and $n^{-j}e_j=1/(N)_j$.
Thus only the last term in the expansion
\begin{eqnarray*}
\displaystyle
\lambda({\bm \pi}) = \sum_{j = q({\bm \pi})}^{r({\bm \pi})}c_je_j
\end{eqnarray*}
makes a contribution: $n^{-r({\bm \pi})} \lambda({\bm \pi}) = c_j / (N)_j$ at $j = r({\bm \pi})$.
For example,
\begin{eqnarray*}
\displaystyle
C_{12, 3}
&=&
\displaystyle
N \left\{ \lambda (1)n^{-1}-\lambda (2)n^{-2} \right\}
\nonumber
\\
&=&
\displaystyle
N \left\{ e_1n^{-1} - \left(e_1-e_2\right) n^{-2} \right\}
\nonumber
\\
&=&
\displaystyle
\left\{ \begin{array}{ll}
\displaystyle
(n)_1 \left( n^{-1}-n^{-2} \right), & \mbox{ at } N=\infty,
\nonumber
\\
\displaystyle
N \left\{ 1/(N)_1+1/(N)_2 \right\}, & \mbox{ at } n=\infty.
\end{array}
\right.
\end{eqnarray*}
These results agree with James (1958) who gave UEs of $\mu_r$ for $r\leq6$, $N=\infty$.

\section{Multivariate extensions}

For the multivariate case one considers the population ${{x_1,\ldots,x_N}}$ to
be in $R^r$ not $R$.  For $X =(X_1,\ldots,X_r)'\sim F$,define its cross-moments as \\
\begin{eqnarray*}
m_{r.i_1, \ldots, i_r}=E \ X_{i_1} \cdots X_{i_r}, \ \
\mu_{r.i_1, \ldots, i_r}=E(X_{i_1}-m_{1.i_1}) \cdots E(X_{i_r}-m_{1.i_r})
\end{eqnarray*}
and cross-cumulants as $\kappa_{r.i_1, \ldots, i_r}$. (The subscript $r.$ serves to
avoid confusion with the univariate notation.)\\

Fisher noted that the UE of $\kappa_{r.1, \ldots, r}$ follows from that of $\kappa_r$
for the univariate problem.  In the same way for $t,X$ in $R^r$,
$\mu_{r.1, \ldots, r}$ is just the coefficient of $t_1 \cdots t_r/r!$ in $\mu_r(Y)$ for $Y=t'X$.
But for $a(x):R^r\rightarrow R$ any function, and $Y=a(X)$, the
UE of $\mu_r(Y)$ is as given previously with $X_j$ replaced by $a(X_j)$.\\
For example $\widetilde{\mu}_2/(1-n^{-1})$ implies
$\widetilde{\mu}_{2.11}=\widehat{\mu}_{2.11}/(1-n^{-1})$ is an UE of $\mu_{2.11}$ where
$\widehat{\mu}_{2.11}=\mu_{2.11}(F)$.\\
The same method gives expectations for products of sample cross-moments and UEs
for products of cross-moments for the multivariate versions of \\
$\mu_2,\mu^2,\mu_3,\mu^3,\mu_4,\mu_2^2,\mu^4,\mu_5, \ldots$ but fails for the
multivariate versions of
\\ \begin{eqnarray}\mu_2\mu,\mu_3\mu,\mu_2\mu^2,\mu_3\mu_2, \ldots\end{eqnarray}\\
For example the coefficient of $2t_1t_2$ in
\begin{eqnarray*}
E\widehat{\mu}^2=C_{1^2.2}\mu_2+C_{1^2.1^2}\mu^2 \mbox{ is }
E\widehat{m}_{1.1}\widehat{m}_{1.2}=C_{1^2.2}\mu_{2.12}+C_{1^2.1^2}m_{1.1}m_{1.2}
\end{eqnarray*}
and in
\begin{eqnarray*}
E(C^{1^2.2}\widehat{\mu}_2+C^{1^2.1^2}\widehat{\mu}^2)=\mu^2
\mbox{ is }
E(C^{1^2.2}\widehat{\mu}_{2.12}+C^{1^2.1^2}\widehat{m}_{1.1}\widehat{m}_{1.2})
=m_{1.1}m_{1.2},
\end{eqnarray*}
but the coefficient of $3!t_1t_2t_3$ in
$E\widehat{\mu}_2\widehat{\mu}=C_{21.3}\mu_3+C_{21.21}\mu^2\mu $
is \\ \begin{eqnarray}E(\sum^3\widehat{\mu}_{2.12}\widehat{m}_{1.3})/3=
                   C_{21.3}\mu_{3.123}+C_{21.21}(\sum^3\mu_{2.12}m_{1.3})/3\end{eqnarray}
\\where $\sum^3\mu_{2.12}m_{1.3}=\mu_{2.12}m_{1.3}+\mu_{2.23}m_{1.1}+\mu_{2.13}m_{1.2}.$
\\ \\
However from (2.2) one might surmise that
\begin{eqnarray*}
E\widehat{\mu}_{2.12}\widehat{m}_{1.3}=C_{21.3}\mu_{3.123}+C_{21.21}\mu_{2.12}m_{1.3}
\end{eqnarray*}
which is correct, as can be shown from the multivariate symmetric function
\begin{eqnarray*}
\langle 21 \rangle ^n_{123}=\sum^{'}X_{i1}X_{i2}X_{j3}/(n)_2 \mbox{ has mean }
 \langle 21 \rangle ^N_{123}=\sum^{'}x_{i1}x_{i2}x_{j3}/(N)_2.
\end{eqnarray*}
Similarly, the analogous results for (2.1) and (3.1) also hold.\\
These can be written for fixed $i_{1}, \cdots , i_{r}$ in
\{ $1, \cdots , r$ \} as
\begin{eqnarray}
E \widehat{m}(\pi)_{r.i_1, \ldots, i_r}=\sum_{\pi'}^r \mathbf{B}_{\pi.\pi'}m(\pi')_{r.i_1, \ldots, i_r}
\end{eqnarray}
\begin{eqnarray}
\mbox{ and }
E \widehat{mu}(\pi)_{r.i_1, \ldots, i_r}=\sum_{\pi'}^r \mathbf{C}_{\pi.\pi'}mu(\pi')_{r.i_1, \ldots, i_r}
\end{eqnarray}
\\
where for $\pi_{1} \leq \pi_{2} \leq \ldots , m (\pi)_{r.i_1, \ldots i_r,}
= m_{\pi_{1}. i_{1} \ldots i_{\pi_1}} m_{\pi_{2}. j_{1} \ldots j_{\pi_2}} \ldots $
where
$j_{k} = i_{k + \pi_{1}} , $and $\mu(\pi)_{r.i_1, \ldots, i_r}$ is similarly defined,
\begin{eqnarray*}
{\mathbf B}_{\pi . \pi^{'}} = P (\pi^{'})^{-1} B_{\pi. \pi^{'}}
\sum^{P(\pi)}
\end{eqnarray*}
and
\begin{eqnarray*}
\mathbf{C}_{\pi.\pi^{'}} = P (\pi^{'})^{-1}  \  C_{\pi. \pi^{'}}
\sum^{P(\pi^{'})}
\end{eqnarray*}
are now operators.
For example
\begin{eqnarray}
m (12)_{3.ijk} = m_{1.i} m_{2.jk}
\end{eqnarray}
\begin{eqnarray}
\mbox{ so } E \widehat{m}_{1.i} \widehat{m}_{2.jk} &=& \mathbf{B}_{12.3} m_{3.ijk}
+ \mathbf{B}_{12.12} m_{1.i} m_{2.jk} + \mathbf{B}_{12.1^3} m_{1.i} m_{2.j} m_{3.k}, \nonumber \\
&=& n^{-2} ( \lambda_{2} N m_{3.ijk} + \lambda_{11}  N^{2} m_{1.i} m_{2.jk} )
\end{eqnarray}
Similarly $ m (1^{3})_{3.ijk} = m_{1.i} m_{2.j} m_{3.k} $  so
\begin{eqnarray}
E \widehat{m}_{1.i} \widehat{m}_{2.j} \widehat{m}_{3.k} =
\mathbf{B}_{1^{3}.3} m_{3.ijk} + \mathbf{B}_{1^{3}.12} m_{1.i} m_{2.jk}+
\mathbf{B}_{1^{3}.1^{3}} m_{1.i} m_{1.j} m_{1.k} \nonumber \\
= n^{-3} ( \lambda_{3} N m_{3.ijk} + \lambda_{21} N^{2} \sum^{3} m_{1.i}
m_{2.jk} + \lambda_{1^{3}} N^{3} m_{1.i} m_{1.j} m_{1.k} ).\label{4.7}
\end{eqnarray}
$ \mbox{ Note that (3.7) involves } \sum^{3} m_{1.i} m_{2.jk} \mbox{ but
(3.6) only } m_{1.i} m_{2.jk}.$

Set  $ m_{(i_1 \ldots i_{r})} $
\begin{eqnarray}
= \{ m (\pi)_{r.j_1 \ldots j_r} : \pi \mbox{
partition of } r , (j_{1} \ldots j_{r}) \mbox{ a permutation of } (i_{1} \ldots
i_{r} ) \} \label{4.8}
\end{eqnarray}
For example
\begin{eqnarray}
m_{(ijk)} = ( m_{3.ijk}, m_{1.i} m_{2.jk}, m_{1.j} m_{2.k i}, m_{1.k} m_{2.ij},m_{1.i}  m_{1.j} m_{1.k} )^{'}  \label{4.9}
\end{eqnarray}
So $ m_{(i_{1} \ldots i_{r})} \mbox{ has dimension } N_{r} =
 \sum_{\pi}^{r} \ P (\pi).$
Then (4.3), (4.4) can be written as
\begin{eqnarray}
E \widehat{m}_{(i_{1} \ldots i_{r})} = \mathbf{B}_{r} m_{(i_{1} \ldots i_{r})}, \
E \widehat{mu}_{(i_{1} \ldots i_{r})} = \mathbf{C}_{r} mu_{(i_{1} \ldots i_{r})} ,
\end{eqnarray}
where $ \mathbf{B}_{r} , \mathbf{C}_{r}$ are $N_{r} \times  N_{r}.$

This form is most useful for $(i_{1} \ldots i_{r}) = (1 \cdots r)$ :
other cases where 2 or more of $i_{1} \ldots i_{r} $ are equal
could have their dimension reduced.  But we shall achieve this
shortly anyway from (4.10).

Because $N_{r}$ is also the number of multivariate symmetric
functions, the proof of the univariate invariance principle given in
Section 7 extends immediately to prove {\bf the multivariate invariance
principle} :
\begin{eqnarray}
\mathbf{B}_{r} (N,n) &=& \mathbf{B}_{r} \mbox{ and } \mathbf{C}_{r} (N,n) = \mathbf{C}_{r} \mbox{ satisfy } \label{3.11} \\
\mathbf{B}_{r} (N,n)^{-1} &=& \mathbf{B}_{r} (n,N) \mbox{ and } \mathbf{C}_{r} (N,n)^{-1} = \mathbf{C}_{r} (n, N). \label{3.12}
\end{eqnarray}
\begin{eqnarray*}
\mbox{ So } \mathbf{B}_{r} (n,N) \widehat{m}_{(i_{1} \ldots i_{r})} \mbox{ is an
UE of } m_{(i_{1} \ldots i_r )}, \\
\mbox{ and } \mathbf{C}_r (n,N) \widehat{mu}_({i_{1} \ldots i_{r})} \mbox{ is an
UE of } m u_{(i_{1} \ldots i_{r})}.
\end{eqnarray*}
That is,
\begin{eqnarray}
\sum_{\pi'}^{r} P (\pi')^{-1} B_{\pi.\pi^{'}} \sum^{P (\pi^{'})} \widehat{m}
(\pi^{'})_{r.i_1, \ldots i_{r}} \mbox{ is an UE of } m (\pi)_{r.i_{1} \ldots i_{r}}
\label{3.13}
\end{eqnarray}
and
\begin{eqnarray}
\sum_{\pi^{'}}^{r} P (\pi^{'})^{-1} C_{\pi.\pi^{'}} \sum^{P (\pi^{'})}
\widehat{mu} (\pi^{'})_{r.i_{1} \ldots i_{r}} \mbox{ is an UE of } m u
(\pi)_{r.i_{1} \ldots i_{r}}.
\label{3.14}
\end{eqnarray}
For example letting $\overline{\lambda}_{\pi}$ denote $\lambda_{\pi}$ with
$(N,n)$ reversed, by (4.6),
\begin{eqnarray*}
N^{-2} ( \overline{\lambda}_{2} n \widehat{m}_{3.ijk} + \overline{
\lambda}_{11} N^{2} \widehat{m}_{1.i} \widehat{m}_{2.jk}) \mbox{ is an UE of }
m_{1.i} m_{2.jk},
\end{eqnarray*}
\begin{eqnarray*}
\mbox{ and } N^{-3} ( \overline{\lambda}_{3} n \widehat{m}_{3.ijk} +
\overline{\lambda}_{21} n^{2} \sum^{3} \widehat{m}_{1.i} \widehat{m}_{2.jk} +
\overline{\lambda}_{1^3} n^{3} \widehat{m}_{1.i} \widehat{m}_{1.j} \widehat{m}_{1.k} )
\end{eqnarray*}
is an UE of $ m_{1.i} \ m_{2.j} \ m_{3.k}.$\\
\begin{note} As in (2.4), we can also write $ \mathbf{B}_{r} = \mathbf{B}_{r}
(n,N)$ in the form
\begin{eqnarray}
\mathbf{B}_{r} = \mathbf{D}_{rn}^{-1} \ \mathbf{A}_{r} \ \mathbf{D}_{rN}
\mbox{ where } \mathbf{D}_{rN} = diag ( N^{q_{1}}, \ldots N^{q_{N_{r}}} )
\label{3.15}
\end{eqnarray}
and $q_{k} = q (\pi_{k})$ if the $kth$ element of $m_{(i_{1} \ldots i_{r})}$ is
$m(\pi_{k})_{r.i_{1}\ldots i_{r}}.$

Taking $m_{(i)} = m_{1.i}, m_{(ij)} = ( m_{2.ij}, m_{1.i} \ m_{1.j}
)^{'}$ and $m_{(ijk)}$ as in (4.9) this gives $\mathbf{D}_{1N} = N ,
\mathbf{D}_{2N} = diag (N, N^{2}), \mathbf{D}_{3N} = diag ( N, N^{2},
N^{2}, N^{2}, N^{3} ),$
\begin{eqnarray*}
\mathbf{A}_{1} = \lambda_{1}, \mathbf{A}_{2} = \left( \begin{array}{cc}
\lambda_{1} &\\
\lambda_{2} & \lambda_{11} \end{array} \right), \mathbf{A}_{3} =
\left[ \begin{array}{lllll}
\lambda_{1}& &&& \\
\lambda_{2} & \lambda_{11}& &&\\
\lambda_{2} & 0 & \lambda_{11} && \\
\lambda_{2} & 0 & 0 & \lambda_{11} &\\
\lambda_{3} & \lambda_{12} & \lambda_{12} & \lambda_{12} & \lambda_{1^3}
\end{array} \right]
\end{eqnarray*}
As $ N_{4} =1 + 4 + 3 + 6 + 1 = 15, $ it is not practical ( or
necessary ) to write out any other $\mathbf{A}_{r}$ in full.
\end{note}$\Box$

If $i_1 = \ldots = i_{r}$ then (4.3), (4.4), (4.13), (4.14) reduce to
the univariate formulas.  However to apply them at the intermediate
level when between 2 and $r-1$ of $i_{1, \cdots ,}i_{r}$ are equal,
one needs to be clear what $ \sum^{P (\pi)} m (\pi)_{r.i_{1}\ldots
i_{r}}$ means, and so what $\sum^{P(\pi)} mu (\pi)_{r.i,\ldots,}$
means.  This can be done by reinterpreting the expressions for $ E
S_{i_{1}} \cdots S_{i_{r}}$ in Appendix E.

For example the contribution from $\pi = (13) $ to $ E S_{i_{1}}
\cdots S_{i_{4}}$ is \\ $\lambda (13) \sum^{4} \ s_{i_{1} + i_{2}
+i_{3}} \ s_{i_{4}}$ , and to $ E S_{i_1}^{2} S_{i_{2}} S_{i_{3}}$ is
$\lambda (13) ( \sum^{2} \ s_{2i_{1} + i_{2}} s_{i_{3}} + 2
s_{i_{1}+i_2+ i_{3}} s_{i_{1}} ),$
\begin{eqnarray*}
\mbox{ so } \sum^{P(13)} m (13)_{4.i_{1}i_{1}i_{2}i_{3}} = \sum^{2}
m_{3.2i_{1} + i_{2}} m_{1.i_{3}} + 2 m_{3.i_{1}+i_{2}+i_{3}} m_{1.i_{1}}.
\end{eqnarray*}
In this way the expressions in Appendix E can be interpreted to give $ E
\widehat{m} ( \pi)_{r.i_{1}\ldots i_{r}}$ or an UE for $ m (\pi)_{r.i_{1} \ldots
i_{r}} :$ \\
Set $ \mathbf{i}_{a} = (i_{1,} \ldots,i_{a}), \mathbf{j}_{b} = (j_{1,}
\ldots,j_{b}), \ldots $.
Then
\begin{eqnarray}
E  n \widehat{m}_{a.\mathbf{i}_{a}} &=& \lambda_{1} N m_{a.\mathbf{i}_{a}},
\nonumber\\
E  n^{2} \widehat{m}_{a.\mathbf{i}_{a}} \widehat{m}_{b.\mathbf{j}_{b}} &=& \lambda_{2} N m_{a+b.\mathbf{i}_{a} \mathbf{j}_{b}} + \lambda_{11} N^{2} m_{a.\mathbf{i}
_{a}} m_{b.\mathbf{j}_{b}} , \nonumber\\
En^{3} \widehat{m}_{a.\mathbf{i}_{a}} \widehat{m}_{b.\mathbf{j}_{b}} \widehat{m}_{c.
\mathbf{k}_{c}} &=& \lambda_3 N m_{a+b+c.\mathbf{i}_{a} \mathbf{j}_{b}
\mathbf{k}_{c}} \nonumber\\
&& + \lambda_{2 1}{N}^2 \sum^{2} m_{a+b.\mathbf{i}_{a}\mathbf{j}_{b}} m_{c.
\mathbf{k}_{c}} + \lambda_{1^{3}} N^{3} m_{a.\mathbf{i}_{a}} m_{b.
\mathbf{j}b} m_{c.\mathbf{k}_{c}}, \nonumber\\
En^{3} \widehat{m}_{a.\mathbf{i}_{a}}^2 \widehat{m}_{b.\mathbf{j}b} &=& \lambda_{3}
N m_{2 m+b.\mathbf{i}_{a}\mathbf{i}_{a}\mathbf{j}_b}\nonumber\\
& &+ \lambda_{2 1} N^{2} (m_{2a.\mathbf{i}_{a}\mathbf{i}_{a}}
m_{b.\mathbf{j}_{b}} + 2 m_{a+b.\mathbf{i}_{a}\mathbf{j}_{b}}
m_{a.\mathbf{i}_{a}} )+ \lambda_{1^{3}} N^{3} m_{a.\mathbf{i}_{a}}^{2}
m_{b.\mathbf{j}_{b}}, \nonumber
\end{eqnarray}
and so on.

We end this section with UEs for multinomial parameters.\\
\begin{example}
Suppose $F$ is multinomial $_k (N, \mathbf{p})$, that is for
$1 \leq i \leq k, Np_{i}$ of \{ $x_{1}, \ldots , x_{N}$ \} equal $e_{ik}$,
the $i$th unit vector in $R^{k}.$

Take $2 \leq r \leq k$.  Then $m_{r.i_{1}\ldots i_{r}} = \left\{ \begin{array}
{lll} 0 & \mbox{ if } & \mbox{ any two of $i_{1,}\ldots,i_{r}$ differ }. \\
p_{i_{1}} & \mbox{ if } & i_{1} = \ldots = i_{r}. \end{array}\right.$

Set $\widehat{p}_{i} = \widehat{m}_{1.i}.$

Taking $a=1,b=2, \ldots, e=5$ in Appendix E, discarding all terms with + as
a subscript, and setting $S_{a} = n \widehat{p}_{a}, s_{/a} = Np_a,$
we have
\begin{eqnarray*}
E \widehat{p}_{1} \cdots \widehat{p}_{r} = n^{-r} \lambda_{1^{r}} N^{r} p_{1} \cdots
p_{r} \mbox{ for } r \geq 1,
\end{eqnarray*}
\begin{eqnarray*}
E \widehat{p}_{1}^{2} \widehat{p}_{2} & \cdots & \widehat{p}_{r-1} = n^{-r} (
\lambda_{1^{ r-2}2} + \lambda_{1^{r}} N p_{1} ) N^{r-1} p_{1} \cdots p_{r-1}
\mbox{ for }  r \geq 3, \\
E \widehat{p}_{1}^{3} \widehat{p}_{2} &\cdots& \widehat{p}_{r-2} = n^{-r} ( \lambda_{1^{
r-3}3} + 3 \lambda_{1^{r-2}2} N p_{1}
+ \lambda_{1^r} N^{2} p_{1}^{2})N^{r-2} p_{1} \cdots p_{r-2}, \\
E \widehat{p}_{1}^{2} \widehat{p}_{2}^2 &=& n^{-4}(\lambda_{22}+\lambda_{1^22}
N(p_1+p_2)+\lambda_{1^4}N^2p_1p_2)N^2p_1p_2,\\
E \widehat{p}_1^4 \widehat{p}_2 &=& n^{-5}( \lambda_{14}+4 \lambda_{1^23}Np_1+3
\lambda_{12^2}Np_1+6 \lambda_{1^32}N^2p_1^2+
\lambda_{1^5}N^3p_1^3)N^2p_1p_2, \\
E \widehat{p}_1^3 \widehat{p}^2_2 &=& n^{-5}  ( \lambda_{23}+ \lambda_{1^23}Np_2+3
\lambda_{12^2}Np_1+ \lambda_{1^32}N^2(p_1^2+
3p_1p_2) \\
&& +\lambda_{1^5}N^3p_1^2p_2)N^2p_1p_2, \\
E  \widehat{p}_{1}^{5} \widehat{p}_{2} &=&
n^{-6}(\lambda_6+\lambda_{15}+5\lambda
_{1^24}Np_1+10\lambda_{123}Np_1+10\lambda_{1^33}N^2p_1^2,\\
&&+15\lambda_{1^22^2}N^2p_1^2+10\lambda_{1^42}N^3p_1^3+\lambda_{1^5}N^4p_1^4)
N^2p_1p_2, \\
E  \widehat{p}_{1}^{4} \widehat{p}_{2}^{2} &=&
n^{-6}(\lambda_{24}+\lambda_{1^24}Np_2
+4\lambda_{123}Np_1+3\lambda_{2^3}Np_1+4\lambda_{1^33}N^2p_1p_2 \\
&&+3\lambda_{1^22^2}N^2p_1(2p_1+p_2)+\lambda_{1^42}N^3p_1^2(p_1+6p_2)+
\lambda_{1^6}N^4 p_1^3p_2)N^2p_1p_2, \\
E  \widehat{p}_{1}^{3} \widehat{p}_{2}^{3} &=&
n^{-6}(\lambda_{3^2}+3\lambda
_{123}N(p_1+p_2)+\lambda_{1^33}N^2(p_1^2+p_2^2)+9\lambda_{1^22^2}N^2p_1p_2
\\
&&+
3\lambda_{1^42}N^3p_1p_2(p_1+p_2)+\lambda_{1^6}N^4p_1^2p_2^2)N^2p_1p_2,\\
E \widehat{p}_{1}^{4} \widehat{p}_{2} \widehat{p}_{3} &=& n^{-6}(\lambda_{1^24}+4
\lambda_{1^33}Np_1+3\lambda_{1^22^2}Np_1+6\lambda_{1^42}N^2p_1^2 \\
&&+\lambda_{1^6} N^3 p_1^3)N^3p_1p_2p_3, \\
E  \widehat{p}_{1}^{3} \widehat{p}_{2}^{2} \widehat{p}_{3} &=& n^{-6}(\lambda_{123}
+\lambda_{1^33}Np_2+3\lambda_{1^22^2}Np_1 +\lambda_{1^42}N^2p_1
(p_1+p_2)\\
&&+\lambda_{1^6}N^3p_1^2p_2)N^3p_1p_2p_3, \\
E  \widehat{p}_{1}^{2} \widehat{p}_{2}^{2} \widehat{p}_{3}^{2} &=& n^{-6}
(\lambda_{2^3} +\lambda_{1^22^2}N(p_1+p_2+p_3)
+\lambda_{1^42}N^2(p_1p_2+p_1p_3+p_2p_3), \\
&&+
\lambda_{1^6}N^3p_1p_2p_3)N^3p_1p_2p_3, \\
E  \widehat{p}_{1}^{2} \widehat{p}_{2}^{2} \widehat{p}_{3} \widehat{p}_{4} &=&
n^{-6}(\lambda_{1^22^2}+\lambda_{1^42}N(p_1+p_2)+\lambda_{1^6}N^2p_1p_2)N^4 p_1p_2p_3p_4.
\end{eqnarray*}
If the expression for $ E \widehat{p}_{\pi_1} \widehat{p}_{\pi_{2}} \cdots$ is $a_{Nn} ( \pi, p),$ then $ a_{nN} ( \pi, \widehat{p})$
 is an UE for $ p_{\pi_1} p_{\pi_2} \cdots .$

This generates a new set of matrices $ A=A (N,n)$ satisfying the invariance
principle $ A(N, n)^{-1} = A (n, N)$ : \\
from $( 1 \ldots r), A = \lambda_{1}^{r};\\
\mbox{ from } (1 \ldots r-1), (1^{2} 2 \ldots r-1), A = \left( \begin{array}
{ll} \lambda_{1^{r-1}}& \\ \lambda_{1^{r-2}2} & , \lambda_{1^{r}} \end{array}
\right) ;\\
 \mbox{ from } (1 \ldots r-2 ), (1^{2} 2 \ldots r-2), (1^{3}2 \ldots r-2 ),
A= \left( \begin{array}{lll} \lambda_{{1}^{ r-2}} & & \\ \lambda_{1^{r-3}3}&,
 \lambda_{1^{r-1}} & \\ \lambda_{1^{r-3}3} & ,  3 \lambda_{1^{r-2}2}, &
\lambda_{1^r} \end{array} \right); \\
 \mbox{ from }  (1 2, 1^2 2, 1^3 2), A = \left( \begin{array}{lll} \lambda_{2} & & \\ \lambda_{12} & \lambda_{1^3} & \\ \lambda_{13} & 3 \lambda_{1^22} &
\lambda_{1^4} \end{array} \right) ; \\
\mbox{ from } (1 2 , 1^2 2, 1 2^2, 1^2 2^2 ), A= \left( \begin{array}{llll}
\lambda_{2} &&& \\ \lambda_{12} & \lambda_{1^3}&& \\ \lambda_{12} & 0 &
\lambda_{1^3} & \\ \lambda_{22} & \lambda_{1^22} & \lambda_{1^22} &
 \lambda_{1^4} \end{array} \right) ; \\
\mbox{ from }  (1 2, 1^2 2, 1^3 2, 1^4 2 ), A = \left( \begin{array}{llll}
\lambda_{2} &&& \\ \lambda_{12} & \lambda_{1^{3}} && \\ \lambda_{13} & 3 \lambda_{1^22} & \lambda_{1^4} & \\ \lambda_{14} & 4 \lambda_{1^23} + \lambda_{12^2} & 6 \lambda_{1^{3}2} & \lambda_{1^{5}} \end{array} \right) ; \\
\mbox{ from }  ( 12, 12^2, 1^2 2, 1^3 2, 1^2 2^2, 1^3 2^2 ), \\
 A = \left( \begin{array}{llllll} \lambda_{2} &&&&& \\
\lambda_{12} & \lambda_{1^3} &&&& \\
\lambda_{12} & 0 & \lambda_{1^3} &&& \\
\lambda_{13} & 0 & 3\lambda_{1^22} & \lambda_{1^4} && \\
\lambda_{22} & \lambda_{1^22} & \lambda_{1^22} & 0 & \lambda_{1^{4}} & \\
\lambda_{23} & \lambda_{1^23} & 3 \lambda_{12^2} & \lambda_{1^32} & 3 \lambda_{1^32} & \lambda_{1^{5}} \end{array} \right ) ; \\
 \mbox{ from } (123, 12^2 3, 1^2 23, 1^2 2^2 3 ), A= \left( \begin{array}{llll}
\lambda_{1^3} &&& \\
\lambda_{1^22} & \lambda_{1^4} && \\
\lambda_{1^22} & 0 & \lambda_{1^4} & \\
\lambda_{12^2}  & \lambda_{1^{3}2} & \lambda_{1^{3}2}  & \lambda_{1^{5}}
\end{array} \right) ; \\$
and so on. $\Box$ \end{example}

Raghunandanan and Srinivasan (1973) give multivariate analogs of
$E\widehat{\mu}_r$ for $5 \leq r \leq 8$ in terms of symmetric functions with
tables to express these in terms of noncentral moments, and
$E(\overline{X}-\mu)^r$ for $r=5,6.$  The latter agree with Sukhatme except
have $-5e_3$ where Sukhatme has $e_2-6e_3$ in the coefficient of $\mu_3^2$
in $E(\overline{X}-\mu)^6.$  Sukhatme's version is the correct one since
$g_9n^6$ is $C_{33}=C_3\otimes C_3$ in the notation of Dwyer and Tracy.
They also give the multivariate analogs of the UE of $E(\overline{X}-\mu)^4$
in terms of $\widehat{\mu_4}$ and
$\sum_{i\neq j}(X_i-\overline{X})^2(X_j-\overline{X})^2$, and the UE of
$E(\overline{X}-\mu)^5$ in terms of $\widehat{\mu}_5$ and
$\sum_{i\neq j}(X_i-\overline{X})^3(X_j-\overline{X})^2$.\\

It would be useful if Pierce's transformations (to obtain results for
finite $N$ from results for $N=\infty$) were extended to cover multivariate
problems like (3.2).

\section{UEs of products of central moments}

Here we express the joint central moments $ \mu (\pi)$ of (1.3) as a linear
combination of $ m (\pi)$ and so obtain
\begin{eqnarray*}
& E \widehat{\mu} (\pi_{1} ) \widehat{\mu} ( \pi_{2}) \ldots \\
\mbox{ and } &   \mbox{ UEs  of } \mu (\pi_{1})  \  \mu (
\pi_{2}) \ldots \mbox{ for total order } \\
& r = r ( \pi_{1}) + r ( \pi_{2} ) \ldots  \leq 6.
\end{eqnarray*}
We derive coefficients $ M_{\pi.\pi \_^{'}} $ such that
\begin{eqnarray}
\mu ( \pi ) = \sum^{r}_{\pi_{-}^{'}}  \  M_{ \pi  . \pi_{-}^{'}}mu(\pi'_{-}).
\label{4.1}
\end{eqnarray}
This enables us to express $ \mu ( \pi_{1}) \mu (\pi_{2}) \ldots$ in the form
\begin{eqnarray}
D = \sum_{\pi \_}^{r} D_{ \pi \_} mu ( \pi \_ ). \label{4.2}
\end{eqnarray}
By Section 2, an UE of D is then
\begin{eqnarray}
&& \widetilde{D}  =  \sum_{\pi_{-}^{'}}^{r}  \ D^{*}_{ \pi_{-}^{'}}  \ \widehat{mu} (
\pi_{-}^{'}) \label{4.3} \\
\mbox{ where } & & D^{*}_{ \pi_{-}^{'}} = \sum_{\pi \_}^{r} D_{ \pi \_}
C^{\pi \_ , \pi \_^{'}} \mbox{ and }
\left ( C^{\pi_{-}.\pi_{-}} \right )=C_{--r}^{-1}\label{4.6}
\end{eqnarray}
We begin with $ E \widehat{\mu} ( \pi )$ and give the UE of $ \mu (\pi)$ in the order
\begin{eqnarray*}
\pi = 1^{r}, \ 2^{r}, \ 3^2, \ 1^{r} 2, \  1^{r} 3, \  1^{r} 4,
\ 1^r2^2, \ 23, \ 123.
\end{eqnarray*}
We omit $ \pi = 15$ as Sukhatme omitted $\mu_{(51)}.$

As noted in Section 2,
\begin{eqnarray*}
\mu (1^{r}) = E ( \overline{X} - \mu )^{r} = \sum_{\pi \_}^{r}
C_{1^{r}.\pi \_} mu (\pi \_),
\end{eqnarray*}
\begin{eqnarray}
\mbox{ so in (5.1), } M_{1^{r} . \pi \_} = C_{1^{r}  .  \pi \_}
\label{4.5}
\end{eqnarray}
The $D^{*}$ needed in (5.3) for the UE of $D = \mu (1^{r})$ are as
follows.\\ \\
\begin{eqnarray}
\mu(1^2)  :  D_{2}^{*} &=& (N-n)N^{-1}(n-1)^{-1},\nonumber\\
\mu(1^3)  :  D_{3}^{*} &=& (N-2n)(N-n)N^{-2}(n-1)_2^{-1},\nonumber\\
\mu(1^4)  :  D_{4}^{*} &=& -(N-n)n^{-1}(2N^2n-4N^2+6Nn-3n^3\nonumber\\
&&-3n^2)(n-1)_3^{-1}N^{-3},\\
D_{2^2}^{*} &=& 3(N-n)n^{-1}(N^2n^2-4N^2n+4N^2-Nn^3-6Nn\nonumber\\
&& +6Nn^2+3n^2 -3n^3)(n-1)_3^{-1}N^{-3},\nonumber\\
\mu(1^5)  :  D_{5}^{*} &=& -(N-n)(N-2n)n^{-1}(9N^2n-15N^2+20Nn-8Nn^2 \nonumber\\
&& -10n^2-2n^3)(n-1)_4^{-1}N^{-4},\\
D_{23}^{*} &=& 10(N-n)(N-2n)n^{-1}(N^2n^2-3N^2n+3N^2-Nn^3\nonumber\\
&& -4Nn -4Nn+4Nn^2-2n^3+2n^2)(n-1)_4^{-1}N^{-4}\nonumber
\end{eqnarray}
\begin{eqnarray}
\mu(1^6) :  D_{6}^{*} &=&  (N-n)n^{-3}(-24N^4+6N^4n^3-84N^4n^2+126N^4n \nonumber\\
&& +120N^3n^3-345N^3n^2+60N^3n+45N^3n^4-80N^2n^2 \nonumber\\
&& +50N^2n^4+400N^2n^3 -130N^2n^5+75Nn^6-225Nn^4 \nonumber\\
&& +60Nn^3-150Nn^5-20n^4+80n^6+5n^7+55n^5)(n-1)_5^{-1}N^{-5},\\
D_{24}^{*} &=&   -15(N-n)n^{-3}(24N^4+2N^4n^4-24N^4n^3+82N^4n^2 \nonumber\\
&& -102N^4n-2N^3n^5-220N^3n^3+285N^3n^2-60N^3n\nonumber\\
&&+57N^3n^4-3N^2n^6+80N^2n^2+220N^2n^4-320N^2n^3\nonumber\\
&& -43N^2n^5+3Nn^7+6Nn^6-90Nn^5-60Nn^3+165Nn^4\nonumber\\
&& +5n^7+10n^6+20n^4-35n^5)(n-1)_5^{-1}N^{-5}\nonumber \\
D_{3^2}^{*} &=& 10(N-n)n^{-3}(24N^4+N^4n^4-12N^4n^3\nonumber\\
&&+57N^4n^2-90N^4n-5N^3n^5 -165N^3n^3+255N^3n^2\nonumber \\
&& -60N^3n+39N^3n^4+8N^2n^6+80N^2n^2+175N^2n^4\nonumber\\
&& -280N^2n^3-47N^2n^5-4Nn^7+24Nn^6-75Nn^5\nonumber\\
&& -60Nn^3+135Nn^4 -5n^7 +10n^6+20n^4-25n^5)(n-1)_5^{-1}N^{-5}  \nonumber\\
D_{2^3}^{*} &=&
15(24N^4+N^4n^4-12N^4n^3+46N^4n^2-66N^4n\nonumber\\
&& -2N^3n^5+30N^3n^4+180N^3n^2 -60N^3n-122N^3n^3\nonumber\\
&& -210N^2n^3+N^2n^6-27N^2n^5+127N^2n^4\nonumber\\
&& +80N^2n^2+120Nn^4-60Nn^3-59Nn^5+9Nn^6\nonumber\\
&& +10n^6+20n^4-30n^5)n^{-2}(N-n)(n-1)_5^{-1}N^{-5}. \nonumber\\
\end{eqnarray}
Now in the expansion for $\mu_{r}$ in terms of $ \{ \mu_{i}^{'} \}$, the
coefficient of $\mu^{'}_r$ is 1.  So in the expansion of $\mu (\pi)$, the
coefficient of $ m (\pi)$ is 1.  So
\begin{eqnarray}
M_{\pi .  r} = C_{\pi  .  r} \label{4.9}
\end{eqnarray}
For example $\mu (2^{2}) = E \widehat{\mu}_{2}^{2} - ( E \widehat{\mu}_{2} )^{2}
= R H S $ (5.1) with
\begin{eqnarray*}
M_{2^2.4} &=& C_{2^2.4} , \\
M_{2^{2}.2^{2}} &=& C_{2^{2} .  2^{2}} - C_{2 . 2}^{2} \\
&=& -N(n-1)(N^2n-3N^2+6N-3n-3)(N-n)n^{-3}(N-1)^{-2}(N-2)_2^{-1},
\end{eqnarray*}
So by (5.4) the UE of $ \mu (2^{2})$ is (5.3) with
\begin{eqnarray}
 D_{4}^{*} &=& 0,\nonumber\\
 D_{2^{2}}^{*} &=& 1.\nonumber \\
\mbox{ Since } \mu (2^{3} ) &=& E \widehat{\mu}_{2}^{3} -3 E \widehat{\mu}_{2}^{2}  \ E \widehat{\mu}_{2} + 2 ( E \widehat{\mu}_{2} )^{3}, \nonumber\\
 M_{2^{3}.24} &=& C_{2^{3}.24} -3 C_{2^{2}.4} C_{2.2} \nonumber\\
 &=& -3N(n-1)(N-n)(
N^4n^2
+5N^4
-6N^4n
\nonumber\\
& &
-2N^3n^3
+5N^3
+13N^3n^2
-53N^2n^2
+20N^2n
\nonumber\\
& &
-45N^2
+2N^2n^3
+55N
+30Nn
+35Nn^2
\nonumber\\
& &
+20Nn^3
-20n
-20n^2
-20n^3
-20
)n^{-5}(N-1)^{-2}(N-2)_4^{-1},\nonumber\\
 M_{2^{3}.3^{2}} &=& C_{2^{3}.3^{2}}, \nonumber\\
 M_{2^{3}.2^{3}} &=& C_{2^{3}.2^{3}} -3 C_{2^{2}.2^{2}} C_{2.2} + 2
C_{2.2}^{3} \nonumber\\
 &=& 2N^2(n-1)(N-n)(N^4n^2+15N^4-12N^4n\nonumber\\
&& -75N^3-2N^3n^3+25N^3n^2+30N^3n\nonumber\\
& & -3N^2n^3+135N^2-139N^2n^2+51N^2n\nonumber\\
& & +167Nn^2-144Nn-105N+56Nn^3+75n-75n^3\nonumber\\
& & -30n^2+30)n^{-5}(N-1)^{-3}(N-2)_4^{-1},
\nonumber\\
\mbox{ So }  & & \nonumber\\
 D_{6}^{*} &=& 0, \nonumber\\
 D_{24}^{*} &=& 0, \nonumber\\
 D_{3^{2}}^{*} &=& 0, \nonumber\\
 D_{2^{3}}^{*} &=& 1.\nonumber \\
\mu (3^{2}) :  M_{3^{2}.3^{2}} &=& C_{3^{2}.3^{2}} - C_{3.3}^{2} \nonumber\\
 &=& -N(n-1)_2(N-n)(20N^5+N^5n^2-12N^5n-100N^4\nonumber\\
&&+4N^4n^2+42N^4n+260N^3-41N^3n^2+24N^3n\nonumber\\
&& +40N^2n^2-460N^2-282N^2n+420Nn+100Nn^2\nonumber\\
&& +440N-240n-160-80n^2)n^{-5}(N-1)^{-2}(N-2)^{-2}(N-3)_3^{-1},\nonumber\\
M_{3^2.2^3}  &=& C_{3^2.2^3}, \nonumber \\
D_{6}^{*} & =& 0,\ D_{24}^{*} = 0, \ D_{3^{2}}^{*} =1, \ D_{2^{3}}^{*} = 0
\nonumber \\
\mbox{ Also }  \mu (1^{r} k) &=& E ( \overline{x} - \mu )^{r} ( \widehat{\mu}_{k} - E \widehat{\mu}_{k} ) \nonumber\\
 &=& \sum_{\pi \_}^{r+k} C_{1^{r}k. \pi \_} m ( \pi \_) - \sum_{\pi \_}^{r}
C_{1^{r} . \pi \_} m (\pi \_) \sum_{\pi \_}^{k}
C_{k . \pi' \_} m (\pi' \_). \nonumber\\
\mbox{ So }  \mu (12) &=& M_{12 . 3} \mu_{3} \mbox{ for } M_{12 . 3} = C_{12.3}, \mbox{ as given by (4.9), and } \nonumber\\
 D_{3}^{*} &=& (N-n)N^{-1}(n-2)^{-1}.\nonumber \\
\mu (1^{2} 2): M_{1^{2}2.2^{2}} &=& C_{1^{2}2.2^{2}} - C_{1^{2}.2} C_{2.2} , \nonumber\\
&=& -N(n-1)(N-n)n^{-3}(N-1)^{-2}(N-2)_2^{-1}(3N^2-4Nn	 \nonumber\\
&& -6N+3+6n),\nonumber\\
 D_{4}^{*}= &=&(N-n)n^{-1}(2N-n-n^2)(n-2)_2^{-1}N^{-2}, \nonumber\\
 D_{2^{2}}^{*} &=&  (N-n)n^{-1}(6N+Nn^2-6Nn-3n+3n^2)(n-2)_2^{-1}N^{-2}.
\nonumber \\
\mu (1^{3} 2): M_{1^{3}2.23 } &=& C_{1^{3}2.23} - C_{1^{3}.3} C_{2.2}
\nonumber \\
  &=& N(n-1)(N-n)(-10N^3+3N^3n-3N^2n^2\nonumber\\
& & +15N^2n+20N^2-9Nn^2-36Nn-10N+24n^2\nonumber\\
&& +12n)n^{-4}(N-1)^{-2}(N-2)_3^{-1},\nonumber\\
D_{5}^{*}= &=& -(N-n)(3N^2n-9N^2+15Nn-3Nn^2-5n^2\nonumber\\
& & -n^3)n^{-1}(n-2)_3^{-1}N^{-3},\nonumber\\
D_{23}^{*} &=&  (N-n)(4N^2n^2-18N^2n+18N^2-5Nn^3-30Nn\nonumber\\
& &+30Nn^2-10n^3+10n^2)n^{-1}(n-2)_3^{-1}N^{-3}.\nonumber \\
\mu (1^{4} 2) : M_{1^{4}2.24} &=& C_{1^{4}2.24} - C_{1^{4}.4} C_{2.2} \nonumber\\
 &=& N(n-1)(N-n)(6N^4n
-15N^4
-18N^3n^2
+70N^3n
\nonumber\\
& &
-15N^3
-84N^2n^2
-120N^2n
+12N^2n^3
+135N^2
+234Nn^2
\nonumber\\
& &
-40Nn
+36Nn^3
-165N
+60n
-60n^2
-120n^3
\nonumber\\
& &
+60
)n^{-5}(N-1)^{-2}(N-2)_4^{-1},
\nonumber\\
M_{1^{4}2.3^{2}} &=& C_{1^{4}2.3^{2}}, \nonumber\\
M_{1^{4}2.2^{3}} &=& C_{1^{4}2.2^{3}} - C_{1^{4}.2^{2}} C_{2.2} \nonumber\\
 &=& -6N^2(n-1)(N-n-1)(N-n)(+3N^2n-5N^2-4Nn^2\nonumber\\
& & -4Nn+15N+10n^2-5n-10)n^{-5}(N-1)^{-2}(N-2)_4^{-1},\nonumber\\
D_{6}^{*} &=& -(N-n)n^{-3}(16N^3+4N^3n^3+8N^3n^2-52N^3n+24N^2n^3\nonumber\\
& & -18N^2n^4+90N^2n^2-24N^2n-56Nn^3+14n^5N-46Nn^4\nonumber\\
&& +16Nn^2+11n^4 -4n^3+16n^5+n^6)(n-2)_4^{-1}N^{-4},\nonumber\\
D_{24}^{*} &=& -(N-n)n^{-3}(240N^3+2N^3n^4-72N^3n^3+340N^3n^2\nonumber\\
& & -540N^3n+6N^2n^5+90N^2n^4+990N^2n^2-360N^2n\nonumber\\
&& -600N^2n^3-600Nn^3-9Nn^6-3Nn^5+300Nn^4+240Nn^2\nonumber\\
&& +105n^4-60n^3-30n^5-15n^6)(n-2)_4^{-1}N^{-4} , \nonumber\\
D_{3^2}^{*} &=& 2(N-n)n^{-3}(80N^3+2N^3n^4-16N^3n^3+70N^3n^2\nonumber\\
& & -140N^3n-6N^2n^5+42N^2n^4+270N^2n^2-120N^2n\nonumber\\
&& -150N^2n^3-160Nn^3+4Nn^6-29Nn^5+85Nn^4+80Nn^2\nonumber\\
&& +25n^4-20n^3-10n^5+5n^6)(n-2)_4^{-1}N^{-4},\nonumber\\
D_{2^{3}}^{*} &=& 3(N^3n^4-16N^3n^3+76N^3n^2-140N^3n+80N^3\nonumber\\
& & -120N^2n-N^2n^5+24N^2n^4-128N^2n^3+240N^2n^2\nonumber\\
&& -150Nn^3+80Nn^2+69Nn^4-9Nn^5-10n^5\nonumber\\
&& +30n^4-20n^3)n^{-2}(N-n)(n-2)_4^{-1}N^{-4}. \nonumber \\
\mu (13): M_{13.2^{2}} &=& C_{13.2^2}, \nonumber\\
 D_{4}^{*} &=& (N-n)(n+1)N^{-1}n^{-1}(n-3)^{-1}, \nonumber\\
 D_{2^{2}}^{*} &=& -3(N-n)(n-1)N^{-1}n^{-1}(n-3)^{-1}. \nonumber \\
\mu (1^{2}3): M_{1^{2}3.23} &=& C_{1^{2}3.23} - C_{1^{2}.2} C_{3.3} \nonumber\\
 &=& -2(N-n)N(n-1)_2(5N^3-7N^2n-10N^2+14Nn\nonumber\\
& & +5N-4n)n^{-4}(N-1)^{-2}(N-2)_3^{-1},\nonumber\\
 D_{5}^{*} &=& (N-n)N^{-2}n^{-1}(6N-n^2-5n)(n-3)_2^{-1}, \nonumber\\
 D_{23}^{*} &=& (N-n)N^{-2}n^{-1}(Nn^2-12Nn+12N+10n^2\nonumber\\
&& -10n)(n-3)_2^{-1}.
\nonumber
\end{eqnarray}
\begin{eqnarray}
\mu (1^{3} 3):  M_{1^{3}3.24} &=& C_{1^{3}3.24} , \nonumber\\
 M_{1^{3}3.3^{2}} &=& C_{1^{3}3.3^{2}} - C_{1^{3}.3} C_{3.3} \nonumber\\
 &=& -(N-n)N(n-1)_2(+10N^5-33N^4n-50N^4\nonumber\\
&&+24N^3n^2+130N^3+162N^3n-303N^2n-108N^2n^2\nonumber\\
&& -230N^2+132Nn^2+220N+270Nn-120n\nonumber\\
&& -80)n^{-5}(N-1)^{-2}(N-2)^{-2}(N-3)_3^{-1},\nonumber\\
 M_{1^{3}3.2^{3}} &=& C_{1^{3}3.2^{3}}, \nonumber\\
 D_{6}^{*} &=& -(N-n)N^{-3}n^{-3}(12N^2-12N^2n^2-27N^2n\nonumber\\
&&+3N^2n^3-3Nn^4+30Nn^3-12Nn+33Nn^2\nonumber\\
&& -11n^3-n^5+4n^2-16n^4)(n-3)_3^{-1}, \nonumber\\
 D_{24}^{*} &=&  3(N-n)N^{-3}n^{-3}(N^2n^4-3N^2n^3-25N^2n^2\nonumber\\
&& +75N^2n-60N^2+60Nn-Nn^5+6Nn^4+40Nn^3\nonumber\\
&&-105Nn^2+35n^3-20n^2-10n^4-5n^5)(n-3)_3^{-1},\nonumber\\
D_{3^{2}}^{*} &=&(N-n)N^{-3}n^{-3}(120N^2+N^2n^4-6N^2n^3+15N^2n^2\nonumber\\
&& -90N^2n-120Nn-2Nn^5+24Nn^4-60Nn^3+150Nn^2\nonumber\\
&& -50n^3+40n^2+20n^4-10n^5)(n-3)_3^{-1},\nonumber\\
D_{2^{3}}^{*} &=&-3(-60N^2+3N^2n^3-28N^2n^2+75N^2n-3Nn^4+35Nn^3\nonumber\\
&& -90Nn^2+60Nn-10n^4+30n^3-20n^2)n^{-2}N^{-3}(N-n)(n-3)_3^{-1} .
\nonumber
\end{eqnarray}
\begin{eqnarray}
\mu (14):  M_{14.23} &=& C_{14.23} , \nonumber\\
 D_{5}^{*} &=& (n^2-n-3)(N-n)n^{-1}N^{-1}(n-4)^{-1}(n-2)^{-1},  \nonumber\\
 D_{23}^{*} &=&  -2(n-1)(2n-3)(N-n)n^{-1}N^{-1}(n-4)^{-1}(n-2)^{-1}.
\nonumber \\
\mu(1^{2}4):  M_{1^{2}4.24} &=& C_{1^{2}4.24} - C_{1^{2}.2} C_{4.4},
\nonumber\\
 &=& -(n-1)N(N-n)(
9N^4n^2-39N^4n+45N^4+39N^3n^2\nonumber\\
&& -81N^3n-10N^3n^3+45N^3+45N^2n^2-405N^2\nonumber\\
&&-6N^2n^3+237N^2n-435Nn^2+495N\nonumber\\
&& +135Nn-180n+270n^2-180-30n^3\nonumber\\
&& +70Nn^3)n^{-5}(N-1)^{-2}(N-2)_4^{-1},\nonumber\\
 M_{1^{2}4.3^{2}} &=& C_{1^{2}4.3^{2}} , \nonumber\\
 M_{1^{2}4.2^{3}} &=& C_{1^{2}4.2^{3}} - C_{1^{2}.2} C_{4.2^{2}}\nonumber\\
 &=& 6(n-1)N(N-n)(2N^4n^2-12N^4n+9N^3n^2+15N^4\nonumber\\
&& -2N^3n^3-60N^3+27N^3n+75N^2-44N^2n^2+15N^2n\nonumber\\
&& -30Nn+9Nn^2-30N+26Nn^3+30n^2\nonumber\\
&& 30n^3 )n^{-5}(N-1)^{-2}(N-2)_4^{-1},\nonumber
\end{eqnarray}
\begin{eqnarray}
 D_{6}^{*} &=& (
-36N^4n^3
+8N^4n^4
+22N^4n^2
+54N^4n
-24N^4
-N^3n^6
-155N^3n^2
\nonumber\\
&&
-118N^3n
+3N^3n^5
+121N^3n^3
-18N^3n^4
+72N^3
-N^2n^6
+N^27n^5
\nonumber\\
&&
+60N^2n
+233N^2n^3
-126N^2n^4
-5N^2n^2
+64Nn
-96Nn^3
+4Nn^4
\nonumber\\
&&
+32Nn^5+12Nn^6
-304Nn^2
-112n^5
-64n^2
+208n^3
+168n^4
\nonumber\\
&&
-8n^6)n^{-3}(N-n)(n-2)_4^{-1}(N-3)^{-1}N^{-4} ,
\label{4.11}
\end{eqnarray}
\begin{eqnarray}
 D_{24}^{*} &=&  (
110N^4n^4
-204N^4n^3
+450N^4n
-21N^4n^5
+N^4n^6
-360N^4
\nonumber\\
&&
+1080N^3
+6N^3n^6
-6N^3n^5
-186N^3n^4
-675N^3n^2
-690N^3n
\nonumber\\
&&
+687N^3n^3
-408N^2n^4
+423N^2n^3
-3N^2n^6
+135N^2n^5
-615N^2n^2
\nonumber\\
&&
+900N^2n
+240Nn^5
+2400Nn^3
-468Nn^4
+960Nn
-108Nn^6
\nonumber\\
&&
-3600Nn^2
-960n^2
+2160n^3
-1320n^4
\nonumber\\
&&
+120n^6)n^{-3}(N-n)\nonumber\\
&& (n-2)_4^{-1}(N-3)^{-1}N^{-4}, \nonumber\\
\end{eqnarray}
\begin{eqnarray}
 D_{3^{2}}^{*} &=&  -2(
10N^4n^4
-57N^4n^3
+65N^4n^2
+90N^4n
-120N^4
-310N^3n^2
\nonumber\\
&&
+39N^3n^4
+360N^3
-50N^3n
-10N^3n^5
+71N^3n^3
-2N^2n^6
\nonumber\\
&&
+2N^2n^5
+300N^2n
+184N^2n^3
-45N^2n^4
-295N^2n^2
+720Nn^3
\nonumber\\
&&
-124Nn^4
+28Nn^6
+320Nn
-64Nn^5
-1040Nn^2
+160n^5
\nonumber\\
&&
-320n^2
+560n^3
-360n^4
)n^{-3}(N-n)(n-2)_4^{-1}(N-3)^{-1}N^{-4}, \nonumber\\
\end{eqnarray}
\begin{eqnarray}
 D_{2^{3}}^{*} &=&  6(N-n)n^{-2}(2N^4n^4-18N^4n^3+64N^4n^2-105N^4n+60N^4\nonumber\\
&& -180N^3+205N^3n-2n^5N^3+15N^3n^4-31N^3n^3\nonumber\\
&& -32N^3n^2+141N^2n^3+90N^2n-228N^2n^2  -27N^2n^4 \nonumber\\
&& +24Nn^5-160Nn+320Nn^2-72Nn^4-72Nn^3 +200n^4 \nonumber\\
&& -40n^5+160n^2-320n^3)(n-2)_4^{-1}(N-3)^{-1}N^{-4} , \label{4.12} \\
\end{eqnarray}
\begin{eqnarray*}
\mbox{ Also } \mu (1^{r} k^{2} ) = E ( \overline{X} - \mu )^{r} ( \widehat{\mu}
_{k}^{2} - 2 \widehat{\mu}_k E \widehat{\mu}_{k} + ( E \widehat{\mu}_{k} )^{2})
\end{eqnarray*}
\begin{eqnarray}
\mu (12^{2}): M_{12^{2}.23} &=& C_{12^{2}.23} - 2 C_{12.3} C_{2.2} \nonumber\\
 &=& -2(n-1)N(N-n)(
3N^3n-5N^3-4N^2n^2\nonumber\\
&& +10N^2+8Nn^2-11Nn-5N+2n^2\nonumber\\
&& +2n)n^{-4}(N-1)^{-2}(N-2)_3^{-1},\nonumber\\
D_{5}^{*} &=& -(n-1)n^{-1}N^{-1}(N-n)(n-4)^{-1}(n-2)^{-1}, \nonumber\\
D_{23}^{*} &=& 2(n^2-3n+1)n^{-1}N^{-1}(N-n)(n-4)^{-1}(n-2)^{-1}. \nonumber
\end{eqnarray}
\begin{eqnarray}
\mu (1^{2}2^{2}):  M_{1^{2}2^{2}.24} &=& C_{1^{2}2^{2}.24} -2 C_{1^{2}2.4}
C_{2.2} \nonumber\\
 &=&  (n-1)N(N-n)(
-12N^4n
+N^4n^2
+15N^4
+31N^3n^2
\nonumber\\
&&
-35N^3n
-N^3n^3
+15N^3
-22N^2n^3
-135N^2
\nonumber\\
&&
-29N^2n^2
+90N^2n
+65Nn
+61Nn^3
+165N
-109Nn^2
\nonumber\\
&&
+10n^2
-60n
-60
+10n^3
)n^{-5}(N-1)^{-2}(N-2)_4^{-1},\nonumber\\
 M_{1^{2} 2^{2}.3^{2}} &=& C_{1^{2}2.3^{2}} , \nonumber\\
 M_{1^{2}2^{2}.2^{3}} &=& C_{1^{2}2^{2}.2^{3}} -2 C_{1^{2}2.2^{2}} C_{2.2}
+ C_{1^{2}.2} C_{2.2}^2 \nonumber\\
 &=&  -N^2(n-1)(N-n)(
-21N^4n
+30N^4
+N^4n^2
+44N^3n^2
\nonumber\\
&&
-N^3n^3
-150N^3
+45N^3n
+270N^2
-25N^2n^3
+93N^2n
\nonumber\\
&&
-194N^2n^2
+127Nn^3
-210N
-237Nn
+178Nn^2
+60
\nonumber\\
&&
+120n
-5n^2
-125n^3
)n^{-5}(N-1)^{-3}(N-2)_4^{-1},\nonumber\\
 D_{6}^{*} &=&-(N-n)(n-1)(
+2Nn^3-4Nn^2-10Nn +8N-4n\nonumber\\
&& +7n^2+8n^3-3n^4)N^{-2}n^{-3}(n-2)_4^{-1},\nonumber\\
 D_{24}^{*} &=&  (N-n)(120N+16Nn^4-102Nn^3+200Nn^2\nonumber\\
&& -150Nn-145n^3+105n^2-60n+n^5-2n^6\nonumber\\
&&+53n^4)N^{-2}n^{-3}(n-2)_4^{-1} ,\nonumber\\
D_{3^{2}}^{*} &=&  2(N-n)(-40N+Nn^5-9Nn^4+31Nn^3-45Nn^2\nonumber\\
&& +30Nn+40n^3-25n^2+20n+8n^5-n^6\nonumber\\
&& -26n^4)N^{-2}n^{-3}(n-2)_4^{-1},\nonumber\\
D_{2^{3}}^{*} &=& (-120N-17Nn^4+92Nn^3-208Nn^2+210Nn+Nn^5\nonumber\\
&& +60n-48n^4-120n^2+6n^5+126n^3)n^{-2}N^{-2}(N-n)(n-2)_4^{-1}. \nonumber\\
\mbox{Also, }
\mu (23): M_{23.23} &=& C_{23.23} - C_{2.2} C_{3.3}, \nonumber\\
&=& -2(N-n)N^2(n-1)_2(2N^2n-5N^2+10N-\nonumber\\
&& -n)n^{-4}(N-1)^{-2}(N-2)_3^{-1},\nonumber\\
D_{5}^{*} &=& 0 ,\ D_{23}^{*}= 1.\nonumber
\end{eqnarray}
\begin{eqnarray*}
\mbox{ and } \mu(123) = E ( \overline{X} - \mu ) ( \widehat{\mu}_{2} \widehat{\mu}_{3} - \widehat{\mu}_{2} E \widehat{\mu}_{3} - \widehat{\mu}_{3} E \widehat{\mu}_{2} + E \widehat{\mu}_{2} E \widehat{\mu}_{3} )
\end{eqnarray*}
\begin{eqnarray}
\mbox{ so }
 M_{123.24} &=& C_{123.24} - C_{13.4} C_{2.2}\nonumber\\
 &=& -(N-n)N(n-1)_2(
-15N^4
+9N^4n
-15N^3
-2N^3n
-13N^3n^2
\nonumber\\
&&
+135N^2
+42N^2n^2
+9N^2n
-165N
-130Nn
-35Nn^2
+30n^2
\nonumber\\
&&
+90n+60
)n^{-5}(N-1)^{-2}(N-2)_4^{-1},\nonumber\\
M_{123.3^{2}} &=& C_{123.3^{2} } - C_{12.3} C_{3.3}\nonumber\\
&=& -(N-n)N(n-1)_2(
-10N^5
+4N^5n
-3N^4n
+50N^4
-5N^4n^2
\nonumber\\
&&
-130N^3
+18N^3n^2
-62N^3n
+230N^2
+13N^2n^2
\nonumber\\
&&
-90Nn^2
+195N^2n
-230Nn
-220N
+80
+120n
\nonumber\\
&&
+40n^2
)n^{-5}(N-1)^{-2}(N-2)^{-2}(N-3)_3^{-1},\nonumber\\
M_{123.2^{3}} &=& C_{123.2^{3}} - C_{13.2^{2}} C_{2.2}\nonumber\\
&=&  3N^2(n-1)_2(N-n)(-10N^2+4N^2n-5Nn^2-Nn\nonumber\\
&& +30N+15n^2-25n-20)n^{-5}(N-1_5^{-1},\nonumber\\
D_{6}^{*} &=& -(1+n)(n-1)^2(N-n)n^{-3}N^{-1}(n-3)^{-1}(n-5)^{-1}, \nonumber\\
D_{24}^{*} &=& (n^4-10n^2-15)n^{-3}N^{-1}(N-n)(n-3)^{-1}(n-5)^{-1}, \nonumber\\
D_{3^{2}}^{*} &=& (n^4-5n^3+5n^2+5n+10)n^{-3}N^{-1}(N-n)(n-3)^{-1}(n-5)^{-1}, \nonumber\\
D_{2^{3}}^{*} &=& -3(N-n)N^{-1}n^{-2}(n^3-5n^2+5n-5)(n-3)^{-1}(n-5)^{-1}.
\nonumber\\
\end{eqnarray}
Now we consider products of {\bf two} central moments $\mu (\pi_{1}) \mu
( \pi_{2} )$, in the order
\begin{eqnarray*}
( \pi_{1} . \pi_{2} ) = ( 1^{2} . 1^{2} ), ( 1^{2} . 1^{3} ),
(1^{2} . 1^{4} ), ( 1^{3} . 1^{3} ),
 ( 1^{2} . 12 ), ( 1^{2} .  1^{2}
2 ), ( 1^{2} .  2^{2} ) , ( 1^{3} .  12 ), ( 12  . 12 ),
\end{eqnarray*}
\begin{eqnarray}
\mu (1^{2})^{2} &=& C_{1^{2}.2}^{2} \mu_{2}^{2} \mbox{ so } \nonumber \\
D_{4}^{*} &=& 0,\nonumber\\
D_{2^{2}}^{*} &=&(N-n)^2(N^2n^2-3N^2n+3N^2-2Nn^2+3N+3n^2 \nonumber \\
&& -3n)n^{-1}(N-1)^{-1}N^{-3}(n-1)_3^{-1} \nonumber\\
\mu (1^{2}) \mu (1^{3}) &=& C_{1^{2}.2} C_{1^{3}.3} \mu_{2} \mu_{3} \mbox{ so }\nonumber\\
 D_{5}^{*} &=& 0, \nonumber\\
D_{23}^{*} &=&  (N-n)^2(N-2n)(N^2n^2-5Nn^2-2N^2n+2N^2+10N+10n^2 \nonumber \\
&&-10n)n^{-1}(N-1)^{-1}N^{-4}(n-1)_4^{-1},\nonumber \\
\mu (1^{2}) \mu (1^{4}) &=& C_{1^{2}.2} C_{1^{4}.4} \mu_{2} \mu_{4} +
C_{1^{2}.2} C_{1^{4}.2^{2}} \mu_{2}^{3} \nonumber\\
\mbox{ so } D_{6}^{*} &=& D_{3^{2}}^{*}=0,\nonumber\\
 D_{24}^{*} &=&-(N-n)^2(
-63N^4n^3
+190N^4n^2
-225N^4n
+8N^4n^4
+60N^4
\nonumber\\
&&
+60N^3
+640N^3n^2
-525N^3n
+98N^3n^4
-12N^3n^5
-345N^3n^3
\nonumber\\
&&
+198N^2n^4
+3N^2n^6
-45N^2n^5
-750N^2n^3
+1140N^2n^2
-240N^2n
\nonumber\\
&&
-9Nn^5
-990Nn^3
+405Nn^4
+18Nn^6
+360Nn^2
-90n^5
-180n^3
\nonumber\\
&&
+315n^4
-45n^6
)N^{-5}(n-1)_5^{-1}(N-1)^{-1}n^{-3},\nonumber\\
 D_{2^{3}}^{*}&=&  (-60N^2+n^5N^2-11N^2n^4+44N^2n^3-82N^2n^2+90N^2n-60N
\nonumber\\
&& +21Nn^4-96Nn^2+180Nn-24Nn^3-3Nn^5+60n-48n^4\nonumber\\
&& -120n^2+6n^5+126n^3)N^{-2}(N-n)(n-2)_4^{-1}(N-1)^{-1}n^{-2},
\nonumber \\
\mu (1^{3})^{2} &=& C_{1^{3}.3}^{2} \mu_{3}^{2} \mbox{ so } D_{6}^{*} =
D_{24}^{*} = D_{2^{3}}^{*} = 0, \nonumber\\
D_{3^{2}}^{*} &=& (N-2n)^2(N-n)^2(
N^3n^4
+25N^3n^2
-40N^3
-10N^3n
-8N^3n^3
\nonumber\\
&& +8N^2n^3
-26N^2n^2
+70N^2n
-4N^2n^4
-120N^2
+14Nn^4
-4Nn^3
\nonumber\\
&&-70Nn^2
-80N
+220Nn
+80n
-100n^2
+40n^3
\nonumber\\
&&-20n^4
)n^{-3}(N-1)_2^{-1}(n-1)_5^{-1}N^{-5}, \nonumber \\
\mu (1^{2}) \mu (12) &=& C_{1^{2}.2} C_{12.3} \mu_{2} \mu_{3} \mbox{ so }
D_{5}^{*} = 0, \nonumber\\
D_{23}^{*} &=&(N-n)^2(N^2n^2-2N^2n+2N^2+10N-5Nn^2+10n^2\nonumber\\
&&-10n)n^{-1}(N-1)^{-1}
N^{-3}(n-2)_3^{-1}, \nonumber \\
\mu (1^{2}) \mu (1^{2} 2) &=& C_{1^{2}.2} C_{1^{2} 2.4} \mu_{2} \mu_{4} +
C_{1^{2}.2} C_{1^{2}2.2^{2}} \mu_{2}^{3} \nonumber\\
\mbox{ so } D_{6}^{*} &=& D_{3^{2}}^{4} = 0, \nonumber\\
D_{24}^{*} &=&-(N-n)^2(60N^3+2N^3n^4-21N^3n^3+70N^3n^2-105N^3n \nonumber \\
&& -59N^2n^3+60N^2-285N^2n+17N^2n^4+190N^2n^2-n^5N^2 \nonumber \\
&& -135Nn^3-120Nn+3Nn^4
-6n^5N+330Nn^2+60n^2+15n^5\nonumber \\
&&-105n^3 +30n^4)(n-2)_4^{-1}N^{-4}(N-1)^{-1}n^{-3},\nonumber\\
 D_{2^{3}}^{*} &=& (60N^3+N^3n^4-12N^3n^3+50N^3n^2-90N^3n+60N^2+30N^2n^2\nonumber\\
&&-120N^2n+3N^2n^3-3Nn^4-57Nn^3-120Nn+210Nn^2 \nonumber\\
&& -90n^3+60n^2+30n^4)
(N-n)^2(n-2)_4^{-1}N^{-4}(N-1)^{-1}n^{-2},
\nonumber \\
\mu (1^{2}) \mu (2^{2}) &=& C_{1^{2}.2} C_{2^{2}.4} \mu_{2} \mu_{4} + C_{1^{2}.2} C_{2^{2}.2^{2}}\mu_{2}^{3} , \nonumber\\
\mbox{ so } D_{6}^{*} &=& D_{3^{2}}^{*} = 0 , \nonumber\\
D_{24}^{*} &=&-(N-n)^2(-60N+2Nn^5-17Nn^4+49Nn^3-55Nn^2\nonumber\\
&&+45Nn-60+
n^4-145n^2+105n+53n^3\nonumber\\
&& -2n^5)N^{-2}(N-1)^{-1}(n-2)_4^{-1}n^{-3},\nonumber\\
D_{2^{3}}^{*} &=&(-60N^2+N^2n^5-11N^2n^4+44N^2n^3-82N^2n^2+90N^2n-60N\nonumber\\
&&+21Nn^4-96Nn^2+180Nn-24Nn^3-3n^5N+60n-48n^4\nonumber\\
&&-120n^2+6n^5+126n^3)
N^{-2}(N-n)(n-2)_4^{-1}(N-1)^{-1}n^{-2}, \nonumber
\end{eqnarray}
\begin{eqnarray}
\mu (1^{3}) \mu (12) &=& C_{1^{3}.3} C_{12.3} \mu_{3}^{2} \nonumber\\
\mbox{ so } D_{6}^{*} &=& D_{24}^{*} = D_{2^{3}}^{*} = 0, \nonumber\\
D_{3^{2}}^{*} &=&(N-2n)(N-n)^2(
N^3n^4
-40N^3
+25N^3n^2
-8N^3n^3
\nonumber\\
&&-10N^3n
-26N^2n^2
+70N^2n
+8N^2n^3
-120N^2
-4N^2n^4
-70Nn^2
\nonumber\\
&&-80N
+14Nn^4
+220Nn
-4Nn^3
-20n^4
+40n^3
\nonumber\\
&&-100n^2
+80n
)n^{-3}(N-1)_2^{-1}N^{-4}(n-2)_4^{-1},\nonumber
\end{eqnarray}
\begin{eqnarray}
\mu (12)^{2} &=& C_{12.3}^{2} \mu_{3}^{2} \nonumber\\
\mbox{ so } D_{6}^{*} &=& D_{24}^{*} = D_{2^{3}}^{*} = 0 , \nonumber\\
D_{3^{2}}^{*}&=& (n-1)(N-n)^2(
+N^3n^4
-40N^3
+25N^3n^2
-8N^3n^3
-10N^3n
\nonumber\\
&&+8N^2n^3
-26N^2n^2
+70N^2n
-120N^2
-4N^2n^4
-4Nn^3
-70Nn^2
\nonumber\\
&&-80N
+220Nn
+14Nn^4
+40n^3
-100n^2
+80n
\nonumber\\
&&-20n^4
)N^{-3}n^{-3}(N-1_2^{-1}(n-2)_4^{-1}.
\nonumber
\end{eqnarray}
Finally, we have one product of 3 central moments of order 6:
\begin{eqnarray}
\mu (1^{2})^{3} &= & C_{1^{2}.2}^{3} \mu_{2}^{3}\mbox{ has } D_{6}^{*} =
D_{24}^{*} = D_{3^{2}}^{*} = 0, \nonumber \\
D_{2^{3}}^{*} &=&(N-n)^3(N-2)(N^3n^3+9Nn^3-15n^3-3N^2n^3-9N^3n^2 \nonumber \\
&& -9N n^2+9N^2n^2 \nonumber \\
&& +45n^2-30n+29N^3n-6N^2n-45Nn \nonumber \\
&& -30N^3+30N)n^{-1}(N-1)^{-2}(n-1)_5^{-1} N^{-5}.\nonumber
\end{eqnarray}
When $N=\infty$ the nonzero coefficients $D_{\pi^{'}_{-}}^{*}$ for the
UEs of $\mu(\pi_1)\mu(\pi_2) \cdots$ are as follows.
\begin{eqnarray}
\mu (1^{2})^{2} &=& C_{1^{2}.2}^{2} \mu_{2}^{2} \mbox{ so } D_{4}^{*} = 0,\nonumber\\
D_{2^{2}}^{*} &=& n^{-1}(n^2-3n+3)(n-1)_3^{-1}, \nonumber \\
\mu (1^{2}) \mu (1^{3}) &=& C_{1^{2}.2} C_{1^{3}.3} \mu_{2} \mu_{3} \mbox{ so } D_{5}^{*} = 0, \nonumber\\
D_{23}^{*} &=& n^{-1}(n^2-2n+2)(n-1)_4^{-1}, \nonumber \\
\mu (1^{2}) \mu (1^{4}) &=& C_{1^{2}.2} C_{1^{4}.4} \mu_{2} \mu_{4} +
C_{1^{2}.2} C_{1^{4}.2^{2}} \mu_{2}^{3} \nonumber\\
\mbox{ so } D_{6}^{*} &=& D_{3^{2}}^{*}=0, \nonumber\\
D_{24}^{*} &=& -(8n^4-63n^3+190n^2-225n+60)n^{-3}(n-1)_5^{-1},\nonumber\\
 D_{2^{3}}^{*} &=&  (n^4-9n^3+26n^2-30n+30)n^{-2}(n-3)_3^{-1}, \nonumber \\
\mu (1^{3})^{2} &=& C_{1^{3}.3}^{2} \mu_{3}^{2} \mbox{ so } D_{6}^{*} =
D_{24}^{*} = D_{2^{3}}^{*} = 0, \nonumber\\
D_{3^{2}}^{*} &=& n^{-3}(n^4-8n^3+25n^2-10n-40)(n-1)_5^{-1}, \nonumber \\
\mu (1^{2}) \mu (12) &=& C_{1^{2}.2} C_{12.3} \mu_{2} \mu_{3} \mbox{ so }
D_{5}^{*} = 0, \nonumber\\
D_{23}^{*} &=& n^{-1}(n^2-2n+2)(n-2)_3^{-1}, \nonumber \\
\mu (1^{2}) \mu (1^{2} 2) &=& C_{1^{2}.2} C_{1^{2} 2.4} \mu_{2} \mu_{4} +
C_{1^{2}.2} C_{1^{2}2.2^{2}} \mu_{2}^{3} \nonumber\\
\mbox{ so } D_{6}^{*} &=& D_{3^{2}}^{4} = 0, \nonumber\\
D_{24}^{*} &=&-(2n^4-21n^3+70n^2-105n+60)(n-2)_4^{-1}n^{-3},\nonumber\\
D_{2^{3}}^{*} &=&  (n^3-10n^2+30n-30)n^{-2}(n-3)_3^{-1}, \nonumber \\
\mu (1^{2}) \mu (2^{2}) &=& C_{1^{2}.2} C_{2^{2}.4} \mu_{2} \mu_{4} + C_{1^{2}.2} C_{2^{2}.2^{2}}\mu_{2}^{3} , \nonumber\\
\mbox{ so } D_{6}^{*} &=& D_{3^{2}}^{*} = 0 , \nonumber\\
D_{24}^{*} &=& -(2n^5-17n^4+49n^3-55n^2+45n-60)(n-2)_4^{-1}n^{-3} ,\nonumber\\
 D_{2^{3}}^{*} &=& (n^4-9n^3+26n^2-30n+30)n^{-2}(n-3)_3^{-1},\nonumber \\
\mu (1^{3}) \mu (12) &=& C_{1^{3}.3} C_{12.3} \mu_{3}^{2} \nonumber\\
\mbox{ so } D_{6}^{*} &=& D_{24}^{*} = D_{2^{3}}^{*} = 0, \nonumber\\
D_{3^{2}}^{*} &=& n^{-3}(n^4-8n^3+25n^2-10n-40)(n-2)_4^{-1},
\nonumber
\end{eqnarray}
\begin{eqnarray}
\mu (12)^{2} &=& C_{12.3}^{2} \mu_{3}^{2} \nonumber\\
\mbox{ so } D_{6}^{*} &=& D_{24}^{*} = D_{2^{3}}^{*} = 0 , \nonumber\\
D_{3^{2}}^{*} &=& (n-1)n^{-3}(n^4-8n^3+25n^2-10n-40)(n-2)_4^{-1}.\nonumber
\end{eqnarray}
Finally, we have one product of 3 central moments of order 6:
\begin{eqnarray}
\mu (1^{2})^{3} &=& C_{1^{2}.2}^{3} \mu_{2}^{3} \mbox{ has } D_{6}^{*} =
D_{24}^{*} = D_{3^{2}}^{*} = 0,\nonumber \\
 D_{2^{3}}^{*} &=& n^{-1}(n^2-7n+15)(n-3)_3^{-1}(n-1)^{-1}. \nonumber
\end{eqnarray}

\section{UEs of products of cumulants}

Here we give $ E \widehat{\kappa} (\pi_{1}) \widehat{\kappa} (\pi_{2}) \ldots$ and UEs of
$ \kappa (\pi_{1}) \kappa (\pi_{2}) \ldots$ for joint cumulants
$\{ \kappa(\pi_i) \}$ up to total order 6 not covered by Section 2,3,5.

We first give coefficients $ K_{\pi .  \pi \_}$ such that
\begin{eqnarray}
\kappa ( \pi) = \sum_{\pi \_}^{r} K_{\pi .  \pi \_} mu (\pi \_)
\mbox{ for }
\pi \mbox{ of order } r>1.
\label{5.1}
\end{eqnarray}
Its UE is then (5.3) for $ D_{\pi \_}^{*}$ of (5.4) with $D_{\pi \_} =
K_{\pi .  \pi \_}.$
\begin{eqnarray}
\kappa(1^{4}) : \mbox{ Since } \kappa_{1^{4}} &=& \mu_{1^{4}} - 3 \mu_{1^{2}}^{2} ,
\nonumber\\
K_{1^{4}.4} &=& M_{1^{4}.4} =C_{1^{4}.4} , \nonumber\\
K_{1^{4}.2^{2}} &=& M_{1^{4}.2^{2}} - 3 M_{1^{2}.2}^{2}= C_{1^{4}.2^{2}} - 3 C_{1^{2}.2}^{2} \nonumber\\
&=& -3(N-n)(N^3+4Nn^2-4N^2n-2N^2+N\nonumber\\
&& +6Nn-6n^2)n^{-3}(N-1)^{-2}(N-2)_2^{-1},\nonumber
\end{eqnarray}
so
\begin{eqnarray}
D^{*}_4&=&K_{1^4.4}*C^{4.4}+K_{1^4.22}C^{22.4}\nonumber\\
&=& (N-n)(nN^3+N^3-6n^2N^2-nN^2-7N^2+6n^2N+6Nn^3\nonumber\\
&&+12nN-6n^3-6n^2)(n-1)_3^{-1}N^{-3}(N-1)^{-1}n^{-1}\nonumber\\
D_{22}^{*} &=& -3(N-n)(nN^3-N^3-4n^2N^2-nN^2+7N^2+6n^2N\nonumber\\
&&+4Nn^3-12nN
-6n^3+6n^2)(n-1)_3^{-1}N^{-3}(N-1)^{-1}n^{-1}. \nonumber
\end{eqnarray}
$ \kappa (1^{5})$: similarly
\begin{eqnarray}
K_{1^5.5} &=& C_{1^5.5},\\
K_{1^5.23}&=&C_{1^5.23}-10C_{1^2.2}C_{1^3.3}\nonumber \\
&=& -10(N-n)(N-2n) (N^3- 6 n N^2  + 6 n^2  N  - 2 N^2 \nonumber \\
&& + 12 n N - 12 n^2 N) n^{-4}(N-1)^{-2}(N-2)_3,\\
D_{5}^{*} &=&
(-N+2n)(-N+n)(5N^3+nN^3-12n^2N^2-nN^2-65N^2+12n^3N\nonumber \\
&& +12n^2N+120nN-12n^3-60n^2) n^{-1}(n-1)_4^{-1}N^{-4}(N-1)^{-1}, \nonumber\\
D_{23}^{*} &=&
 -10(N-2n)(N-n)(nN^3-N^3-6n^2N^2-nN^2+13N^2+6n^3N\nonumber \\
&& + 12n^2N-24nN-12n^3+12n^2)
 n^{-1}(n-1)_4^{-1}N^{-4}(N-1)^{-1},
\end{eqnarray}
\begin{eqnarray}
\kappa (1^{6}):  K_{1^{6}.6} &=& C_{1^{6}.6} , \nonumber\\
K_{1^{6}.24} &=& C_{1^{6}.24} - 15 C_{1^{4}.4} C_{1^{2}.2} \nonumber\\
&=& -15(N-n)(
+N^5
-16N^4n
+N^4
+30N^3n
+64N^3n^2
\nonumber\\
&& -9N^3
-96N^2n^3
+10N^2n
+11N^2
-150N^2n^2
+48Nn^4
\nonumber\\
&&+240Nn^3
-10Nn^2
-4N
-120n^4
)n^{-5}(N-1)
^{-2}(N-2)_4^{-1},\nonumber\\
K_{1^{6}.3^{2}} &=& C_{1^{6}.3^{2}} - 10 C_{1^{3}.3}^{2}\nonumber\\
&=& -10(N-n)(
N^6
-5N^5
-12N^5n
+48N^4n^2
+13N^4
\nonumber\\
&&+54N^4n
-72N^3n^3
-234N^3n^2
-23N^3
-66N^3n
+36N^2n^4
\nonumber\\
&&+360N^2n^3
+306N^2n^2
+22N^2
-180Nn^4
-8N
-480Nn^3
\nonumber\\
&&+240n^4
)n^{-5}(N-1)^{-2}(N-2)^{-2}(N-3)_3^{-1},\nonumber\\
K_{1^{6}.2^{3}} &=& C_{1^{6}.2^{3}} -15 C_{1^{4}.2^{2}}C_{1^2.2}
 + 30 C_{1^{2}.2}^{3} \nonumber\\
&=& 30(N-n)(
N^6
-5N^5
-11N^5n
+9N^4
+52N^4n
+39N^4
\nonumber\\
&&-7N^3
-56N^3n^3
-176N^3n^2
-71N^3n
+30N^2n
+28N^2n^4
\nonumber\\
&&+248N^2n^3
+2N^2
+191N^2n^2
-240Nn^3
-30Nn^2
-124Nn^4
\nonumber\\
&&+120n^4
n^2
)n^{-5}(N-1)^{-3}(N-2)_4^{-1},\nonumber
\end{eqnarray}
\begin{eqnarray}
D_{6}^{*} &=& (N-n)(N^7n^3-2880Nn^3+3120n^2N^2-1440nN^3+960n^4+248N^4
\nonumber\\
&& +600n^7N-
3000n^6N^2+6150n^5N^3+13800n^5N+6000n^6N\nonumber\\
&&-16860n^5N^2-3840n^6-240n^7
-2640n^5-6300N^4n^4\nonumber\\
&&+18240N^3n^4-20280N^2n^4+8400Nn^4+16050n^3N^3\nonumber\\
&&-
10020n^3N^2-8402n^3N^4+5280n^2N^3-7592n^2N^4\nonumber\\
&&-922nN^4+2945N^5n^3-140
N^5-104N^6+1430N^5n^2\nonumber\\
&&+2005N^5n-454N^6n^2-134N^6n16N^7n^2-
94N^6n^3\nonumber\\
&& + 150N^5n^5-30N^6n^4-4N^7-480n^7N^2+120N^3n^7\nonumber\\
&&-240N^4n^6+960
n^6N^3-600N^4n^5+11N^7n\nonumber\\
&&+210n^4N^5)N^{-5}(n-1)_5^{-1}(N-1)_2^{-1}(N-1)^{-1}n^{-3}\nonumber
\end{eqnarray}
\begin{eqnarray}
D_{24}^{*} &=&  -15(N-n)(N^7n^3+2880Nn^3-3120n^2N^2+1440nN^3-960n^4-248N^4
\nonumber \\
&& +456n^7N
-960n^6N^2+722n^5N^3+120n^5N+912n^6N-1596n^5N^2\nonumber \\
&&-480n^6  -240n^7+1680
n^5-302N^4n^4+1512N^3n^4+3360N^2n^4 \nonumber \\
&& -5520Nn^4 - 6650n^3N^3 +6900n^3N^2- 626n^3N^4-3840n^2N^3\nonumber \\
&& +5236n^2N^4+674nN^4 +115N^5n^3+140N^5+104N^6-44
N^5n^2 \nonumber \\
&& -1745N^5n-34N^6n^2+118N^6n+2N^7n^2 -4N^6n^3+64N^5n^5
-16N^6 n^4 \nonumber \\
&& +4N^7-264n^7N^2+48N^3n^7-96N^4n^6+528n^6N^3-342N^4n^5-7N^7n \nonumber \\
&& +78n^4N^5)n^{-3}(n-1)_5^{-1}(N-2)_2^{-1}(N-1)^{-1}N^{-5},\nonumber
\end{eqnarray}
\begin{eqnarray}
D_{3^{2}}^{*} &=& -10(N-n)(N^7n^3-2880Nn^3+3120n^2N^2-1440nN^3+960n^4
\nonumber \\
&& +248N^4 + 408n^7 N-672n^6N^2+102n^5N^3-600 n^5 N-48n^6N \nonumber \\
&& +1284n^5N^2+480n^6 -240n^7-1200 n^5+354N^4n^4-1704N^3n^4 \nonumber \\
&& -1920N^2n^4+4080Nn^4+5070n^3N^3  -5340n^3N^2+790n^3N^4 \nonumber \\
&& +3120n^2N^3-4496n^2N^4-550nN^4 -241N^5n^3-140N^5 -104N^6-4N
^5n^2 \nonumber \\
&& +1615N^5n+62N^6n^2-110N^6n-2N^7n^2 -4N^6n^3+48N^5n^5
- 12N^6n^4 \nonumber \\
&& -4N^7-216n^7N^2+36N^3n^7-72N^4n^6+432n^6N^3
-282N^4n^5+5N^7n \nonumber \\
&& +66n^4N^5)n^{-3}N^{-5}(n-1)_5^{-1}(N-1)^{-1}_2(N-1)^{-1},\nonumber
\end{eqnarray}
\begin{eqnarray}
D_{2^{3}}^{*} &=& 30(N-n)(N^7n^2+1440Nn^3-1560n^2N^2+720nN^3-480n^4-124N^4
\nonumber \\
&& -180n^6
N^2+360n^5N^3-384n^5N+384n^6N-588n^5N^2-240n^6+720n^5\nonumber \\
&&-239N^4n^4-64
N^3n^4+2292N^2n^4-1680Nn^4-2796n^3N^3+900n^3N^2\nonumber \\
&&+563n^3N^4 +540n^2N^
3+1278n^2N^4-594nN^4+44N^5n^3+130N^5-8N^6 \nonumber \\
&& -322N^5n^2 -167N^5n+11N^6
n^2+44N^6n-11N^6n^3+2N^7+28n^6N^3 \nonumber \\
&& -56N^4n^5 -3N^7n +39n^4N^5)
n^{-2}N^{-5}(n-1)_5^{-1}(N-1)_2^{-1}(N-1)^{-1}.
\nonumber
\end{eqnarray}
$\kappa (1^{4}2) : \mbox{ since } \kappa_{1^42} = \mu_{1^42}-6 \mu_{1^{2}} \mu_{1^{2} 2}  - 4  \mu_{1^{3}} \mu_{12},$
\begin{eqnarray*}
K_{1^{4}2.6} &=& M_{1^{4}2.6},\\
K_{1^{4}2.24} &=& M_{1^{4}2.24} - 6 C_{1^{2}.2} C_{1^{2}2.4}\nonumber \\
&=& -(N-n)N(n-1)(
15N^4
-118N^3n
+15N^3
+240N^2n^2
\nonumber\\
&&-135N^2
+186N^2n
+160Nn
-144Nn^3
-540Nn^2
+165N
-60n
\nonumber\\
&&
-60n^2
+360n^3
-60
)n^{-5}(N-1)^{-2}(N-2
)_4^{-1},\nonumber\\
K_{1^{4}2.3^{2}} &=& M_{1^{4}2.3^{2}} - 4 C_{1^{3}.3} C_{12.3} \nonumber\\
&=& -2(N-n)N(n-1)(
5N^5
-31N^4n
-25N^4
+60N^3n^2
+65N^3
+144N^3n
\nonumber\\
&&-36N^2n^3
-211N^2n
-288N^2n^2
-115N^2
+90Nn
+180Nn^3
+372Nn^2
\nonumber\\
&&+110N
-40n
-240n^3
-40
)n^{-5}(N-1)^{-2}(N-2)^{-2}(N-3)_3^{-1},\nonumber\\
K_{1^{4}2.2^{3}} &=& M_{1^{4}2.2^{3}} - 6 C_{1^{2}.2} M_{1^{2}2.2^{2}}
\nonumber\\
&=&  6(N-n)N(n-1)(
+5N^5
-28N^4n
-25N^4
+45N^3
+131N^3n
+50N^3n^2
\nonumber\\
&&-228N^2n^2
-28N^2n^3
-35N^2
-178N^2n
+75Nn
+262Nn^2
+124Nn^3
\nonumber\\
&&+10N
-60n^2
-120n^3
)n^{-5}(N-1)^{-3}(N-2)_4^{-1},\nonumber
\end{eqnarray*}
%so $D_{6}^{*}$ is given by (4.10), {\bf NO!}\\
\begin{eqnarray}
D_{6}^{*}& =& (N-n)(N^6n^3-1416Nn^3+480n^2N+1244n^2N^2-400nN^2-298nN^3
\nonumber \\
&&+528n^4-
192n^3-92N^4+120N^3+96n^6N^2-144n^5N^3-1176n^5N-120n^6N\nonumber \\
&&+540n^5N^2-
24N^5+768n^5+48n^6+92N^4n^4-1090N^3 n^4+2884N^2n^4 \nonumber \\
&& -2376Nn^4-2390n^3
N^3+2932n^3N^2+907n^3N^4-1934n^2N^3+696n^2N^4 \nonumber \\
&& +689nN^4-58N^5n^3
+16N^6n^2+11N^6n-14N^5n^4-24n^6N^3+36N^4 n^5 \nonumber \\
&& -98N^5n-238N^5n^2-4
N^6)n^{-3}(n-2)_4^{-1}(N-1)_2^{-1}(N-1)^{-1}N^{-4},\nonumber \\
D_{24}^{*} &=&  -(N-n)(15N^6n^3+14040Nn^3-7200n^2N-12660n^2N^2+6000nN^2
\nonumber \\
&& + 2670nN^3-5040 n^4+2880n^3+1380N^4-1800N^3+792n^6N^2 \nonumber \\
&&-1260n^5 N^3  -2376n^5N-1368n^6N+
2244n^5N^2+360N^5+1440n^5 \nonumber \\
&& +720n^6+564N^4n^4-1370N^3n^4+3948N^2n^4+360
Nn^4-2706n^3N^3 \nonumber \\
&&-10980n^3N^2+723n^3N^4 +16610n^2N^3-504n^2N^4-8235n
N^4 \nonumber \\
&& -84N^5n^3+30N^6n^2-105N^6n-118N^5n^4 -144n^6N^3+240N^4
n^5 \nonumber \\
&& +1110N^5n-236N^5n^2+60N^6)n^{-3}(n-2)_4^{-1}(N-1)_2^{-1}(N-1)^{-1}N^{-4}
\nonumber
\end{eqnarray}
\begin{eqnarray}
D_{3^{2}}^{*} &=& -2(N-n)(+5N^6n^3-3480Nn^3+2400n^2N+3220n^2N^2-2000nN^2
\nonumber \\
&& -590nN^3+1200
n^4-960n^3-460N^4+600N^3+216n^6N^2-348n^5N^3 \nonumber \\
&&+168n^5N -408n^6N+492n^5 N^2-120N^5-480n^5+240n^6+163N^4 n^4 \nonumber \\
&&+124N^3n^4-1276N^2n^4 +360Nn^4+
1130n^3N^3+2420n^3N^2-391n^3 N^4 \nonumber \\
&& -4480n^2N^3+57n^2N^4+2395n N^4 -8N^5
n^3-10N^6n^2+25N^6n-31N^5n^4 \nonumber \\
&&-36n^6N^3+60N^4n^5-310 N^5n+
133N^5n^2 \nonumber \\
&& -20N^6)n^{-3}(n-2)_4^{-1}(N-1)_2^{-1}(N-1)^{-1}N^{-4}\nonumber \\
D_{2^{3}}^{*} &=& 6(N-n)(5N^6n^2+1320Nn^3-1200n^2N-180n^2N^2+1000nN^2
\nonumber \\
&& -720nN^3-720n^
4+480n^3+350N^4-300N^3-28n^5N^3-384n^5N+180n^5N^2 \nonumber \\
&& -60N^5+240n^5+50N^
4n^4-298N^3n^4+420N^2n^4+504Nn^4+316n^3N^3 \nonumber \\
&& -2116n^3N^2+116n^3N^4+
1798n^2N^3-496n^2N^4-404nN^4-28N^5n^3-15N^6n \nonumber \\
&& +139N^5n+21N
^5n^2+10N^6)n^{-2}(n-2)_4^{-1}(N-1)_2^{-1}(N-1)^{-1}N^{-4}.
\nonumber
\end{eqnarray}
Similarly we can write down UEs for products of cumulants up to order 6
not covered by Section 2, Section 4.  These are
\begin{eqnarray}
\mbox{ of order } & 4 & : \mu^{2} \mu (1^{2}), \mu \mu (12), m (2) \mu
(1^{2}), \nonumber\\
\mbox{ of order } & 5 & : \mu^{3} \mu (1^{2}), \mu^{2} \mu (12), \mu^{2}
\kappa (1^{3}), \mu m(2) \mu (1^{2}), \nonumber\\
& & \mu \mu (13), \mu \mu (2^{2}), m (2) \mu (12), m (2) \mu (1^{3}),
m (3) \mu (1^{2}), \nonumber\\
\mbox{ of order } & 6 &: \mu^{4} \mu
(1^{2}), \mu^{3} \mu(12), \mu^{3} \mu(1^{3}), \mu^{2} m (2)
\mu(1^{2}), \nonumber\\
&& \mu^{2} \mu (13), \mu^{2} \mu (2^{2}),
\mu^{2} \kappa (1^{2} 2), \mu^{2} \kappa (1^{4}), \nonumber\\
&& \mu m
(2) \mu (12), \mu m (2) \mu (1^{3}), \mu \mu (14), \mu \mu (23),
\nonumber\\
&& m (2) \kappa (1^{4}), m (2) \mu (1^{2}2), m (2) \mu
(2^{2}), m (2) \mu (13),\nonumber\\
&& m (2) \mu (1^{2})^{2}, m
(2)^{2} \mu (1^{2}), m (3) \mu (12). \nonumber
\end{eqnarray}
Similarly we can write down an UE for $m (2) \mu (1^{4}). $

\section{Proof of the inversion principles}

We begin by following a false trail that leads in Section 8 to the solution of the eigenfunction
problem stated at the end of the abstract.\\
A proof for $B_r$ would follow if
\begin{eqnarray}
B_{r} = H \Lambda H^{-1} \mbox{ where } \Lambda = diag ( \lambda_{1,} \lambda_{2,} \ldots ), \mbox{ the evalues of } B, \label{6.1}
\end{eqnarray}
and $\Lambda$ satisfies the inversion principle, that is
\begin{eqnarray}
\lambda_{i} ( N, n) = \lambda_{i} \mbox{ satisfies } \lambda_{i}^{-1} = \lambda_{i} ( n, N). \label{6.2}
\end{eqnarray}
\begin{eqnarray*}
\mbox{ Writing } H = ( u_{1}, u_{2}, \ldots ) \mbox{ and } H^{-1} = \left(
\begin{array}{c}
 v_{1}^{'}  \\  v_{2}^{'} \\  \vdots   \end{array}
\right), \mbox{ (6.1) is equivalent to }
\end{eqnarray*}
\begin{eqnarray}
Bu_{i} = \lambda_{i} u_{i} , v_{i}^{'} B = \lambda_{i} v_{i}^{'} , v_{i}^{'} u_{j} = \delta_{ij} , B = \sum_{i} \lambda_{i} u_{i} v_{i}^{'}. \label{6.3}
\end{eqnarray}
This implies that $ a_{i} (F) = v_{i}^{'} \mu_{(r)} $ satisfies
\begin{eqnarray}
E a_{i} ( \widehat{F}) = \lambda_{i} a_{i} (F). \label{6.4}
\end{eqnarray}
If $v_{i}$ does not depend on $(n, N)$ we call $a_{i} (F)$ an
{\bf r-e function} or r-eigen function).
Since  $m_{r} = E X^{r} = \sum_{j=0}^{\infty}
\left ( \begin{array}{c}
r\\ j
\end{array} \right ) \mu^{r-j} \mu_{j} = v_{1}^{'} m_{(r)}$ say
and $E \overline{X}^{r} = E X^{r},$ for all $r$ there is always an {\bf r-e}
function  with $\lambda=1.$\\

Other examples we have seen are $\mu_{r}$ (with  $\lambda = C_{r.r} $)
for $ r = 2$ or 3.  In each case $\lambda_{i}$ satisfies (7.2). \\

One can show that for $r \leq 3$, one can choose $H$ in (7.1) to be
{\bf constant}, so that the number of efunctions equals the number
of partitions of $r$ say $n_r$ that is, the dimension of $m_{(r)}$.
Unfortunately this is false for $r \geq 4$.  In fact in Section 8 we show that
the {\bf only} $r-$efunctions apart from $m_{r}$ (up to a constant multiple)
are $ E \prod_{i=1}^{s} \left( X^{a_{1}} - E X^{a_{i}} \right)$,
with $\lambda = C_{s.s}$ for $s=2$ or $3$ where $ \sum_{1}^{s} a_{i} = r.$\\

 For $ r=1,\ H= \Lambda = 1$ so $ a_{1} = \mu.$\\

 For $r=2$, we can take
 $ \lambda_{1} = 1, \lambda_{2} = C_{2.2},\\
 H = \left ( \begin{array}{lr} 0 & -1 \\ 1 & 1 \end{array} \right ),
H^{-1} = \left (\begin{array}{rr} 1 & 1\\ -1 & 0 \end{array} \right ),  $ so $ a_{1} = \mu_{2} + \mu^{2} = m_2$ and $a_{2} = - \mu_{2}.$\\

For $r=3$, we can take $\lambda_{1} = C_{3.3,} \lambda_{2} = C_{2.2}, \lambda_{3} = 1,$
\begin{eqnarray*}
 H = \left(
\begin{array}{rrl}
1 & 0 & 0 \\
- 1/2 & -1 & 0 \\
1/2 & 3 & 1 \end{array} \right) \mbox{ and } H^{-1} \left(
 \begin{array}{rrl}
1 & 0 & 0 \\
- 1/2 & -1 & 0 \\
1 & 3 & 1 \end{array} \right ) ,
\end{eqnarray*}
so $ a_{1} = \mu_{3}, a_{2} = - ( \mu_{3} + 2 \mu_{2} \mu ) / 2
= - ( m_{3} - \mu m_{2} ) /2, a_{3} = \mu_{3} + 3 \mu_{2} \mu + \mu^{3} = m_{3}.$

This was derived using the following result.
\begin{theorem} \  Suppose for $i = 1, 2,$ that $C_i = H_{i} \Lambda_{i} H_{i}^{-1}$
is $p_{i} \times p_{i}$ with $ \Lambda_{i} = diag ( \lambda_{i1},\lambda_{i2}, \ldots)$
and $C_{1}, C_{2}$ have no evalues in common.
Then for $D$ any $ p_{2} \times p_{1}$ matrix
\begin{eqnarray} \left( \begin{array}{ll}
C_{1} & 0 \\
D & C_{2} \end{array} \right) = H \left(
\begin{array}{ll}
\Lambda_{1} & 0 \\
0 & \Lambda_{2}
\end{array} \right ) H^{-1} \mbox{ where } \label{6.5}
\end{eqnarray}
\begin{eqnarray*}
H = \left (
\begin{array}{cc}
H_{1} &  0 \\
H_{2}X_{0} & H_{2}
\end{array} \right ) , H^{-1} = \left (
\begin{array}{ccc}
H_{1}^{-1} &  0 \\
-X_{0}H_{1}^{-1} & H_{2}^{-1}
\end{array} \right ) ,
\end{eqnarray*}
\begin{eqnarray}
( X_{0})_{ij} = ( H_{2}^{-1} D H_{1} )_{ij} / ( \lambda_{1j} - \lambda_{2i}).
\label{6.6}
\end{eqnarray}
\end{theorem}
The proof uses the following easily proved lemma.
\begin{lemma}  \  \  Suppose  for $i = 1, 2$ that $C_{i} = H_{i} \Lambda_{i} H_{i}^{-1}$ is $p_{i} \times p_{i}$ and $X$ is $p_{2} \times p_{1}$ satisfying
\begin{eqnarray}
X C_{1} - C_{2} X = D. \label{6.7}
\end{eqnarray}
Then (6.5) holds with $ H = \left ( \begin{array}{cc} H_{1} & 0 \\ X H_{1} &
H_{2}
\end{array} \right ), H^{-1} = \left( \begin{array}{cc} H_{1}^{-1} & 0\\
- H_{2}^{-1}X & H_{2}^{-1} \end{array} \right ).$
\end{lemma}
{\bf PROOF OF THEOREM 7.1.}  Put $X_{0} = H_{2}^{-1} X H_{1}$ so (6.7) becomes
\begin{eqnarray*}
X_{0} \Lambda_{1} - \Lambda_{2} X_{0} = H_{2}^{-1} D H_{1} , \mbox{ which is just } (6.6).
\Box
\end{eqnarray*}
The theorem  can be extended to the case where an evalue of$ \left (
\begin{array}{ll} C_1 & 0 \\  0 & C_{2} \end{array} \right )$ has multiplicity
greater than 1.  If $B_{r}$ has repeated evalue, then in place of (6.1), one
has
\begin{eqnarray*}
B_{r} = H \Lambda H^{-1} \mbox{ where } \Lambda_{ij} = \left \{
\begin{array}{lcl}
\lambda_{i} & \mbox{ for } & j=i \\
1 \mbox{ or } 0 & \mbox{ for } & j=i+1 \\
0 && \mbox{ otherwise }. \end{array} \right.
\end{eqnarray*}

But $ \left ( \begin{array}{cc}
\lambda & 0\\
 1  & \lambda \end{array} \right )^{-1}
= \left ( \begin{array}{cc}
\lambda^{-1} & 0\\
- \lambda^{-2} & \lambda^{-1} \end{array} \right )$
so $\Lambda$ does not satisfy the inversion
principle, so a different proof is needed. \\

Recall now the symmetric functions
\begin{eqnarray*}
[ \pi_{1} \pi_{2} \cdots ]_{n} = \sum_{n}^{'} X_{i_{1}}^{\pi_{1}} X_{i_{2}}^{\pi_{2}} \cdots, \ \ [ \pi_{1} \pi_{2} \cdots ]_{N} = \sum_{N}^{'} x_{i_{1}}^{\pi_{1}} x_{i_{2}}^{\pi_{2}} \cdots
\end{eqnarray*}
\begin{eqnarray*}
\mbox{ where } \sum_{n}^{'} \mbox{ sums over } i_{1}, i_{2}, \cdots
\mbox{ distinct in $ 1, \cdots , n$ }.
\end{eqnarray*}
The standardised function $  \langle  \pi_{1} \cdots \pi_{j}  \rangle _{n} = [ \pi_{1} \cdots
\pi_{j} ]_{n} / (n)_{j}$ satisfies the invariance principle
\begin{eqnarray}
E  \langle  \pi_{1} \cdots \pi_{j}  \rangle _{n} &=&  \langle  \pi_{1} \cdots \pi_{j}  \rangle _{N}.
\nonumber\\
\mbox{ So for every partition $\pi$ of $r$, } && E  \langle  \pi  \rangle _{n} =  \langle  \pi  \rangle _{N}.
\label{6.8}
\end{eqnarray}
\begin{eqnarray}%{rl}
\mbox{ Now }  [1]_{N}^{2} = \sum x_{i} x_{j} = \sum_{i \neq j} x_{i} x_{j} + \sum x_{i}^{2} = [1^{2}]_{N} + [2]_{N} \nonumber \\
\mbox{ so }  [1^{2}]_{N} = [1]_{N}^{2} - [2]_{N}. \nonumber
\end{eqnarray}
Similarly we can write
\begin{eqnarray}
 S(\pi) &=& \sum_{\pi^{'}}^r V^{\pi.\pi^{'}} [\pi]_N \mbox{ and }
 [ \pi ]_{N} = V (\pi)^{'} s_{(r)} = \sum_{\pi^{'}}^{r}  \  V_{\pi. \pi^{'}}
\  s(\pi^{'}) \label{6.9} \\
\mbox{ where }  s_{(r)} &=& \{ s_{1}^{r_{1}} s_{2}^{r_{2}} \cdots : 1.r_{1} +
2.r_{2} + \cdots = r \}=\{ s(\pi):r (\pi)=r \}, \nonumber\\
 s(\pi) &=& s_{\pi_1}s_{\pi_2} \cdots,
s_{i} = [ i ]_{N} = \sum_{j=1}^{N} x_{j}^{i} = N m_{i} \mbox{ and }m_{i} =
E X^{i}. \nonumber
\end{eqnarray}
The constant column vector $V (\pi) = \{ V_{\pi.\pi^{'}} : r ( \pi^{'} ) = r
\}$ is given by the columns of Appendix Table 10 of Stuart and Ord (1987) for
$r \leq 6$.  The rows of this table give $\{V^{\pi.\pi^{'}}\}$.  We can write
(7.9) as
\begin{eqnarray*}
([\pi]_N:r (\pi)=r)=Vs_{(r)},n_r \times 1, \mbox{ where }
V = \{ V (\pi)^{'} \} \mbox{ is } n_{r} \times n_{r}.
\end{eqnarray*}
Now $s_{1}^{r_{1}} s_{2}^{r_{2}} \cdots = N^{q} m_{1}^{r_{1}} m_{2}^{r_{2}} \cdots \mbox{ where } q = r_{1} + r_{2} + \cdots = q (1^{r_{1}} 2^{r_{2}} \cdots ).$
That is, as before $q (\pi)$ is the number of partitions in $\pi.$
\begin{eqnarray}%{rlll}
\mbox{ Also }   \langle  \pi  \rangle _{N} &=& [ \pi]_{N} / (N)_{q(\pi)}, \nonumber\\
\mbox{ and }  s_{(r)} &=& D_{N} m_{(r)} \mbox{ where } D_{N} = diag
(N^{q(\pi)} : r (\pi) = r ), \nonumber\\
 m_{(r)} &=& \{ m_{1}^{r_{1}} m_{2}^{r_{2}} \cdots : 1.r_{1} + 2.r_{2} +
\cdots = r \}. \nonumber\\
\mbox{ Also }  mu_{(r)} &=& G m_{(r)} \mbox{ where } G \mbox{ is a constant } n_{r} \times n_{r}
\mbox{ matrix }. \nonumber\\
\mbox{ Set } \overline{D}_{N} &=& diag \{  (N)_{q (\pi)} : r (
\pi) = r   \}, \ n_r \times n_r.\nonumber\\
\end{eqnarray}
So by (7.9), $ \{  \langle  \pi  \rangle _{N} : r (\pi) = r \} =E_{N}  \ m_{(r)}
=F_Nmu_{(r)}$
\begin{eqnarray*}
\mbox{ where } E_N= \overline{D}_{N}^{-1} \ V D_{N} \mbox{ and }
F_N=E_NG^{-1}
\end{eqnarray*}
and (7.8) can be written
\begin{eqnarray*}
E \ F_n \widehat{mu}_{(r)}=F_N \widehat{mu}_{(r)}, \ \ E  \ E_{n} \widehat{m}_{(r)} =  \  \  E_{N} m_{(r)},\mbox{ so that }
\end{eqnarray*}
\begin{eqnarray*}
B_{r} = E_{n}^{-1}  E_{N}, \ B_{r}^{-1} =E_{N}^{-1} \ E_{n},
C_r=F_n^{-1}F_N \mbox{ and } C_r^{-1}=F_N^{-1}F_n.
\end{eqnarray*}
This proves the univariate inversion principles.
We now prove the multivariate inversion principle.

For $\pi$ an ordered partition of $r, [\pi]_{1\ldots r}$ now depends on the
 order of $\pi_{1} \pi_{2} \ldots $.  For example
\begin{eqnarray*}
[21]_{123}^{N} = \sum^{'} x_{i1} x_{i2} x_{j3} \neq [12]_{123}^{N}
= \sum^{'} x_{i1} x_{j2} x_{j3} = [21]_{231} = [21]_{321}
\end{eqnarray*}
Put another way, for $\pi$ an unordered partition, we need to consider not just $[\pi]_{1 \ldots r}$ but all distinct $[\pi]_{a}$ for $a$ any permutation of $1 \ldots r.$

For $r=3$, this gives $[3]_{123,}  [21]_{123,} [21]_{231,}
[21]_{312,} [1^{3}]_{123,}$ so the dimension of $\mu_{(3) 123} = \{ \mu (\pi)_{123} \}$ is 5 whereas that for $ \mu_{(3)}$ was only 3.

As in (7.9) there exists a constant vector $V(\pi)$ such that\\ $[ \pi]_{1 \ldots r}^{N} = V_{1 \ldots r} (\pi)' S_{1 \ldots r} (\pi), $
that is $ ([ \pi]_{1 \ldots r}^{N} ) = V s_{(r)1 \ldots r}$ say.
Also
$s_{(r)1..r}=D_Nm_{(r)1..r}$ where $D_{N}$ as before but with dimensions $d_{r}$ say. \\
Redefining $\overline{D}_{N} $ similarly.
\begin{eqnarray*}
\mbox{ So }  E \overline{D}_{n}^{-1} V D_{n} \widehat{m}_{(r)1 \ldots r} &=& E \overline{D}_{n}^{-1} V S_{(r) 1 \ldots r} = E \overline{D}_{n}^{-1} ([\pi]_{1 \ldots r}^{n} ) \\
 &=& \overline{D}_{N}^{-1} ([ \pi]_{1 \ldots r}^{N}) = \overline{D}_{N}^{-1} Vs_{(r)1..r}=\overline{D}_{N}^{-1} VD_{N} m_{(r) 1 \ldots r}.
\\
\mbox{ So }  E \widehat{m}_{(r)1 \ldots r} &=& E_{n}^{-1} E_{N} m_{(r) 1 \ldots r}
\mbox{ where } E_{n} = \overline{D}_{n}^{-1} V D_{n} \end{eqnarray*}
and $ E_{N}^{-1} E_{n} \widehat{m}_{(r) 1 \ldots r} \mbox{ is an {\bf UE} of } m_{(r)1..r}.$
This proves the multivariate inversion principle for $m_{(r) 1 \ldots r}$.
That for $\mu_{(r) 1 \ldots r}$ follows since $\mu_{(r) 1 \ldots r} = G m_{(r) 1 \ldots r}$ where $G$ is a constant non singular matrix.\\
Note 7.1  The evalues of $B_{r}$ satisfy
\begin{eqnarray*}
| \lambda V \Lambda_{nN} - \overline{\Lambda}_{nN} V | = 0 \mbox{ where } \Lambda_{nN} = D_{n} D_{N}^{-1}, \overline{\Lambda}_{nN} = \overline{D}_{n} \overline{D}_{N}^{-1}.
\end{eqnarray*}
Since $\Lambda{_nN}$ and $ \overline{\Lambda}_{nN}$ are diagonal and $V$ is
lower-triangular, it follows that the evalues of $B_{r}$ are $\{ \nu_{i} =
e_{i} e_{1}^{-i}, \  1 \leq i \leq r \}$ where $ e_{i} = (n)_{i} / (N)_{i}$,
with $ \nu_{i}$ having multiplicity equal to the number of partitions $ \pi $
of $r$ with $ q (\pi) = i.$

\begin{example}
$ \mu_{(4)}^{'} = (\mu_{4}, \ \mu_{2}^{2}, \ \mu \mu_{3}, \ \mu^{2}
\mu_{2}, \ \mu^{4}) \mbox{ so }$
\begin{eqnarray*}  B_{4} = \left (
\begin{array}{ll}
A_{4} & 0 \\ B_{21} & B_{22} \end{array} \right ) \mbox{ has evalues }
\nu_{1} \ \nu_{2}^2 \ \nu_{3} \  \nu_{4}, \mbox{ that is }
\{ \nu_{1} = 1, \nu_{2}, \nu_{2}, \nu_{3}, \nu_{4} \}.
\end{eqnarray*}
\begin{eqnarray*}
\mbox{But } B_{22} &=& \left( \begin{array}{ccc}
C_{3.3} & 0 & 0 \\
3 C_{12.3} & C_{2.2} & 0 \\
10 C_{1^{3}.3} & 10 C_{1^{2}.2} & 1
\end{array} \right ) \mbox{  has evalues } \nu_{1} = 1, \ \nu_{2} = C_{2.2}, \ \nu_{3} = C_{3.3}. \end{eqnarray*}
So $A_{4}$ has evalues $\nu_{2},  \  \nu_{4} \ $ , as was confirmed using
MAPLE.

To see that $H$ for $B_{4}$ is not constant, it suffices to show that $H$ for $A_{4}$ is not constant.  Set  $ \left ( \begin{array}{cc} a & b \\ d & c
\end{array} \right ) = A_{4}.$
Then $A u = \lambda u \Longleftrightarrow (a - \lambda) u_{1} + b u_{2} = 0, \ d u_{1} + ( c - \lambda) u_{2} = 0 ,$
\begin{eqnarray*}
\Longrightarrow \mbox{ one can take $u$ proportional to } \left ( \begin{array}{c}
b \\ \lambda - a \end{array} \right ) \mbox{ or } \left ( \begin{array}{c}
\lambda - c \\ d \end{array} \right ).
\end{eqnarray*}
So if $H$ is constant for $ A_{4}$ then $( \lambda - a ) /b$ and $ (\lambda - c ) / d$ are constant.\\
   For $ \lambda = \nu_{4}, (\lambda-a)/ b = -1/3$ and $( \lambda - c ) / d = - 3, $\\
but for $ \lambda = \nu_{2}, ( \lambda-a)/b = (Nn -n-1-N) / (2 Nn - 3n +
3 -3N)$\\
and $ (\lambda - c) / d = (2 Nn - 3n+3-3N) / (Nn-n-1-N)$
 depend on $n$ and $N$. \   \  \end{example}$\Box$
\begin{example}  \ $\mu_{(5)}^{'} =
( \mu_{5,} \mu_{2} \mu_{3}, \mu \mu_{4}, \mu_{2}^{2} \mu,
\ \mu^{2} \mu_{3}, \mu^{3} \mu_{2}, \mu^{5})$ \\
so $B_{6}$
has evalues $\nu_{1} \ \nu_{2}^2 \ \nu_{3}^2 \ \nu_{4} \ \nu_{5},$
 that is $\{ \nu_{1} = 1, \ \nu_{2}, \ \nu_{2}, \ \nu_{3}, \ \nu_{3}, \nu_{4},
\nu_{5} \}.$
\begin{eqnarray*}
\mbox{ But by } (2.7) B_{22} = \left[
\begin{array}{ll} A_{4} & 0 \\ C & D \end{array} \right ]
\end{eqnarray*}
where D is lower triangular with diagonal $(\nu_{3} = C_{3.3} , \  \nu_{2} = C_{2.2} , \ \nu_{1}=1 ).$

So $A_{5}$ has evalues $\nu_{3}, \ \nu_{5}.  \Box $ \end{example}
\begin{example}  $ B_{6}$ has evalues $\nu_{1} \nu_{2}^{3} \nu_{3}^{3} \nu_{4}^{2} \nu_{5} \nu_{6}.$
\begin{eqnarray*}
\mbox{ But for } r=6,  B_{22} = \left [
\begin{array}{lllll}
A_{5} & 0 & 0 & 0 & 0 \\
.     & A_{4} & 0 & 0 & 0 \\
. & . & C_{3.3} & 0 & 0 \\
. & . & . & C_{2.2} & 0 \\
. & . & . & . & 1 \end{array} \right ]
\end{eqnarray*}
so has evalues $\nu_{1}  \nu_{2}^{2} \nu_{3}^{2}  \nu_{4}  \nu_{5}.$
So $A_{6}$ has evalues $\nu_{2} \nu_{3} \nu_{4} \nu_{6}.  \  \Box $
\end{example}
The relation between the evalues of $A_{r}$ and the number of $\pi\_$
with $q (\pi\_) = i$ is not known.  $\Box$

\section{Symmetric functions and efunctions}

Here we extend the usual expressions for symmetric functions in terms of
products of power sums, and use these relations to show that the only
efunctions are the noncentral moments and the 2nd and 3rd order generalised central
moments.

The symmetric function $ [ \pi_{1} \pi_{2} \cdots ] = \sum_{N}^{'} \
x_{i_{1}}^{\pi_{1}} \  x_{i_{2}}^{\pi_{2}} \cdots$ is well defined for
$\pi_{1}, \  \pi_{2}, \ \ldots$ any real or complex numbers.

The constants $V ( \pi)$ in (7.9) tabled in Stuart and Ord can be derived as follows.  We use an obvious notation.
\begin{eqnarray}
{[ a] [b]} &=& [a b] + [a+b] \mbox{ so } [ab] = [a][b] - [a+b], \nonumber\\
{[a][b][c]} &=& [abc] + \sum^{3} [a, b+c] + [a+b+c] \nonumber\\
\mbox{ so }[abc] &=& [a][b][c] - \sum^{3} [a] [b+c] + 2 [a+b+c], \nonumber\\
{[a] \cdot [d]} &=&[a \cdots d] + \sum^{6} [a,b,c+d] + \sum^{3} [a+b,c+d]
\nonumber\\
&& + \sum^{4}[a,b+c+d]+[a+b+c+d], \nonumber\\
\mbox{ so } [a \ldots d]&=&[a] \ldots [d] - \sum^{6} [a][b][c+d] + \sum^{3}
[a+b][c+d]  \nonumber\\
 && + 2\sum^{4} [a][b+c+d] -6 [a+b+c+d], \nonumber \\
{[a] \cdots [e]} &=& [a...e] + \sum^{10} [abc,d+e] + \sum^{15} [a,b+c,d+e]
\nonumber\\
&&+ \sum^{10} [ab,c+d+e] + \sum^{10} [a+b,c+d+e]  \nonumber\\
 && + \sum^{5} [a,b+c+d+e]+ [a + .. + e],  \nonumber\\
so [a...e] &=& [a]...[e] - \sum^{10} [a][b][c][d+e] + \sum^{15} [a][b+c][d+e]
\nonumber\\
&&  + 2 \sum^{10} [a][b][c+d+e]-2 \sum^{10} [a+b][c+d+e]\nonumber\\
&&  - 6 \sum^{5} [a][b+c+d+e]+  24 [a + .. +e].\nonumber
\end{eqnarray}
Setting  $ S_{m} (1^{r_{1}} 2^{r_{2}} \cdots) = \sum^{m} [a_{1} \ldots
a_{r_{1}}, \
a_{r_{1}+1}+a_{r_{1}+2}, \ldots , a_{r_{1} +2 r_{2}-1} + a_{r_{1}+2 r_{2}} $,
\begin{eqnarray}
\left. a_{r_{1}+2 r_{2}+1} + a_{r_{1} +2 r_{2} +2} + a_{r_{1}+2 r_{2} + 3},
\ldots \right] \label{7.1}
\end{eqnarray}
and $T_{m} (1^{r_{1}} 2^{r_{2}} \ldots )= \sum^{m} [a_{1}] \ldots [a_{r_{1}}]
[a_{r_{1}+1} + a_{r_{1}+2}] \ldots [a_{r_{1}+2r_{2}-1} + a_{r_{1}+2 r_{2}}]$
\begin{eqnarray}
[ a_{r_{1}+2r_{2}+1} \ + \ a_{r_{1}+2r_{2}+2} \ +  \ a_{r_{1}+2r_{2} +3}]
\ldots , \label{7.2}
\end{eqnarray}
where $m = P(1^{r_1}2^{r_2}...)$,
we can write these pairs of equations more compactly as
\begin{eqnarray}
T_{1} (1^{2}) &=&  S_{1} (1^{2}) + S_{1} (2), \label{7.3} \\
S_{1} (1^{2}) &=& T_{1} (1^{2}) - T_{1} (2) , \label{7.4} \\
T_{1} (1^{3}) &=& S_{1} (1^{3}) + S_{3} (12) + S_{1} (3), \label{7.5} \\
S_{1} (1^{3}) &=& T_{1} (1^{3}) - T_{3} (12) + 2 T_{1} (3), \label{7.6} \\
T_{1} (1^{4}) &=& S_{1} (1^{4}) + S_{6} (1^{2}2)+S_{3} (2^{2}) + S_{4} (13) +
S_{1} (4), \label{7.7}\\
S_{1} (1^{4}) &=& T_{1} (1^{4}) - T_{6} (1^{2}2) + T_{3} (2^{2}) + 2T_{4}
(13) - 6 T_{1} (4), \label{7.8} \\
T_{1}(1^{5}) &=& S_{1} (1^{5})+S_{10} (1^{3}2) + S_{15} (12^{2}) + S_{10}
(1^{2}3) \nonumber\\
&& + S_{10} (23) + S_{5} (14) + S_{1} (5), \label{7.9} \\
S_{1} (1^{5})&=&T_{1} (1^{5}) - T_{10} (1^{3}2) + T_{15} (12^{2}) + 2 T_{10}
(1^{2} 3) \nonumber\\
&& - 2T_{10} (23) - 6 T_{5} (14) + 24 T_{1} (5). \label{7.10}
\end{eqnarray}
In this way we replace and extend the last row and column of the
$ r^{th} $ table in Appendix Table 10 of Stuart and Ord by a pair of equations.

The general expression for $T_1(1^r)$ is \\
\begin{eqnarray*}
T_{1} (1^{r}) = \sum_{i=1}^{r}  \  B_{ri}
\end{eqnarray*}
where $B_{ri}$ can be written down from the expression for the partial
exponential Bell polynomial $B_{ri} (x)$ tabled on p307 of Comtet (1974).
For example

$B_{64} (x) = 20 x_{1}^{3} x_{3} + 4 5 x_{1}^{2} x_{2}^{2}$ gives $ B_{64} =
S_{20} (1^{3} 3) + S_{45} (1^{2} 2^{2}). $

In this way one obtains
\begin{eqnarray}
T_{1} (1^{6}) &=& S_{1} (1^6) + S_{15} (1^4 2) + S_{20} (1^3 3) + S_{45} (1^2 2^2) + S_{15} (1^2 4) \nonumber\\
&& + S_{60} (123) + S_{15} (2^3) + S_{15} (24) + S_6 (15) + S_{10} (3^2) + S_1 (6).
\label{7.11}
\end{eqnarray}
The coefficients of the reverse series are most easily obtained from
the coefficients of the terms in the expansion of $( 1^{r})$ in
Appendix Table 10.  In this way we obtain
\begin{eqnarray}
S_{1}(1^6) &=& T_1 (1^6) - T_{15} (1^4 2) + 2 T_{20} (1^3 3) + T_{45} (1^2
2^2) -6 T_{15} (1^2 4) \nonumber\\
&&- 2T_{60} (123) - T_{15} (2^3) + 6 T_{15} (24) + 24 T_{6} (15) + 4 T_{10}
(3^2) - 120 T_1 (6). \label{7.12} \nonumber
\end{eqnarray}
(More generally reading their $rth$ table horizontally gives $\{ S_m(\pi)/m
\}$ in terms of $ \{ T_{m'}(\pi')/m' \},$ while reading vertically gives the
reverse.  For example
reading line 2 of their 3rd table gives $S_3(12)/3=1.T_1(3)/1+1.T_3(12)/3$,
that is $\sum^3[a,b+c]=3[a+b+c]+\sum^3[a][b+c]$.)\\
These expressions for $\{ S_{1} (1^i) \}$ can be used to prove our assertion
in Section 6 that the only {\bf $r-e$} functions are  $ m_r , ($ with $\lambda= 1),
\mu_{ab}= m_{a+b} - m_a m_b $ where $a+b = r ,$ with $\lambda = C_{2.2}, $ and
$ \mu_{abc} = m_{a+b+c} - \sum^{3} m_{a} m_{b+c} + 2 m_{a} m_b m_c $ where $
a+b+c = r, $ (with $\lambda = C_{3.3}).$
\begin{eqnarray}%{rlcl}
\mbox{ We have }   \langle  a+b \rangle _{N} &=& [a+b]_{N} / N = m_{a+b}, \nonumber\\
\mbox{ and }   \langle ab \rangle _{N} &=& [ab]_{N} / (N)_{2} = ([a]_{N} [b]_{N} - [a+b]_{N}
) / (N)_{2} \nonumber\\
 &=& (N m_{a} m_{b} - N m_{a+b}) / (N-1). \nonumber
\end{eqnarray}
So for $  \langle ab \rangle _{N} + \nu  \langle a+b._{N}$ to be proportional to an efunction for some constant $\nu$, it must equal $N(m_{a} m_{b} + \theta m_{a+b}) / (N-1)$ for
some constant $\theta$, in which case since $E  \langle  \pi  \rangle _{n} =  \langle  \pi  \rangle _{N},$
\begin{eqnarray*}
E \alpha ( \widehat{F} ) = \lambda \alpha (F) \mbox{ where } \alpha (F) = m_a  m_b + \theta m_{a+b} \mbox{ and } \lambda = (1 -n^{-1}) / (1 - N^{-1}) = C_{2.2}.
\end{eqnarray*}
So we need to solve $ (N m_{a} m_{b} - N m_{a+b}) / (N-1) = N (m_{a} m_{b} + \theta m_{a+b}) / (N-1).$

This gives $\theta = \nu = -1$, which proves that (to within a multiplicative
constant), $ \mu_{ab} = m_{a+b} - m_{a} m_{b} $ is the only efunction of this
type.

Taking $a+b=r$ makes it an $r-e$ function.

A similar argument shows that for $ \langle abc \rangle _{N} + \nu_{1} \sum^{3}  \langle a,b+c \rangle _{N} \\
+ \nu_{2}  \langle  a+b+c \rangle _{N}$ to be proportional to an efunction, the efunction
must be $\mu_{abc}$ and that $\lambda = C_{3.3}.$

However when one looks for a linear combination of say $ \langle abcd \rangle _{N}, \ldots , \\
 \langle a+b+c+d \rangle _{N}$ proportional to an efunction, we find that the only solutions have
efunctions of the form $ m_{a}^{'}$ or $ \mu_{a^{'}b^{'}} $ or $ \mu_{a^{'}b^{'}c^{'}}.$

We now give a much more general meaning to the symmetric function
relations (8.3) - (8.11).  The functions $S (\pi) = S_{m(\pi)} (\pi) =
S_{m} (\pi)$ and $T (\pi) = T_{m(\pi)} (\pi) = T_{m} (\pi) $ are
defined in terms of the functions $[a_{1} \ldots a_{r} ]$ and $
[a_{1}] \ldots [a_{r}]$ respectively, where $r = r (\pi).$

Let $t (x_1 \ldots x_r)=t (x_{1}, \ldots , x_{r})$ be an arbitrary function.
Replace $ [a_{1}
\ldots a_{r}]$ in the definition of $S (\pi)$ by $\sum^{'} t (x_{i_{1}}, \ldots,
x_{i_{r}})$ and $[a_{1}] \ldots [a_{r}]$ in the definition of $ T (\pi)$ by
$\sum t (x_{i_{1}}, \ldots , x_{i_{r}}), $ where $ \sum$ sums over
$i_{1}, \ldots ,i_{r}$ in $1, \ldots, N$ and $\sum^{'}$ sums over distinct
$ i_{1}, \ldots ,i_{r}$ in $1, \ldots, N.$

Then {\bf the relations (8.3) - (8.11) remain true.}

For example (8.5) gives
\begin{eqnarray}%{llcl}
 \sum t (x_{i} x_{j} x_{k}) &=& T (1^3) = S (1^3) + S (12) + S(3)
\nonumber\\
\mbox{ where }  S (1^3) &=& \sum^{'} t ( x_{i} x_{j} x_k), \nonumber\\
 S(12) &=& \sum^{3} \sum^{'} t (x_{i} x_{j} x_{j} ) = \sum^{'} \{ t (x_{i}
x_{j} x_{j}) + t (x_{j} x_{i} x_{j}) + t ( x_{j} x_{j} x_{i}) \} \nonumber\\
\mbox{ and } S(3) &=& \sum t (x_{i} x_{i} x_{i} ).
\nonumber
\end{eqnarray}
Similarly (8.6) gives
\begin{eqnarray*}
\sum^{'} t (x_{i} x_{j} x_{k}) &=& S(1^3) = T (1^3) - T (12) + 2 T (3) \\
\mbox{ where } T (1^3) &=& \sum t (x_{i} x_{j} x_{k} ), \\
T (12) &=& \sum^3 \sum t ( x_{i} x_{j} x_{j} ), \\
\mbox{ and } T (3) &=& \sum t ( x_{i} x_{i} x_{i}) = S (3).
\end{eqnarray*}
If we choose $t ( x_{1} x_{2} \ldots ) = x_{1}^{a_{1}} x_{2}^{a_{2}} \ldots,$
we obtain our original (8.3) - (8.11).  However the population need no longer
be real numbers : $x_{1,} \ldots , x_{N}$ can now belong to any space $\Omega$,
and $t ( x_{1}, \ldots , x_{r})$ is any real ( or complex ) function on
$ \Omega^{r}$.  The preceding results on r-efunctions may be similarly
extended.
Let us say that $a (F)$ is an {\bf efunction} with {\bf evalue} $ \lambda_{n N}$ if
\begin{eqnarray*}
E a (\widehat{F}) = \lambda_{nN} a (F).
\end{eqnarray*}
Then for any functions $t_{r} (x_{1,} \ldots , x_{r})$ set
\begin{eqnarray}
 a_{1} (F) &=& E t_{1} (X_{1}^{'}), \label{7.13} \\
 a_{2} (F) &=& E \{ t_{2} (X_{1}^{'} , X_{1}^{'}) - t_{2} ( X_{1}^{'},
X_{2}^{'} ) \}, \label{7.14}
\end{eqnarray}
\begin{eqnarray}
\mbox{ and } a_{3} (F) = E \{ t_{3} (X_{1}^{'} X_{1}^{'} X_{1}^{'}) - \sum^{3} t_{3} (X_{1}^{'} X_{2}^{'} X_{2}^{'}) + 2 t_{3} (X_{1}^{'} X_{2}^{'} X_{3}^{'}) \}
\label{7.15}
\end{eqnarray}
where $ X_{1}^{'}, X_{2}^{'}, X_{3}^{'}, \ldots$ are {\bf independent }with
distribution $F$.

Then for $1 \leq i \leq 3,  \  a_{i} (F)$ is an efunction with $\lambda_{nN} = C_{i.i}.$\\
We call $a_2$ and $a_3$ {\bf generalised 2nd and 3rd order central moments.}
\begin{eqnarray*}
\mbox{For if }  t_{r} (x_{1} \ldots x_{r}) = f_1(x_1) \ldots f_{r} (x_{r}) \mbox{ and } Y_{i} = f_{i} (X)
\end{eqnarray*}
$\mbox{ then }  a_{2} = \mu (Y_{1}, Y_{2}) = \mbox{ covar } (Y_{1}, Y_{2} )
\mbox{ and }  a_{3} = \mu (Y_{1}, Y_{2}, Y_{3} ). $\\
By the same argument as before there are no other efunctions of the form
\begin{eqnarray*}
\sum_{r=1}^{s} \ \sum_{i_{1} \ldots i_r=1} \  \nu_{i_{1} \ldots i_{r}}  \
E t_{r} (X_{i_{1}}^{'} , \ldots , X_{i_{r}}^{'} ).
\end{eqnarray*}
But any smooth functional $T(F)$ can be expanded about a fixed distribution
$F_{0} $ as
\begin{eqnarray*}
T (F) = T (F_{0} ) + \sum_{r=1}^{\infty} E \ T_{F_{0}} (X_{1}^{'}, \ldots, X_{r}^{'}) / r !
\end{eqnarray*}
where $T_{F_{0}} (x_{1}, \ldots , x_{r} ) $ is the $r$th functional (von Mises)
derivative of $T (F_{0}).$

So the only efunctions that exist are those of form (8.12) - (8.14).\\

\section*{Appendix A}

Here, we give the coefficients $C_{{\bm \pi}, {\bm \pi}'}$ obtainable from Sukhatme (1944)
needed for ${\bf C}_{--r}$ of (3.4) and ${\bf C}_r$ of (1.7).
The factorisations were obtained using MAPLE.
Recall that $C_{{\bm \pi}_{-}, {\bm \pi}_{+}'} = 0$ so expressions like $C_{2, 1^2}=0$ are not included here.

We first give $\lambda \big({\bm \pi}\big)$ for ${\bm \pi}$ a partition of $r \leq 6$:
$\lambda\big(1^i\big)=e_i,
\
\lambda\big(1^i2\big)=e_{i+1}-2e_{i+2},
\\
\lambda\big(1^i3\big)=e_{i+1}-3e_{i+2}+2e_{i+3},
\
\lambda\big(1^i4\big)=e_{i+1}-7e_{i+2}+12e_{i+3}-6e_{i+4},
\\
\lambda\big(1^i5\big)=e_{i+1}-15e_{i+2}+50e_{i+3}-60e_{i+4}+24e_{i+5},
\\
\lambda\big(6\big)=e_1-31e_2+180e_3-390e_4+360e_5-120e_6,
\\
\lambda\big(2^2\big)=e_2-2e_3+e_4,
\\
\lambda\big(23\big)=e_2-4e_3+5e_4-2e_5,
\\
\lambda\big(24\big)=e_2-8e_3+19e_4-18e_5+6e_6,
\\
\lambda\big(3^2\big)=e_2-6e_3+13e_4-12e_5+4e_6,
\\
\lambda\big(2^3\big)=e_3-3e_4+3e_5-e_6$.

Set $f_i=\lambda\big(i\big)n^{-i},
\
F_i=Nf_i$, so
$F_1=1,
\\
F_2=\big(N-n\big)n^{-1}\big(N-1\big)^{-1},
\\
F_3=\big(N-2n\big)\big(N-n\big)n^{-2}\big(N-1\big)_2^{-1},
\\
F_4=\big(N-n\big)\big(N^2-6Nn+N+6n^2\big)n^{-3}\big(N-1\big)_3^{-1},
\\
F_5=\big(N-2n\big)\big(N-n\big)\big(N^2-12Nn+5N+12n^2\big)n^{-4}\big(N-1\big)_4^{-1},
\\
F_6 = \big(N-n\big)\big(N^4-30N^3n+16N^3-240Nn^3+150N^2n^2+90Nn^2-90N^2n
\\
\quad
-4N+11N^2+120n^4\big)n^{-5}\big(N-1\big)_5^{-1},
\\
g_1=\lambda\big(12\big)n^{-3}=n^{-2}\big(n-1\big)\big(N-n\big)N^{-1}\big(N-1\big)_2^{-1},
\\
g_2=\lambda\big(2^2\big)n^{-4}=\big(n-1\big)\big(N-n\big)_2n^{-3}N^{-1}\big(N-1\big)_3^{-1},
\\
g_3=\big(\lambda\big(13\big)+3\lambda\big(2^2\big)\big)n^{-4}=\big(n-1\big)\big(4N-5n-2\big)\big(N-n\big)n^{-3}N^{-1}\big(N-1\big)_3^{-1},
\\
g_4=\lambda\big(23\big)n^{-5}=\big(n-1\big)\big(N-n\big)_2\big(N-2n\big)n^{-4}N^{-1}\big(N-1\big)_4^{-1},
\\
g_5=\lambda\big(1^2\big)n^{-2}=n^{-1}\big(n-1\big)N^{-1}\big(N-1\big)^{-1},
\\
g_6=\big(\lambda\big(14\big)+6\lambda\big(23\big)\big)n^{-5},
\\
\quad
=\big(n-1\big)\big(N-n\big)\big(7N^2-24Nn+6n-N+18n^2\big)n^{-4}N^{-1}\big(N-1\big)_4^{-1},
\\
g_7=\lambda\big(24\big)n^{-6},
\\
\quad
=\big(n-1\big)\big(N-n\big)_2\big(N^2-6Nn+3N-4+6n^2\big)n^{-5}N^{-1}\big(N-1\big)_5^{-1},
\\
g_8=\lambda\big(13\big)n^{-4}=\big(n-1\big)\big(N-2n+1\big)\big(N-n\big)n^{-3}N^{-1}\big(N-1\big)_3^{-1},
\\
g_9=\lambda\big(3^2\big)n^{-6},
\\
\quad
=\big(n-1\big)\big(N-n\big)_2\big(N^2-4Nn-N+4+4n^2\big)n^{-5}N^{-1}\big(N-1\big)_5^{-1},
\\
g_{10}=\big(2\lambda\big(13\big)+3\lambda\big(2^2\big)\big)= \big(n-1\big)\big(5N-7n-1\big)\big(N-n\big)n^{-3}N^{-1}\big(N-1\big)_3^{-1},
\\
g_{11}=\big(6\lambda\big(13\big)+7\lambda\big(2^2\big)\big)n^{-4}
\\
\quad
=\big(n-1\big)\big(13N-19n-1\big)\big(N-n\big)n^{-3}N^{-1}\big(N-1\big)_3^{-1},
\\
g_{12}=\big(\lambda\big(14\big)+6\lambda\big(23\big)\big)n^{-5},
\\
\quad
=\big(n-1\big)\big(N-n\big)\big(7N^2+18n^2-24Nn+6n-N\big)n^{-4}N^{-1}\big(N-1\big)_4^{-1},
\\
g_{13}=\big(2\lambda\big(13\big)+9\lambda\big(2^2\big)\big)n^{-4}
\\
\quad
=\big(n-1\big)\big(11N-13n-7\big)\big(N-n\big)n^{-3}N^{-1}\big(N-1\big)_3^{-1},
\\
g_{14}=\big(2\lambda\big(13\big)+\lambda\big(2^2\big)\big)n^{-4},
\\
\quad
=\big(n-1\big)\big(3N-5n+1\big)\big(N-n\big)n^{-3}N^{-1}\big(N-1\big)_3^{-1},
\\
g_{15}=\big(\lambda\big(13\big)+\lambda\big(2^2\big)\big)n^{-4}=\big(n-1\big)\big(2N-3n\big)\big(N-n\big)n^{-3}N^{-1}\big(N-1\big)_3^{-1},
\\
g_{16}=\big(\lambda\big(13\big)+6\lambda\big(2^2\big)\big)n^{-4}=\big(n-1\big)\big(7N-8n-5\big)\big(N-n\big)n^{-3}N^{-1}\big(N-1\big)_3^{-1},
\\
g_{17}=\big(\lambda\big(14\big)+6\lambda\big(23\big)\big)n^{-5}
\\
\quad
=\big(n-1\big)\big(N-n\big)\big(7N^2+18n^2-24nN+6n-N\big)n^{-4}N^{-1}\big(N-1\big)_4^{-1},
\\
g_{18}=\lambda\big(3^2\big)n^{-6}
\\
\quad
=\big(n-1\big)\big(N-n\big)_2\big(N^2+4n^2-4nN-N+4\big)n^{-5}N^{-1}\big(N-1\big)_5^{-1},
\\
h_1=\lambda\big(12^2\big)n^{-5}=\big(n-1\big)_2\big(N-n\big)_2n^{-4}N^{-1}\big(N-1\big)_4^{-1},
\\
h_2=\lambda\big(2^3\big)n^{-6}=\big(n-1\big)_2\big(N-n\big)_3n^{-5}N^{-1}\big(N-1\big)_5^{-1},
\\
h_3=\lambda\big(1^22\big)n^{-4}=n^{-3}\big(n-1\big)_2\big(N-n\big)N^{-1}\big(N-1\big)_3^{-1},
\\
h_4=\lambda\big(1^3\big)n^{-3}=n^{-2}\big(n-1\big)_2N^{-1}\big(N-1\big)_2^{-1}$,
\\
and $F_r=Nf_r$, $G_r=N^{2}g_r$, $H_r=N^{3}h_r$.
Then,
\\
$C_{2, 2}=F_1-F_2 = {N\big(n-1\big)n^{-1}\big(N-1\big)^{-1}},
\\
C_{1^2, 2}=F_2 ={\big(N-n\big)n^{-1} \big(N-1\big)^{-1}},
\\
C_{1^2, 1^2}=1,
\\
C_{3, 3}=F_1-3F_2+2F_3 ={N^{2} \big(n-1\big)_2n^{-2} \big(N-1\big)_2^{-1}},
\\
C_{12, 3}=F_2-F_3 ={N \big(n-1\big) \big(N-n\big)n^{-2} \big(N-1\big)_2^{-1}},
\\
C_{12, 12}=C_{2, 2},
\\
C_{12, 1^3}=0,
\\
C_{1^3, 3}=F_3 ={ \big(N-2n\big) \big(N-n\big)n^{-2} \big(N-1\big)_2^{-1}},
\\
C_{1^3, 12}=3C_{1^2, 2},
\\
C_{1^3, 1^3}=1,
\\
C_{4, 4}=F_1-4F_2+6F_3-3F_4
\\
\quad
=N \big(n-1\big) \big(n^{2}N^{2}-2n^{2}N+3n^{2}-3n-3nN^{2}+3N^{2}+3N\big)n^{-3} \big(N-1\big)_3^{-1},
\\
C_{4, 2^2}=6G_1-9G_2={3N \big(n-1\big) \big(2nN-3n+3-3N\big) \big(N-n\big)n^{-3} \big(N-1\big)_3^{-1}},
\\
C_{2^2, 4}=F_2-2F_3+F_4={N\big(n-1\big) \big(nN-n-1-N\big) \big(N-n\big)n^{-3} \big(N-1\big)_3^{-1}},
\\
C_{2^2, 2^2}=G_5-2G_1+3G_2
\\
\quad
=N \big(n-1\big) \big(n^{2}N^{2}-3n^{2}N+3n^{2}+3n-2nN^{2}+3N^{2}-3N\big)n^{-3} \big(N-1\big)_3^{-1},
\\
C_{13, 4}=F_2-3F_3+2F_4={\big(n-1\big)_2\big(N+1\big)_2 \big(N-n\big)n^{-3} \big(N-1\big)_3^{-1}},
\\
C_{13, 2^2}=-3G_1+6G_2=-{ 3\big(n-1\big)_2 \big(N-n\big)Nn^{-3} \big(N-2\big)^{-1}_2},
\\
C_{13, 13}=C_{3, 3},
\\
C_{13, 1^22}=C_{13, 1^4}=0,
\\
C_{1^22, 4}=F_3-F_4={- \big(n-1\big) N\big(N-2n+1\big) \big(N-n\big)n^{-3} \big(N-1\big)_3^{-1}},
\\
C_{1^22, 2^2}=G_1-3G_2={ \big(n-1\big)N \big(nN-3N+3\big) \big(N-n\big)n^{-3} \big(N-1\big)_3^{-1}},
\\
C_{1^22, 13}=2C_{12, 3},
\\
C_{1^22, 1^22}=C_{2, 2},
\\
C_{1^22, 1^4}=0,
\\
C_{1^4, 4}=F_4={ \big(N-n\big) \big(N^{2}+6n^{2}-6nN+N\big)n^{-3} \big(N-1\big)_3^{-1}},
\\
C_{1^4, 2^2}=3G_2=3{N \big(n-1\big)\big(N-n\big)_2n^{-3} \big(N-1\big)_3^{-1}},
\\
C_{1^4, 13}=4C_{1^3, 3},
\\
C_{1^4, 1^22}=6C_{1^2, 2},
\\
C_{1^4, 1^4}=1,
\\
C_{5, 5}=F_1-5F_2+10F_3-10F_4+4F_5
\\
\quad
=N^{2} \big(n-1\big)_2 \big(n^{2}N^{2}-5n^{2}N+10n^{2}-10n-2nN^{2}+2N^{2}
\\
\quad
+10N\big)n^{-4} \big(N-1\big)_4^{-1},
\\
C_{5, 23}=10G_1-10G_3+40G_4
\\
\quad
=10N^{2} \big(n-1\big)_2 \big(nN-2n+2-2N\big) \big(N-n\big)n^{-4} \big(N-1\big)_4^{-1},
\\
C_{23, 5}=F_2-4F_3+5F_4-2F_5
\\
\quad ={N^{2} \big(n-1\big)_2 \big(nN-n-5-N\big) \big(N-n\big)n^{-4} \big(N-1\big)_4^{-1}},
\\
C_{23, 23}=G_5-7G_1+5G_3-20G_4
\\
\quad
=N^{2} \big(n-1\big)_2 \big(n^{2}N^{2}-2n^{2}N+2n^{2}+10n-5nN^{2}-10N
\\
\quad
+10N^{2}\big)n^{-4}\big(N-1\big)^{-1}_4,
\\
C_{14, 5}=F_2-4F_3+6F_4-3F_5 \\
\quad
=N \big(n-1\big) \big(N-n\big) \big(n^{2}N^{2}-n^{2}N+6n^{2}-3nN^{2}-9nN-6n+15N
\\
\quad
+3N^{2}\big)n^{-4}\big(N-1\big)_4^{-1},
\\
C_{14, 23}=-4G_1+6G_3-30G_4
\\
\quad
=-2N \big(n-1\big) \big(N-n\big) \big(2n^{2}N^{2}+n^{2}N-6n^{2}+6n-12nN^{2}+9nN
\\
\quad
+15N^{2}-15N\big)n^{-4}\big(N-1\big)_4^{-1},
\\
C_{14, 14}=C_{4, 4},
\\
C_{14, 12^2}=C_{4, 2^2},
\\
C_{14, 1^23}=C_{14, 1^32}=C_{14, 1^5}= 0,
\\
C_{12^2, 5}=F_3-2F_4+F_5
\\
\quad
=-N \big(n-1\big) \big(N-n\big) \big(-N^{2}n+2n^{2}N-2n^{2}-2n-3nN+5N
\\
\quad
+N^{2}\big)n^{-4} \big(N-1\big)_4^{-1},
\\
C_{12^2, 23}=2G_1-2G_3+10G_4
\\
\quad
=2N \big(n-1\big) \big(N-n\big) \big(N^{2}n^{2}-2Nn^{2}+2n^{2}+3Nn-4N^{2}n+2n-5N
\\
\quad
+5N^{2}\big) n^{-4}\big(N-1\big)_4^{-1},
\\
C_{12^2, 14}=C_{2^2, 4},
\\
C_{12^2, 12^2}=C_{2^2, 2^2},
\\
C_{12^2, 1^23}=C_{12^2, 1^32}=C_{12^2, 1^5}=0,
\\
C_{1^23, 5}=F_3-3F_4+2F_5=-{N \big(n-1\big)_2 \big(-N^{2}+2Nn+4n-5N\big) \big(N-n\big)n^{-4} \big(N-1\big)_4^{-1}},
\\
C_{1^23, 23}=G_1-3G_3+20G_4
\\
\quad
=N \big(n-1\big)_2 \big(N^{2}n-10N^{2}+8Nn-8n+10N\big) \big(N-n\big)n^{-4} \big(N-1\big)_4^{-1},
\\
C_{1^23, 14}=2C_{13, 4},
\\
C_{1^23, 12^2}=2C_{13, 2^2},
\\
C_{1^23, 1^23}=C_{3, 3},
\\
C_{1^23, 1^32}=C_{1^23, 1^5}=0,
\\
C_{1^32, 5}=F_4-F_5={N \big(n-1\big) \big(N-n\big) \big(N^{2}+6n^{2}-6Nn-6n+5N\big)n^{-4} \big(N-1\big)_4^{-1}},
\\
C_{1^32, 23}=G_3-10G_4
\\
\quad
=-N \big(n-1\big) \big(N-n\big) \big(-4N^{2}n+5Nn^{2}+12n-12Nn-10N
\\
\quad
+10N^{2}\big)n^{-4} \big(N-1\big)_4^{-1},
\\
C_{1^32, 14}=3C_{1^22, 4},
\\
C_{1^32, 12^2}=3C_{1^22, 2^2},
\\
C_{1^32, 1^23}=3C_{12, 3},
\\
C_{1^32, 1^32}=C_{2, 2},
\\
C_{1^32, 1^5}=0,
\\
C_{1^5, 5}=F_5={ \big(N-2n\big) \big(N-n\big) \big(N^{2}+12n^{2}-12nN+5N\big)n^{-4} \big(N-1\big)_4^{-1}},
\\
C_{1^5, 23}=10G_4=10N\big(n-1\big)\big(N-n\big)_2\big(N-2n\big)n^{-4}\big(N-1\big)_4^{-1},
\\
C_{1^5, 14}=5C_{1^4, 4},
\\
C_{1^5, 12^2}=5C_{1^4, 2^2},
\\
C_{1^5, 1^23}=10C_{1^3, 3},
\\
C_{1^5, 1^32}=10C_{1^2, 2},
\\
C_{1^5, 1^5}=1,
\\
C_{6, 6}=F_1-6F_2+15F_3-20F_4+15F_5-5F_6
\\
\quad
=N \big(n-1\big) \big(n^{4}N^{4}-9n^{4}N^{3}+31n^{4}N^{2}-39Nn^{4}+40n^{4}-5n^{3}N^{4}+30n^{3}N^{3}
\\
\quad
-80n^{3} -95N^{2}n^{3}+30n^{3}N+120n^{2}N^{2}+20n^{2}+100n^{2}N+10n^{2}N^{4}
\\
\quad
-10N^{3}n^{2}-75nN-10nN^{4}-75nN^{3}-100nN^{2}+20n-20N+80N^{3}
\\
\quad
+55N^{2}+5N^{4}\big)n^{-5} \big(N-1\big)_5^{-1},
\\
C_{6, 24}=15\big(G_1-4G_2+G_6-5G_7\big)
\\
\quad
=15N \big(n-1\big) \big(N-n\big) \big(n^{3}N^{3}+29n^{3}N-40n^{3}-8N^{2}n^{3}-50n^{2}N+80n^{2}
\\
\quad
+16n^{2}N^{2}-4N^{3}n^{2}-10nN-20n-nN^{2}+7nN^{3}-5N^{3}+35N-20
\\
\quad
-10N^{2}\big)n^{-5} \big(N-1\big)_5^{-1},
\\
C_{6, 3^2}=10\big(-2G_8+6G_4-5G_9\big)
\\
\quad
=10N \big(n-1\big) \big(N-n\big) \big(-2N^{3}n^{2}-24n^{3}N+40n^{3}+4N^{2}n^{3}-80n^{2}-2n^{2}N^{2}
\\
\quad
+40n^{2}N
+20n-11nN^{2}+6nN^{3}+5nN-25N+10N^{2}-5N^{3}
\\
\quad
+20\big)n^{-5} \big(N-1\big)_5^{-1},
\\
C_{6, 2^3}=15\big(3H_1-5H_2\big)
\\
\quad
=15N^{2} \big(n-1\big)_2 \big(3nN-5N+10-10n\big) \big(N-n\big)_2n^{-5} \big(N-1\big)_5^{-1},
\\
C_{24, 6}=F_2-5F_3+10F_4-9F_5+3F_6
\\
\quad
=N \big(n-1\big) \big(N-n\big) \big(n^{3}N^{3}+11n^{3}N-8n^{3}-4N^{2}n^{3}-4n^{2}N^{2}-8n^{2}-8n^{2}N
\\
\quad
-4N^{3}n^{2}
+12n+42nN^{2}+12nN+6nN^{3}-33N-48N^{2}-3N^{3}
\\
\quad
+12\big)n^{-5} \big(N-1\big)_5^{-1},
\\
C_{24, 24}=G_5-5G_1+6G_{10}-9G_6+45G_7
\\
\quad
=N\big(n-1\big)\big(n^4N^4-9n^4N^3+53n^4N^2-135Nn^4+120n^4-12n^3N^3
\\
\quad
+41N^2n^3-5n^3N^4-60n^3N+120n^3+3N^3n^2+51n^2N^2+30n^2N^4
\\
\quad
-210n^2N-180n^2-36nN^3+495nN-63nN^4-180n+45N^4+180N
\\
\quad
-315N^2+90N^3\big)n^{-5}\big(N-1\big)_5^{-1},
\\
C_{24, 3^2}=-4G_1+4G_3-36G_4+30G_9
\\
\quad
=-2N \big(n-1\big) \big(N-n\big) \big(2n^{3}N^{3}-40n^{3}-14N^{2}n^{3}+40n^{3}N-8N^{3}n^{2}
\\
\quad
+22n^{2}N^{2}-40n^{2}-10n^{2}N+18nN^{3}+60n-33nN^{2}+15nN-15N^{3}
\\
\quad
-75N+30N^{2}+60\big)n^{-5} \big(N-1\big)_5^{-1},
\\
C_{24, 2^3}=3\big(2H_3-9H_1+15H_2\big)
\\
\quad
=3N^{2} \big(n-1\big)_2 \big(N-n\big) \big(2n^{2}N^{2}+-9n^{2}N+10n^{2}+24nN-9nN^{2}-45N
\\
\quad
+15N^{2}+30\big)n^{-5} \big(N-1\big)_5^{-1},
\\
C_{3^2, 6}=F_2-6F_3+13F_4-12F_5+4F_6
\\
\quad
=N \big(n-1 \big)_2 \big(N-n\big) \big(n^{2}N^{3}-4n^{2}-2n^{2}N^{2}+5n^{2}N
\\
\quad
-12nN^{2}+3Nn-12n-3nN^{3}+2N^{3}-8+22N+32N^{2} \big)n^{-5}\big(N-1\big)_5^{-1},
\\
C_{3^2, 24}=-6G_1+3G_{11}-12G_{12}+60G_7
\\
\quad
=-3N \big(n-1 \big)_2\big(N-n\big) \big(2n^{2}N^{3}+20n^{2}-5n^{2}N^{2}-5n^{2}N-9n
N^{3}+12nN^{2}\
\
\quad
-15Nn+60n+40-70N+20N^{2}+10N^{3} \big)n^{-5} \big(N-1\big)_5^{-1},
\\
C_{3^2, 3^2}=G_5-6G_1+2G_{13}-48G_4+40G_9
\\
\quad
=N\big(n-1\big)_2\big(N^4n^3-10Nn^3-40n^3+25N^2n^3-8N^3n^3-26N^2n^2
\\
\quad
-4N^4n^2+8N^3n^2+70Nn^2-120n^2-4N^3n+220Nn-70N^2n-80n
\\
\quad
+14N^4n+80N-100N^2-20N^4+40N^3\big)n^{-5}\big(N-1\big)_5^{-1},
\\
C_{3^2, 2^3}=3\big(3H_3-12H_1+20H_2\big)
\\
\quad
=3N^{2} \big(n-1 \big)_2\big(N-n\big) \big(3n^{2}N^{2}+20n^{2}-15n^{2}N-12nN^{2}+32Nn
\\
\quad
+40+20N^{2}-60N \big)n^{-5} \big(N-1\big)_5^{-1},
\\
C_{2^3, 6}=F_3-3F_4+3F_5-F_6
\\
\quad
=-N \big(n-1 \big) \big(N-n \big) \big(-n^{2}N^{3}+2n^{3}N^{2}-6n^{3}N+4n^{3}-4n^{2}N^{2}+4n^{2}
\\
\quad
+n^{2}N+14nN^{2}+4n+2nN^{3}+4Nn-N^{3}+4-16N^{2}-11N \big)n^{-5} \big(N-1\big)_5^{-1}
\\
C_{2^3, 24}=3\big(G_1-G_{14}+G_6-5G_7\big)
\\
\quad
=3N \big(n-1 \big) \big(N-n \big)\big(n^{3}N^{3}-20n^{3}-7n^{3}N^{2}+20n^{3}N-3n^{2}N^{3}+2n^{2}N^{2}
\\
\quad
+15n^{2}N-20n^{2}-10Nn-nN^{2}+7nN^{3}-20n+35N-5N^{3}-20
\\
\quad
-10N^{2} \big)n^{-5} \big(N-1\big)_5^{-1},
\\
C_{2^3, 3^2}=2\big(-3G_2+6G_4-5G_9\big)
\\
\quad
=-2N \big(n-1 \big) \big(N-n \big)_2\big(3n^{2}N^{2}+20n^{2}-15n^{2}N-6nN^{2}+10Nn-5N
\\
\quad
+5N^{2}+20 \big)n^{-5} \big(N-1\big)_5^{-1},
\\
C_{2^3, 2^3}=H_4-3H_3+9H_1-15H_2
\\
\quad
=N^{2} \big(n-1 \big)_2  \big(n^{3}N^{3}-30n^{3}-9n^{3}N^{2}+29n^{3}N+9n^{2}N^{2}-6 n^{2}N-3n^{2}N^{3}
\\
\quad
+30n-9nN^{2}-45Nn+9nN^{3}-30N+45N^{2}-15N^{3} \big)n^{-5} \big(N-1\big)_5^{-1}$.

This completes $\left\{ C_{{\bm \pi}_{-}, {\bm \pi}_{-}} \right\}$ for ${\bf C}_6$.
We now give $\left\{ C_{{\bm \pi}_{+}, {\bm \pi}_{-}} \right\}$ for ${\bf C}_6$:
\\
$C_{15, 6}=F_2-5F_3+10F_4-10F_5+4F_6
\\
\quad
=\big(N-n\big)\big(N+1\big)_2\big(n-1\big)_2\big(n^2N^2-5n^2N+16n^2-2nN^2-10nN
\\
\quad
-12n+2N^2+30N-8\big)n^{-5}\big(N-1\big)_5^{-1},
\\
C_{15, 24}=5\big(-G_1+6G_2-2G_{17}+12G_7\big)
\\
\quad
=-5N\big(n-1\big)_2\big(N-n\big)\big(n^2N^3-6n^2N^2+29n^2N-48n^2-4nN^3+36n
\\
\quad
-8nN+12N^2+6N^3+24-42N\big)n^{-5}\big(N-1\big)_5^{-1},
\\
C_{15, 3^2}=10\big(G8-4G_4+4G_{18}\big)
\\
\quad
=-10N\big(n-1\big)_2\big(N-n\big)\big(-nN^3+2n^2N^2-10n^2N+16n^2-12n+5nN
\\
\quad
+2N^3-4N^2-8+10N\big)n^{-5}\big(N-1\big)_5^{-1},
\\
C_{15, 2^3}=30\big(-H_1+2H_2\big)
\\
\quad
=-30N^2\big(n-1\big)_2\big(N-n\big)_2\big(nN-3n+4-2N\big)n^{-5}\big(N-1\big)_5^{-1},
\\
C_{1^24, 6}=F_3-4F_4+6F_5-3F_6
\\
\quad
=-N\big(n-1\big)\big(N-n\big)\big(-N^3n^2+22n^3+2N^2n^3-10N^2n^2-47Nn^2
\\
\quad
-38n^2+3N^3n+87Nn+12n+42N^2n-33N-3N^3+12
\\
\quad
-48N^2\big)n^{-5}\big(N-1\big)_5^{-1},
\\
C_{1^24, 24}=G_1-12G_2+6G_6-45G_7
\\
\quad
=N\big(n-1\big)\big(N-n\big)\big(-12N^3n^2+47Nn^3-90n^3+N^3n^3+330n^2-132Nn^2
\\
\quad
-24N^2n^2-105Nn-180n+42N^3n+99N^2n+315N-90N^2-45N^3
\\
\quad
-180\big)n^{-5}\big(N-1\big)_5^{-1},
\\
C_{1^24, 3^2}=-4G_8+24G_4-30G_9
\\
\quad
=2N\big(n-1\big)\big(N-n\big)\big(-2N^3n^2-12Nn^3+20n^3+4N^2n^3+62Nn^2-20N^2n^2
\\
\quad
-100n^2+60n-15Nn+3N^2n+12N^3n-75N+30N^2-15N^3
\\
\quad
+60\big)n^{-5}\big(N-1\big)_5^{-1},
\\
C_{1^24, 2^3}=18H_1-45H_2
\\
\quad
=9N^2\big(n-1\big)_2\big(2Nn-5n+10-5N\big)\big(N-n-1\big)\big(N-n\big)n^{-5}\big(N-1\big)_5^{-1},
\\
C_{123, 6}=F_3-4F_4+5F_5-2F_6
\\
\quad
=-\big(N+1\big)_2\big(n-1\big)_2\big(N-n\big)\big(-N^2n-2n^2+2Nn^2-5Nn-6n+N^2
\\
\quad
+15N-4\big)n^{-5}\big(N-1\big)_5^{-1},
\\
C_{123, 24}=G_1-6G_{15}+5G_6-30G_7
\\
\quad
= N\big(n-1\big)_2\big(N-n\big)\big(N^3n^2-25Nn^2+30n^2+6N^2n^2-20Nn
\\
\quad
-10N^3n+90n-105N+30N^2+15N^3+60\big)n^{-5}\big(N-1\big)_5^{-1},
\\
C_{123, 3^2}=G_1-G_{16}+20G_4-20G_9
\\
\quad
=N\big(n-1\big)_2\big(N-n\big)\big(N^3n^2+15Nn^2-20n^2-4N^2n^2-5N^3n+25Nn
\\
\quad
-60n-40+50N+10N^3-20N^2\big)n^{-5}\big(N-1\big)_5^{-1},
\\
C_{123, 2^3}=-3H_3+15H_1-30H_2
\\
\quad
=-3N^2\big(n-1\big)_2\big(N-n\big)\big(N^2n^2-4Nn^2+5n^2-5N^2n+5n+10Nn
\\
\quad
+20-30N+10N^2\big)n^{-5}\big(N-1\big)_5^{-1},
\\
C_{1^33, 6}=F_4-3F_5+2F_6
\\
\quad
=N\big(n-1\big)_2\big(N-n\big)\big(N^3+18n^2+6Nn^2-6n-36Nn-6N^2n-4
\\
\quad
+11N+16N^2\big)n^{-5}\big(N-1\big)_5^{-1},
\\
C_{1^33, 24}=3G_2-3G_6+30G_7
\\
\quad
=-3N\big(n-1\big)_2\big(N-n\big)\big(-N^3n+9Nn^2-10n^2+N^2n^2-5Nn-12N^2n
\\
\quad
+30n+20-35N+10N^2+5N^3\big)n^{-5}\big(N-1\big)_5^{-1},
\\
C_{1^33, 3^2}=G_8-12G_4+20G_9
\\
\quad
=-N\big(n-1\big)_2\big(N-n\big)\big(-N^3n+6Nn^2+2N^2n^2+45Nn-60n-24N^2n
\\
\quad
-40-20N^2+50N+10N^3\big)n^{-5}\big(N-1\big)_5^{-1},
\\
C_{1^33, 2^3}=-9H_1+30H_2
\\
\quad
= -3N^2\big(n-1\big)_2\big(N-n\big)_2\big(3Nn-5n+20-10N\big)n^{-5}\big(N-1\big)_5^{-1},
\\
C_{1^22^2, 6}=F_4-2F_5+F_6
\\
\quad
=N\big(n-1\big)\big(N-n\big)\big(N^3n-6n^3+6Nn^3-6N^2n^2-12Nn^2-6n^2
\\
\quad
+29Nn+4n+14N^2n-16N^2-11N-N^3+4\big)n^{-5}\big(N-1\big)_5^{-1},
\\
C_{1^22^2, 24}=G_{14}-2G_6+15G_7
\\
\quad
=-N\big(n-1\big)\big(N-n\big)\big(-3N^3n^2-9Nn^3+10n^3+5N^2n^3+21Nn^2-22N^2n^2
\\
\quad
+10n^2-35Nn+33N^2n-60n+14N^3n+105N-30N^2-15N^3
\\
\quad
-60\big)n^{-5}\big(N-1\big)_5^{-1},
\\
C_{1^22^2, 3^2}=2G_2-8G_4+10G_9
\\
\quad
=2N\big(n-1\big)\big(N-n\big)_2\big(N^2n^2-Nn^2-4N^2n-5N+20
\\
\quad
+5N^2\big)n^{-5}\big(N-1\big)_5^{-1},
\\
C_{1^22^2, 2^3}=H_3-6H_1+15H_2
\\
\quad
=N^2\big(n-1\big)_2\big(N-n\big)\big(N^2n^2-3Nn^2+5n^2-6N^2n+6Nn+15n+15N^2
\\
\quad
+30-45N\big)n^{-5}\big(N-1\big)_5^{-1},
\\
C_{1^42, 6}=F_5-F_6
\\
\quad
=N\big(n-1\big)\big(N-n\big)\big(N-2n+1\big)\big(N^2+12n^2-12Nn-12n-4+15N\big)n^{-5}\big(N-1\big)_5^{-1},
\\
C_{1^42, 24}=G_6-15G_7
\\
\quad
=N\big(n-1\big)\big(N-n\big)\big(7N^3n+18Nn^3+60n^2-24N^2n^2-54Nn^2-60n
\\
\quad
+69N^2n-40Nn+105N-30N^2-15N^3-60\big)n^{-5}\big(N-1\big)_5^{-1},
\\
C_{1^42, 3^2}=2\big(2G_4-5G_9\big)
\\
\quad
=2N\big(n-1\big)\big(N-n\big)_2\big(2N^2n-4Nn^2+10Nn-20-5N^2
\\
\quad
+5N\big)n^{-5}\big(N-1\big)_5^{-1},
\\
C_{1^42, 2^3}=3\big(H_1-5H_2\big)
\\
\quad
= 3N^2\big(n-1\big)_2\big(N-n\big)_2\big(Nn+10-5N\big)n^{-5}\big(N-1\big)_5^{-1},
\\
C_{1^6, 6}=F_6
\\
\quad
=\big(N-n\big)\big(-30N^3n+120n^4-240Nn^3+150N^2n^2+90Nn^2-90N^2n
\\
\quad
-4N+11N^2+16N^3+N^4\big)n^{-5}\big(N-1\big)_5^{-1},
\\
C_{1^6, 24}=15G_7
\\
\quad
=15N\big(n-1\big)\big(N-n\big)_2\big(N^2+6n^2-6Nn+3N-4\big)n^{-5}\big(N-1\big)_5^{-1},
\\
C_{1^6, 3^2}=10G_9
\\
\quad
=10N\big(n-1\big)\big(N-n\big)_2\big(N^2+4n^2-4Nn-N+4\big)n^{-5}\big(N-1\big)_5^{-1},
\\
C_{1^6, 2^3}=15H_2
\\
\quad
= 15N^2\big(n-1\big)_2\big(N-n\big)_3n^{-5}\big(N-1\big)_5^{-1}$.

Finally, $\left\{ C_{{\bm \pi}_{+}, {\bm \pi}_{+}} \right\}$ for ${\bf C}_6$ are obtainable using (3.7).

\section*{Appendix B}

This gives $C_{{\bm \pi}, {\bm \pi}'}$ for ${\bf C}_{--r}$ of (3.4) and ${\bf C}_r$ of (1.7) when $N =\infty$.

The results of Appendix A hold with
\\
$F_1=1,
\
F_2=n^{-1},
\
F_3=n^{-2},
\
F_4=n^{-3},
\
F_5=n^{-4},
\
F_6=n^{-5},
\
G_1 = \big(n-1\big)n^{-2},
\
G_2 = \big(n-1\big)n^{-3},
\
G_3 = 4G_2,
\
G_4 = \big(n-1\big)n^{-4},
\\
G_5 = \big(n-1\big)n^{-1},
\
G_6 = 7G_4,
\
G_7 = \big(n-1\big)n^{-5},
\\
G_8 = G_2,
\
G_9 = G_7,
\
G_{10} = 5G_2,
\\
G_{11} = 13G_2,
\
G_{12} = 7G_4,
\
G_{13} = 11G_2,
\
G_{14} = 3G_2,
\
G_{15} = 2G_2,
\
G_{16} = 7G_2,
\\
H_1 = \big(n-1\big)_2n^{-4},
\
H_2 = \big(n-1\big)_2n^{-5},
\
H_3 = \big(n-1\big)_2n^{-3},
\
H_4 = \big(n-1\big)_2n^{-2},
\\
C_{2, 2} = \big(n-1\big)n^{-1},
\
C_{1^2, 2} = n^{-1},
\
C_{1^2, 1^2} = 1,
\\
C_{3, 3} = \big(n-1\big)_2n^{-2},
\
C_{12, 3} = \big(n-1\big)n^{-2},
\
C_{12, 12} = \big(n-1\big)n^{-1},
\
C_{12, 1^3} = 0,
\\
C_{1^3, 3} = n^{-2},
\
C_{1^3, 12} = 3n^{-1},
\
C_{1^3, 1^3} = 1,
\\
C_{4, 4} = \big(n-1\big)\big(n^2-3n+3\big)n^{-3},
\
C_{4, 2^2} = 3\big(n-1\big)\big(2n-3\big)n^{-3},
\\
C_{2^2, 4} = \big(n-1\big)^2n^{-3},
\
C_{2^2, 2^2} = \big(n-1\big)\big(n^2-2n+3\big)n^{-3},
\\
C_{13, 4} = \big(n-1\big)_2n^{-3},
\
C_{13, 2^2} = -3\big(n-1\big)_2n^{-3},
\
C_{13, 13} = C_{3, 3},
\\
C_{13, 1^22} = C_{13, 1^4} = 0,
\\
C_{1^22, 4} = \big(n-1\big)n^{-3},
\
C_{1^22, 2^2} = \big(n-1\big)\big(n-3\big)n^{-3},
\
C_{1^22, 13} = 2\big(n-1\big)n^{-2},
\\
C_{1^22, 1^22} = \big(n-1\big)n^{-1},
\
C_{1^22, 1^4} = 0,
\\
C_{1^4, 4} = n^{-3},
\
C_{1^4, 2^2} = 3\big(n-1\big)n^{-3},
\
C_{1^4, 13} = 4n^{-2},
\
C_{1^4, 1^22} = 6n^{-1},
\
C_{1^4, 1^4} = 1,
\\
C_{5, 5} = \big(n-1\big)_2\big(n^2-2n+2\big)n^{-4},
\
C_{5, 23} = 10\big(n-1\big)\big(n-2\big)^2n^{-4},
\\
C_{23, 5} = \big(n-2\big)\big(n-1\big)^2n^{-4},
\
C_{23, 23} = \big(n-1\big)_2\big(n^2-5n+10\big)n^{-4},
\\
C_{14, 5} = \big(n-1\big)\big(n^2-3n+3\big)n^{-4},
\
C_{14, 23} = -2\big(n-1\big)\big(2n^2-12n+15\big)n^{-4},
\\
C_{14, 14} = C_{4, 4},
\
C_{14, 12^2} = C_{4, 2^2},
\
C_{14, 1^23} = C_{14, 1^32} = C_{14, 1^5} = 0,
\\
C_{12^2, 5} = \big(n-1\big)^2n^{-4},
\
C_{12^2, 23} = 2\big(n-1\big)\big(n^2-4n+5\big)n^{-4},
\
C_{12^2, 14} = C_{2^2, 4},
\\
C_{12^2, 12^2} = C_{2^2, 2^2},
\
C_{12^2, 1^23} = C_{12^2, 1^32} = C_{12^2, 1^5} = 0,
\\
C_{1^23, 5} = \big(n-1\big)_2n^{-4},
\
C_{1^23, 23} = \big(n-10\big)\big(n-1\big)_2n^{-4},
\
C_{1^23, 14} = 2C_{31, 4},
\\
C_{1^23, 12^2} = 2C_{13, 2^2},
\
C_{1^23, 1^23} = C_{3, 3},
\
C_{1^23, 1^32} = C_{1^23, 1^5} = 0,
\\
C_{1^32, 5} = \big(n-1\big)n^{-4},
\
C_{1^32, 23} = 2\big(n-1\big)\big(2n-5\big)n^{-4},
\
C_{1^32, 14} = 3C_{1^22, 4},
\\
C_{1^32, 12^2} = 3C_{1^22, 2^2},
\
C_{1^32, 1^23} = 3C_{12, 3},
\
C_{1^32, 1^32} = C_{2, 2},
\
C_{1^32, 1^5} = 0,
\\
C_{1^5, 5} = n^{-4},
\
C_{1^5, 23} = 10\big(n-1\big)n^{-4},
\
C_{1^5, 14} = 5C_{1^4, 4},
\
C_{1^5, 12^2} = 15C_{1^4, 2^2},
\\
C_{1^5, 1^23} = 10C_{1^3, 3},
\
C_{1^5, 1^32} = 10C_{1^2, 2},
\
C_{1^5, 1^5} = 1,
\\
C_{6, 6} = \big(n-1\big)\big(n^4-5n^3+10n^2-10n+5\big)n^{-5},
\\
C_{6, 24} = 15\big(n-1\big)\big(n^3-4n^2+7n-5\big)n^{-5},
\\
C_{6, 3^2} = -10\big(n-1\big)\big(2n^2-6n+5\big)n^{-5},
\\
C_{6, 2^3} = 15\big(3n-5\big)\big(n-1\big)_2n^{-5},
\\
C_{24, 6} = \big(n^2-3n+3\big)\big(n-1\big)^2n^{-5},
\\
C_{24, 24} = \big(n-1\big)\big(n^4-5n^3+30n^2-63n+45\big)n^{-5},
\\
C_{24, 3^2} = -2\big(n-1\big)\big(2n^3-8n^2+18n-15\big)n^{-5},
\\
C_{24, 2^3} = 3\big(n-1\big)_2\big(2n^2-9n+15\big)n^{-5},
\\
C_{3^2, 6} = \big(n-1\big)^2\big(n-2\big)^2n^{-5},
\\
C_{3^2, 24} = -3\big(n-1\big)\big(2n-5\big)\big(n-2\big)^2n^{-5},
\\
C_{3^2, 3^2} = \big(n-1\big)\big(n^2-2n+10\big)\big(n-2\big)^2n^{-5},
\\
C_{3^2, 2^3} = 3\big(n-1\big)_2\big(3n^2-12n+20\big)n^{-5},
\\
C_{2^3, 6} = \big(n-1\big)^3n^{-5},
\\
C_{2^3, 24} = 3\big(n^2-2n+5\big)\big(n-1\big)^2n^{-5},
\\
C_{2^3, 3^2} = -2\big(n-1\big)\big(3n^2-6n+5\big)n^{-5},
\\
C_{2^3, 2^3} = \big(n-1\big)_2\big(n^3-3n^2+9n-15\big)n^{-5},
\\
C_{1^24, 6}= \big(n-1\big)\big(n^2-3n+3\big)n^{-5},
\\
C_{1^24, 24}= \big(n-1\big)\big(n-3\big)\big(n^2-9n+15\big)n^{-5},
\\
C_{1^24, 3^2}= -2\big(n-1\big)\big(2n^2-12n+15\big)n^{-5},
\\
C_{1^24, ^3}= 9\big(2n-5\big)\big(n-1\big)_2n^{-5},
\\
C_{123, 6}= \big(n-2\big)\big(n-1\big)^2n^{-5},
\\
C_{123, 24}= \big(n-1\big)_2\big(n^2-10n+15\big)n^{-5},
\\
C_{123, 3^2}= \big(n-1\big)_2\big(n^2-5n+10\big)n^{-5},
\\
C_{123, 2^3}= -3\big(n-1\big)_2\big(n^2-5n+10\big)n^{-5},
\\
C_{1^33, 6}= \big(n-1\big)_2n^{-5},
\\
C_{1^33, 24}= 3\big(n-5\big)\big(n-1\big)_2n^{-5},
\\
C_{1^33, 3^2}= \big(n-10\big)\big(n-1\big)_2n^{-5},
\\
C_{1^33, 2^3}= -3\big(3n-10\big)\big(n-1\big)_2n^{-5},
\\
C_{1^22^2, 6}= \big(n-1\big)^2n^{-5},
\\
C_{1^22^2, 24}= \big(3n-5\big)\big(n-1\big)\big(n-3\big)n^{-5},
\\
C_{1^22^2, 3^2}= 2\big(n-1\big)\big(n^2-4n+5\big)n^{-5},
\\
C_{1^22^2, 2^3}= \big(n-1\big)_2\big(n^2-6n+15\big)n^{-5},
\\
C_{1^42, 6}= \big(n-1\big)n^{-5},
\\
C_{1^42, 24}= \big(n-1\big)\big(7n-15\big)n^{-5},
\\
C_{1^42, 3^2}= 2\big(n-1\big)\big(2n-5\big)n^{-5},
\\
C_{1^42, 2^3}= 3\big(n-5\big)\big(n-1\big)_2n^{-5},
\\
C_{1^6, 6}= n^{-5},
\\
C_{1^6, 24}= 15\big(n-1\big)n^{-5},
\\
C_{1^6, 3^2}= 10\big(n-1\big)n^{-5},
\\
C_{1^6, 2^3}= 15\big(n-1\big)_2n^{-5}$.

\section*{Appendix C}

Here, we give for $r\leq6$,
${\bf C}_{--r}^{-1} {\bf C}_{1, 1}^{-1} = {\bf C}_r^{1, 1} = \left( C^{{\bm \pi}_{-}, {\bm \pi}'_{-}} \right)$
for the case $N = \infty$.
By (3.4) this is needed for the UE of ${\bm \mu}_{-(r)}$.
By the inversion principle, these elements are
$C^{{\bm \pi}, {\bm \pi}'} (\infty, n) = C_{{\bm \pi}, {\bm \pi}'} (n, \infty)$.
We also give the other elements of ${\bf C}_r^{-1}$,
$C_r^{2, 1}$ and $C_r^{2, 2}$ for $r\leq5$ defined by (2.9).
These are needed for the UE of ${\bm \mu}_{+(r)}$ in (3.8).
Set ${\bf B}_r (N, n) = {\bf B}_r$, ${\bf C}_r(N, n) = {\bf C}_r$.
Then ${\bf B}_r^{-1} = {\bf B}_r (n, N)$, ${\bf C}_r^{-1} = {\bf C}_r (n, N)$,

\noindent
$C_2^{1, 1}  = C_{2, 2}^{-1}=C^{2, 2}=n\big(n-1\big)^{-1},
\\
C_2^{2, 1}=C^{1^2, 2}=-\big(n-1\big)^{-1},
\\
C_2^{2, 2}=C^{1^2, 1^2} = 1$.

\noindent
$C_3^{1, 1}=C_{3, 3}^{-1}=C^{3, 3}=\big(n-1\big)_2^{-1}n^2,
{\bf C}_3^{2, 1} = \left(
\begin{array}{c}
C^{12, 3}
\\
C^{1^3, 3}
\end{array} \right),
\\
C^{12, 3}= -n\big(n-1\big)_2^{-1},
\\
C^{1^3, 3}= 2\big(n-1\big)_2^{-1},
\\
{\bf C}_3^{2, 2}=\left(
\begin{array}{cc}
C^{12, 12} & C^{12, 1^3}
\\
C^{1^3, 12} & C^{1^3, 1^3}
\end{array} \right),
\\
C^{12, 12}=n\big(n-1\big)^{-1},
\\
C^{12, 1^3}=0,
\\
C^{1^3, 12}=-3\big(n-1\big)^{-1},
\\
C^{1^3, 1^3}=1$.

\noindent
${\bf C}_4^{11}=\left(
\begin{array}{cc}
C^{4, 4} & C^{4, 2^2}
\\
C^{2^2, 4} & C^{2^2, 2^2}
\end{array} \right),
\\
C^{4, 4}=n\big(n^2-2n+3\big)\big(n-1\big)_3^{-1},
\\
C^{4, 2^2}=-3n\big(2n-3\big)\big(n-1\big)_3^{-1},
\\
C^{2^2, 4}=-n\big(n-2\big)_2^{-1},
\\
C^{2^2, 2^2}=n\big(n^2-3n+3\big)\big(n-1\big)_3^{-1},
\\
{\bf C}_4^{2, 1}=\left( \begin{array}{cc}
C^{13, 4} & C^{13, 2^2}
\\
C^{1^22, 4} & C^{1^22, 2^2}
\\
C^{1^4, 4} & C^{1^4, 2^2}
\end{array} \right),
\\
C^{13, 4}=-n\big(1+n\big)\big(n-1\big)_3^{-1},
\\
C^{13, 2^2}=3n\big(n-2\big)_2^{-1},
\\
C^{1^22, 4}=2n\big(n-1\big)_3^{-1},
\\
C^{1^22, 2^2}=-n^2\big(n-1\big)_3^{-1},
\\
C^{1^4, 4}=-6\big(n-1\big)_3^{-1},
\\
C^{1^4, 2^2}=3n\big(n-1\big)_3^{-1},
\\
{\bf C}_4^{2, 2} =\left( \begin{array}{ccc}
C^{13, 13} & C^{13, 1^22} & C^{13, 1^4}
\\
C^{1^22, 13} & C^{1^22, 1^22} & C^{1^22, 1^4}
\\
C^{1^4, 13} & C^{1^4, 1^22} & C^{1^4, 1^4}
\end{array} \right),
\\
C^{13, 13}=n^2\big(n-1\big)_2^{-1},
\\
C^{13, 1^22}=0,
\\
C^{13, 1^4}=0,
\\
C^{1^22, 13}=-2n\big(n-1\big)_2^{-1},
\\
C^{1^22, 1^22}=n\big(n-1\big)^{-1},
\\
C^{1^22, 1^4}=0,
\\
C^{1^4, 13}=8\big(n-1\big)_2^{-1},
\\
C^{1^4, 1^22}=-6\big(n-1\big)^{-1},
\\
C^{1^4, 1^4}=1$.

\noindent
${\bf C}_5^{1, 1} =\left( \begin{array}{cc}
C^{5, 5} & C^{5, 23}
\\
C^{23, 5} & C^{23, 23}
\end{array} \right),
\\
C^{5, 5}=n^2\big(n^2-5n+10\big)\big(n-1\big)_4^{-1},
\\
C^{5, 23}=-10n^2\big(n-1\big)^{-1}\big(n-3\big)_2^{-1},
\\
C^{23, 5}=-n^2\big(n-2\big)_3^{-1},
\\
C^{23, 23}=n^2\big(n^2-2n+2\big)\big(n-1\big)_4^{-1},
\\
{\bf C}_5^{2, 1} = \left( \begin{array}{cc}
C^{14, 5} & C^{14, 23}
\\
C^{12^2, 5} & C^{12^2, 23}
\\
C^{1^23, 5} & C^{1^23, 23}
\\
C^{1^32, 5} & C^{1^32, 23}
\\
C^{1^5, 5} & C^{1^5, 23}
\end{array} \right),
\\
C^{14, 5}=-n\big(n^2-n+6\big)\big(n-1\big)_4^{-1},
\\
C^{14, 23}=2n\big(n+2\big)\big(2n-3\big)\big(n-1\big)_4^{-1},
\\
C^{12^2, 5}=2n\big(n-2\big)_3^{-1},
\\
C^{12^2, 23}=-2n\big(n^2-2n+2\big)\big(n-1\big)_4,
\\
C^{1^23, 5}=2n\big(n+2\big)\big(n-1\big)_4^{-1},
\\
C^{1^23, 23}=-n\big(n^2+8n-8\big)\big(n-1\big)_4^{-1},
\\
C^{1^32, 5}=-6n\big(n-1\big)_4^{-1},
\\
C^{1^32, 23}=5n^2\big(n-1\big)_4^{-1},
\\
C^{1^5, 5}=24\big(n-1\big)_4^{-1},
\\
C^{1^5, 23}=-20n\big(n-1\big)_4^{-1},
\\
{\bf C}^{2, 2}_5=\left(\begin {array}{ccccc}
C^{14, 14} & C^{14, 12^2} & C^{14, 1^23} & C^{14, 1^32} & C^{14, 1^5}
\\
C^{12^2, 14} & C^{12^2, 12^2} & C^{12^2, 1^23} & C^{12^2, 1^32} & C^{12^2, 1^5}
\\
C^{1^23, 14} & C^{1^23, 12^2} & C^{1^23, 1^23} & C^{1^23, 1^32} & C^{1^23, 1^5}
\\
C^{1^32, 14} & C^{1^32, 12^2} & C^{1^32, 1^23} & C^{1^32, 1^32} & C^{1^32, 1^5}
\\
C^{1^5, 14} & C^{1^5, 12^2} & C^{1^5, 1^23} & C^{1^5, 1^32} & C^{1^5, 1^5}
\end {array} \right),
\\
C^{14, 14} =n\big(n^2-2n+3\big)\big(n-1\big)_3^{-1},
\\
C^{14, 12^2} =-3n\big(2n-3\big)\big(n-1\big)_3^{-1},
\\
C^{14, 1^23} =0,
\\
C^{14, 1^32} =0,
\\
C^{14, 1^5}= 0,
\\
C^{12^2, 14} =-n\big(n-2\big)_2^{-1},
\\
C^{12^2, 12^2} =n\big(n^2-3n+3\big)\big(n-1\big)_3^{-1},
\\
C^{12^2, 1^23} =0,
\\
C^{12^2, 1^32} =0,
\\
C^{12^2, 1^5}=0,
\\
C^{1^23, 14} =-2n\big(n+1\big)\big(n-1\big)_3^{-1},
\\
C^{1^23, 12^2} =6n\big(n-2\big)_2^{-1},
\\
C^{1^23, 1^23} =n^2\big(n-1\big)_2^{-1},
\\
C^{1^23, 1^32} =0,
\\
C^{1^23, 1^5} =0,
\\
C^{1^32, 14} =6n\big(n-1\big)_3^{-1},
\\
C^{1^32, 12^2} =-3n^2\big(n-1\big)_3^{-1},
\\
C^{1^32, 1^23} =-3n\big(n-1\big)_2^{-1},
\\
C^{1^32, 1^32} =n\big(n-1\big)^{-1},
\\
C^{1^32, 1^5}= 0,
\\
C^{1^5, 14} =-30\big(n-1\big)_3^{-1},
\\
C^{1^5, 12^2} =15n\big(n-1\big)_3^{-1},
\\
C^{1^5, 1^23} =20\big(n-1\big)_2^{-1},
\\
C^{1^5, 1^32} =-10\big(n-1\big)^{-1},
\\
C^{1^5, 1^5}=1$.

\noindent
${\bf C}_6^{1, 1} = \left( \begin{array}{cccc}
C^{6, 6} & C^{6, 24}& C^{6, 3^2}& C^{6, 2^3}
\\
C^{24, 6} & C^{24, 24}& C^{24, 3^2}& C^{24, 2^3}
\\
C^{3^2, 6} & C^{3^2, 24}& C^{3^2, 3^2}& C^{3^2, 2^3}
\\
C^{2^3, 6} & C^{2^3, 24}& C^{2^3, 3^2}& C^{2^3, 2^3}\end{array} \right),
\\
C^{6, 6} =n\big(n^4-9n^3+31n^2-39n+40\big)\big(n-1\big)_5^{-1},
\\
C^{6, 24}=-15n\big(n^3-8n^2+29n-40\big)\big(n-1\big)_5^{-1},
\\
C^{6, 3^2}=-40n\big(n^2-6n+10\big)\big(n-1\big)_5^{-1},
\\
C^{6, 2^3}=15n^2\big(3n-10\big)\big(n-1\big)_5^{-1},
\\
C^{24, 6} =-n\big(n^2-3n+8\big)\big(n-2\big)_4^{-1},
\\
C^{24, 24}=n\big(n^4-9n^3+53n^2-135n+120\big)\big(n-1\big)_5^{-1},
\\
C^{24, 3^2}=4n\big(n^2-5n+10\big)\big(n-1\big)^{-1}\big(n-3\big)_3^{-1},
\\
C^{24, 2^3}=-3n^2\big(2n-5\big)\big(n-1\big)^{-1}\big(n-3\big)_3^{-1},
\\
C^{3^2, 6} =-n\big(n^2-n+4\big)\big(n-2\big)_4^{-1},
\\
C^{3^2, 24}=3n\big(2n^3-5n^2-5n+20\big)\big(n-1\big)_5^{-1},
\\
C^{3^2, 3^2}=n\big(n^4-8n^3+25n^2-10n-40\big)\big(n-1\big)_5^{-1},
\\
C^{3^2, 2^3}=-3n^2\big(3n^2-15n+20\big)\big(n-1\big)_5^{-1},
\\
C^{2^3, 6} =2n\big(n-3\big)_3^{-1},
\\
C^{2^3, 24}=-3n\big(n^2-5n+10\big)\big(n-1\big)^{-1}\big(n-3\big)_3^{-1},
\\
C^{2^3, 3^2}=-2n\big(3n^2-15n+20\big)\big(n-1\big)_5^{-1},
\\
C^{2^3, 2^3}=n^2\big(n^2-7n+15\big)\big(n-1\big)^{-1}\big(n-3\big)_3^{-1},
\\
C^{1^24, 6}= 2n\big(n^2+11\big)\big(n-1\big)_5^{-1},
\\
C^{1^24, 24}= -n\big(n^3+47n-90\big)\big(n-1\big)_5^{-1},
\\
C^{1^24, 3^2}= -4\big(n+5\big)n\big(n-3\big)^{-1}_3\big(n-1\big)^{-1},
\\
C^{1^24, 2^3}= 9n^2\big(2n-5\big)\big(n-1\big)_5^{-1},
\\
C^{123, 6}= 2\big(n+1\big)n\big(n-2\big)_4^{-1},
\\
C^{123, 24}= -n\big(n^3+6n^2-25n+30\big)\big(n-1\big)_5^{-1},
\\
C^{123, 3^2}= -n\big(n^3+15n-4n^2-20\big)\big(n-1\big)_5^{-1},
\\
C^{123, 2^3}= 3n^2\big(n^2-4n+5\big)\big(n-1\big)_5^{-1},
\\
C^{1^33, 6}= -6n\big(3+n\big)\big(n-1\big)_5^{-1},
\\
C^{1^33, 24}= 3\big(n+10\big)n\big(n-2\big)_4^{-1},
\\
C^{1^33, 3^2}= 2n^2\big(3+n\big)\big(n-1\big)_5^{-1},
\\
C^{1^33, 2^3}= -3n^2\big(3n-5\big)\big(n-1\big)_5^{-1},
\\
C^{1^22^2, 6}= -6n\big(n-2\big)_4^{-1},
\\
C^{1^22^2, 24}= n\big(5n^2-9n+10\big)\big(n-1\big)_5^{-1},
\\
C^{1^22^2, 3^2}= 2n^2\big(n-2\big)_4^{-1},
\\
C^{1^22^2, 2^3}= -n^2\big(n^2-3n+5\big)\big(n-1\big)_5^{-1},
\\
C^{1^42, 6}= 24n\big(n-1\big)_5^{-1},
\\
C^{1^42, 24}= -18n^2\big(n-1\big)_5^{-1},
\\
C^{1^42, 3^2}= -8n^2\big(n-1\big)_5^{-1},
\\
C^{1^42, 2^3}= 3n^3\big(n-1\big)_5^{-1},
\\
C^{1^6, 6}= -120\big(n-1\big)_5^{-1},
\\
C^{1^5, 24}= 90n\big(n-1\big)_5^{-1},
\\
C^{1^5, 3^2}= 40n\big(n-1\big)_5^{-1},
\\
C^{1^5, 2^3}= -15n^2\big(n-1\big)_5^{-1}$.

\section*{Appendix D}

Here, we give ${\bf a} = {\bf a} ({\bm \pi})$, ${\bf b} = {\bf b} ({\bm \pi})$ needed in (1.11), the explicit
formula for the UEs of products of cumulants up to order six:
\begin{eqnarray*}
\kappa_4:
{\bf a}
&=&
\big(1,-3\big)
\\
{\bf a} {\bf C}_{--4}^{-1}
&=&
n\big(N-1\big)N^{-3}\big(n-1\big)_3^{-1}\big(J_1,J_2\big),
\\
J_1
&=&
N^2n^2-6Nn^2+6n^2+6n+N^2n-6N,
\\
J_2
&=&
-3\big(N^2n^2-N^2n-4Nn^2+6n^2+6N-6n\big).
\\
\kappa_5:
{\bf a}
&=&
\big(1,-10\big),
\\
{\bf a} {\bf C}_{--5}^{-1}
&=&
n^2\big(N-1\big)_2N^{-4}\big(n-1\big)_4^{-1}\big(J_1,J_2\big),
\\
J_1
&=&
N^2n^2-12Nn^2+12n^2+5N^2n+60n-60N,
\\
J_2
&=&
-10\big(N^2n^2-N^2n-6Nn^2+12n^2+12N-12n\big).
\\
\kappa_4\mu:
{\bf b}
&=&
\big(0,0,1,-3,0,0,0\big),
\\
{\bf b} {\bf C}_5^{-1}
&=&
n\big(N-1\big)N^{-4}\big(n-1\big)_4^{-1}\big(J_1,J_2,J_3,J_4,0,0,0\big),
\\
J_1
&=&
-\big(N-n\big)\big(N^2n^2+5N^2n-18Nn-6Nn^2-12N+6n^2+30n\big),
\\
J_2
&=&
2\big(N-n\big)\big(5N^2n^2-5N^2n+18Nn-24Nn^2+12N+30n^2
\\
&&
-30n\big),
\\
J_3
&=&
N\big(n-4\big)\big(N^2n^2-6Nn^2+6n^2+6n+N^2n-6N\big),
\\
J_4
&=&
-3N\big(n-4\big)\big(N^2n^2-N^2n-4Nn^2+6n^2+6N-6n\big).
\\
\kappa_6:
{\bf a}
&=&
\big(1,-15,-10,30\big),
\\
{\bf a} {\bf C}_{--6}^{-1}
&=&
n\big(N-1\big)\big(n-1\big)_5^{-1}N^{-5}\big(J_1,J_2,J_3,J_4\big),
\\
J_1
&=&
1890N^2n^2-1800Nn-2400Nn^2+600N^2n-90N^3n^3+480N
\\
&&
-480n+480N^2-120N^3+1320n^2+1920n^3-30N^3n^4
\\
&&
+N^4n^4+480N^2n^3-1800Nn^3-390N^3n^2-90N^3n+150N^2n^4
\\
&&
-240Nn^4+16N^4n^3+11N^4n^2-4N^4n+120n^4,
\\
J_2
&=&
-15\big(6N^2n^2+1320Nn-360N^2n-480N+480n-480N^2
\\
&&
+120N^3-840n^2 +240n^3-16N^3n^4+N^4n^4+36N^2n^3
\\
&&
-96Nn^3-26N^3n^2+90N^3n+78N^2n^4-168Nn^4+2N^4n^3
\\
&&
-7N^4n^2+4N^4n+120n^4\big),
\\
J_3
&=&
-10\big(-54N^2n^2-1080Nn+120Nn^2+240N^2n+480N-480n
\\
&&
+480N^2 -120N^3+600n^2-240n^3-12N^3n^4+N^4n^4
\\
&&
-60N^2n^3+144Nn^3+54N^3n^2  -90N^3n+66N^2n^4-144Nn^4
\\
&&
-2N^4n^3+5N^4n^2-4N^4n+120n^4\big),
\\
J_4
&=&
30n\big(N-2\big)\big(9N^2n^2-180Nn-36Nn^2+36N^2n+N^3n^3
\\
&&
+120N-120n +180n^2-60n^3-9N^2n^3+36Nn^3-3N^3n^2+2N^3n\big).
\end{eqnarray*}
\begin{eqnarray*}
\kappa_5 \mu:
{\bf b}
&=&
\big(0^4,1,-10,0^5\big),
\\
{\bf b} {\bf C}_6^{-1}
&=&
n\big(N-1\big)_2N^{-4}\big(n-1\big)_4^{-1}\big(J_1,J_2,J_3,J_4,J_5,J_6,0,0,0,0,0\big),
\\
J_1
&=&
-\big(n+1\big)\big(N-n\big)\big(n^2N^2-12n^2N+12n^2+15nN^2-60nN
\\
&&
+180n-4N^2-72N-48\big)N^{-1}\big(n-5\big)^{-1},
\\
J_2
&=&
15\big(N-n\big)\big(n^3N^2-7nN^2+4N^2+2n^2N^2-8n^3N+72N
\\
&&
-16nN+24n^2+12n^3+48-84n\big)N^{-1}\big(n-5\big)^{-1}
\\
J_3
&=&
10\big(N-n\big)\big(n^3N^2+5nN^2-4N^2-2n^2N^2+30nN-6n^3N-72N
\\
&&
-48-24n^2+12n^3+60n\big)N^{-1}\big(n-5\big)^{-1},
\\
J_4
&=&
-30\big(N-n\big)n\big(n^2N^2-3nN^2+2N^2+12nN+6N-6n^2N
\\
&&
+12n^2-36n+24\big)N^{-1}\big(n-5\big)^{-1},
\\
J_5
&=&
n\big(n^2N^2-12n^2N+12n^2+60n+5nN^2-60N\big),
\\
J_6
&=&
-10n\big(n^2N^2-nN^2-6n^2N+12N-12n+12n^2\big).
\\
\kappa_4\mu^2:
{\bf b}
&=&
\big(0^6,1,-3,0^3\big),
\\
{\bf b} {\bf C}_6^{-1}
&=&
n\big(N-1\big)N^{-5}\big(n-1\big)_5^{-1}\big(J_1,J_2,J_3,J_4,J_5,J_6,J_7,J_8,0,0,0\big),
\\
J_1
&=&
\big(N-n\big)\big(24+174Nn-28N^2n^2+84Nn^2-83N^2n+24N-66n-56N^2
\\
&&
+4N^3-96n^2-6n^3+6Nn^3-N^2n^3+2N^3n^2+18N^3n\big),
\\
J_2
&=&
-\big(N-n\big)\big(-210Nn-90N^2n^2+198Nn^2-69N^2n+N^3n^3-360N+630n
\\
&&
+360N^2-60N^3-180n^2-90n^3+84Nn^3-21N^2n^3+15N^3n^2+20N^3n
\\
&&
-360\big),
\\
J_3
&=&
-2\big(N-n\big)\big(-26N^2n^2-30Nn+6Nn^2+59N^2n+120N-150n-100N^2
\\
&&
+20N^3+60n^2-30n^3-5N^2n^3+24Nn^3+7N^3n^2-15N^3n+120\big),
\\
J_4
&=&
3n\big(N-2\big)\big(N-n\big)\big(3N^2n-12Nn^2+12Nn+30n^2-90n+30N+60
\\
&&
-10N^2+N^2n^2\big),
\\
J_5
&=&
-2N\big(n-5\big)\big(N-n\big)\big(N^2n^2+5N^2n-18Nn-6Nn^2-12N+6n^2+30n\big),
\\
J_6
&=&
4N\big(n-5\big)\big(N-n\big)\big(5N^2n^2-5N^2n+18Nn-24Nn^2+12N+30n^2
\\
&&
-30n\big),
\\
J_7
&=&
N^2\big(n-4\big)_2\big(N^2n^2-6Nn^2+6n^2+N^2n+6n-6N\big),
\nonumber
\\
J_8
&=&
-3N^2\big(n-4\big)_2\big(-N^2n-4Nn^2+6n^2+6N-6n+N^2n^2\big).
\nonumber
\end{eqnarray*}
\begin{eqnarray*}
\kappa_4 \kappa_2:
{\bf a}
&=&
\big(0,1,0,-3\big),
\\
{\bf a} {\bf C}_{--6}^{-1}
&=&
n\big(N-1\big)\big(n-1\big)_5^{-1}N^{-5}\big(J_1,J_2,J_3,J_4\big),
\\
J_1
&=&
-\big(N-n\big)\big(24-16N^2n^2+24Nn+84Nn^2-5N^2n+N^3n^3+24N-66n
\\
&&
+4N^2+4N^3-96n^2-6n^3-7N^2n^3+12Nn^3+2N^3n^2-7N^3n\big),
\\
J_2 &=& -93N^2n^2+990Nn-120N^2n+24N^3n^3-360N+360n-360N^2-60N^3
\\
&&
-630n^2-60N^4+180n^3-14N^3n^4+N^4n^4+48N^2n^3-72Nn^3-121N^3n^2
\\
&&
+255N^3n+57N^2n^4-126Nn^4-10N^4n^2+45N^4n+90n^4,
\\
J_3
&=&
2\big(N-n\big)\big(58N^2n^2+30Nn-66Nn^2-85N^2n+2N^3n^3+120 N-150n
\\
&&
+20N^2+20N^3+60n^2-30n^3-17N^2n^3+36Nn^3-5N^3n^2-5N^3n+120\big),
\\
J_4
&=&
-3n\big(N-2\big)\big(9N^2n^2-90Nn-18Nn^2+8N^2n+N^3n^3+60N-60n-20N^3
\\
&&
+90n^2-30n^3-5N^2n^3+18Nn^3-7N^3n^2+20N^3n\big).
\end{eqnarray*}
For $N=\infty$,
\begin{eqnarray*}
\kappa_4:
{\bf a}
&=&
\big(1,-3\big),
\\
{\bf a} {\bf C}_{--4}^{-1}
&=&
n^2\big(n-1\big)_3^{-1}\big(n+1,-3\big(n-1\big)\big).
\\
\kappa_5:
{\bf a}
&=&
\big(1,-10\big),
\\
{\bf a} {\bf C}_{--5}^{-1}
&=&
n^3\big(n-1\big)_4^{-1}\big(n+5,-10\big(n+1\big)\big).
\\
\kappa_4\mu:
{\bf b}
&=&
\big(0,0,1,-3,0,0,0\big),
\\
{\bf b} {\bf C}_5^{-1}
&=&
n^2\big(n-1\big)_4^{-1}\big(J_1,J_2,J_3,J_4,0,0,0\big),
\\
J_1
&=&
-\big(n+5\big),
\\
J_2
&=&
10\big(n-1\big),
\\
J_3
&=&
\big(n+1\big)\big(n-4\big),
\\
J_4
&=&
-3\big(n-1\big)\big(n-4\big).
\\
\kappa_6:
{\bf a}
&=&
\big(1,-15,-10,30\big),
\\
{\bf a} {\bf C}_{--6}^{-1}
&=&
n^2\big(n-1\big)_5^{-1}\big(J_1,J_2,J_3,J_4\big),
\\
J_1
&=&
\big(n+1\big)\big(n^2+15n-4\big),
\\
J_2
&=&
-15\big(n+4\big)\big(n-1\big)^2,
\\
J_3
&=&
-10\big(n^2-n+4\big)\big(n-1\big),
\\
J_4
&=&
30n\big(n-1\big)_2.
\\
\kappa_5 \mu:
{\bf b}
&=&
\big(0^4,1,-10,0^5\big),
\\
{\bf b} {\bf C}_6^{-1}
&=&
\big(J_1,J_2,J_3,J_4,J_5,J_6,0,0,0,0,0\big),
\\
J_1
&=&
-n\big(n+1\big)\big(n^2-4+15n\big)\big(n-1\big)_5^{-1},
\\
J_2
&=&
15n\big(n+4\big)\big(n-1\big)\big(n-2\big)_4^{-1},
\\
J_3
&=&
10\big(n^2-n+4\big)n\big(n-2\big)_4^{-1},
\\
J_4
&=&
-30n^2\big(n-3\big)_3^{-1},
\\
J_5
&=&
n^3\big(n+5\big)\big(n-1\big)^{-1}_4,
\\
J_6
&=&
-10n^3\big(n-2\big)_3^{-1}.
\end{eqnarray*}
These expressions for $\kappa_4$, $\kappa_5$, $\kappa_6$ agree with Fisher's:
see equation (12.28) of Stuart and Ord (1987).
 The UE for $\kappa_4 \mu$ above can be
derived from $\kappa_{41}=\kappa_4 \kappa_1-\kappa_5/n$ at the top of page 6 of
Wishart (1952).

\section*{Appendix E}

Set $ S_{a} = \displaystyle \sum_{1=1}^{n} X_{i}^{a}, s_{a} = \displaystyle \sum_{i=1}^{N} x_{i}^{a}$.
Here, we extend the table of Skellam (1949) for $E S_{a_{1}} \cdots S_{a_{r}}$
in terms of $\left\{ s_{a} \right\}$ up to $r = 6$.
We exclude the cases where $ a_{1}, \ldots, a_{r}$ are all the same or all different.
\begin{eqnarray*}
E S_{a}^{2} S_{b}
&=&
\lambda_{3}  s_{a+b+c} + \lambda_{12} \big( s_{2a} s_{b} + 2 s_{a+b} s_{a}  \big) + \lambda_{1^3} s_{a}^{2} s_{b},
\\
E S_{a}^{2} S_{b} S_{c}
&=&
\lambda_{4} s_{2a+b+c} + \lambda_{13} \big( \sum^{2} s_{2a+b} s_{c} + 2 s_{a+b+c} s_{a} \big)
\\
&&
+ \lambda_{2^2} \big( s_{2a} s_{b+c} + 2 s_{a+b} s_{a+c} \big)
\\
&&
+ \lambda_{1^22} \big(s_{2a} s_{b} s_{c} + s_{b+c} s_{a}^{2} + 2 \sum^{2}
s_{a+b} s_a s_c \big) + \lambda_{1^4} s_{a}^{2} s_b s_c,
\\
E S_a^2 S_b^2
&=&
\lambda_4 s_{2a+2b} + 2\lambda_{13} \big(s_{2a+b} s_b +
 s_{a+2b} s_a \big) + \lambda_{2^2} \big( s_{2a}s_{2b} + 2
s_{a+b}^2  \big)
\\
&&
+ \lambda_{1^22} \big( s_{2a} s_b^2 + s_{2b} s_a^2 + 4 s_{a+b} s_a s_b \big) +
\lambda_{1^4} s_a^2 s_b^2,
\\
E S_a^3 S_b
&=&
\lambda_4 s_{3a+b} + \lambda_{13} \big( s_{3a} s_b + 3 s_{2a+b}
s_a  \big) + 3 \lambda_{2^2} s_{2a} s_{a+b}
\\
&&
+ \lambda_{1^22} \big( 3 s_{2a} s_a s_b + 3 s_{a+b} s_a^2  \big) + \lambda_{1^4} s_a^3 s_b,
\\
E S_a^4 S_b
&=&
\lambda_5 s_{4a+b} + \lambda_{14} \big( 4 s_a s_{3a+b} + s_b
s_{4a} \big) + \lambda_{23} \big( 4s_{a+b} s_{3a} + 6 s_{2a} s_{2a+b} \big)
\\
&&
+ \lambda_{1^23} \big( 4 s_{3a} s_a s_b + 6 s_{2a+b} s_a^2  \big) + \lambda_{12^2}
\big(3 s_{2a}^2  s_b + 12 s_{a+b} s_{2a}  s_a  \big)
\\
&&
+ \lambda_{1^32} \big(6 s_{2a} s_a^2 s_b + 4 s_{a+b} s_a^3  \big) + \lambda_{1^3} s_a^4  s_b,
\\
E S_a^3 S_b^2
&=&
\lambda_5 s_{3 a+2b} + \lambda_{14} \big(3 s_a s_{2a+2b} + 2 s_b
s_{b+3a} \big)
\\
&&
+ \lambda_{23} \big( 3 s_{2a} s_{a+2b} + s_{2b} s_{3b} + 6 s_{a+b} s_{2a+b} \big)
\\
&&
+ \lambda_{1^23} \big(s_{3a}  s_{b}^{2} + 6 s_{2a+b} s_a s_b + 3 s_{a+2b}
s_a^2 \big)
\\
&&
+ \lambda_{12^2} \big( 3 s_{2b} s_{2a} s_a + 6 s_{a+b}^2 s_a + 6 s_{2b}
s_{a+b} s_b  \big)
\\
&&
+ \lambda_{1^32} \big( s_{2b} s_a^3 + 6 s_{a+b} s_a^2 s_b + 3 s_{2a} s_a s_b^2  \big)
 + \lambda_{1^5} s_a^3 s_b^2,
\\
E S_a^3 S_b S_c
&=&
\lambda_5 s_{3a+b+c} + \lambda_{14} \big( 3 s_{2a+b+c} s_a + \sum^2 s_{3a+c} s_b \big)
\\
&&
+ \lambda_{23} \big(s_{3a} s_{b+c} + 3 \sum^2 s_{2a+b} s_{a+c} +3 s_{a+b+c}
s_{2a}  \big)
\\
&&
+ \lambda_{1^23} \big(s_{3a} s_b s_c + 3 \sum^2 s_{2a+b} s_a s_c + 3 s_{a+b+c}
s_a^2 \big)
\\
&&
+ \lambda_{12^2} \big( 3 \sum^2 s_{2a} s_{a+b} s_c + 3 s_{2a} s_{b+c} s_a + 6
s_{a+b} s_{a+c} s_a  \big)
\\
&&
+ \lambda_{1^32} \big(s_{b+c} s_a^3 + 3 \sum^2 s_{a+b} s_c s_a^2 + 3 s_{2a} s_b
s_c \big) + \lambda_{1^5} s_a^3 s_b s_c,
\end{eqnarray*}
\begin{eqnarray*}
E S_a^2 S_b^2 S_c
&=&
\lambda_5 s_{2a+2b+c} + \lambda_{14} \big( s_{2a+2b} s_c +
2 \sum^2 s_{a+c+2b} s_a  \big)
\\
&&
+ \lambda_{23} \big( \sum^2 s_{2a+c} s_{2b} + 2 \sum^2 s_{2a+b} s_{b+c} + 4
s_{a+b+c} s_{a+b} \big)
\\
&&
+ \lambda_{1^23}\big( \sum^2 s_{2b+c}s_a^2+2 \sum^2s_{a+2b}s_as_c+4s_{a+b+c}
s_as_b \big)
\\
&&
+ \lambda_{12^2} \big(s_{2a} s_{2b} s_c + 2\sum^2s_{2a} s_{b+c} s_b + 2
s_{a+b}^2 s_c + 4\sum^2 s_{a+b} s_{a+c} s_b  \big)
\\
&&
+ \lambda_{1^32} \big(\sum^2 s_{2a} s_b^2 s_c + 4 s_{a+b} s_a s_b s_c + 2 \sum^2 s_{a+c} s_a s_b^2 \big) +
\lambda_{1^5} s_a^2 s_b^2 s_c \big),
\\
E S_a^2 S_b S_c S_d
&=&
\lambda_5 s_{2a+b+c+d} + \lambda_{14} \big( \sum^3
s_{2a+b+c} s_d + s_{a+b+c+d} s_a  \big)
\\
&&
+ \lambda_{23}\big(s_{2a}s_{b+c+d}+2 \sum^3s_{a+b}s_{a+c+d+e}+\sum^3s_{b+c}
s_{2a+d} \big)
\\
&&
+ \lambda_{1^23} \big(s_a^2 s_{b+c+d} + 2 \sum^3 s_a s_b s_{a+c+d} + \sum^3 s_b s_c s_{2a+d} \big)
\\
&&
+ \lambda_{12^2} \big( 3 \sum^3 s_a s_{a+b} s_{c+d} + 2 \sum^3 s_b s_{a+c}
s_{a+d} \big)
\\
&&
+ \lambda_{1^32} \big( s_{2a} s_b s_c s_d + 2 \sum^3 s_{a+b} s_a s_c s_d +
\sum^3 s_a^2 s_{b+c} s_d  \big) + \lambda_{1^5} s_a^2 s_b s_c s_d,
\\
E S_{a}^5 S_{b}
&=&
\lambda_6 s_{5a}s_b +
\lambda_{15} \big(5s_{a}s_{4a+b}+s_{b}s_{5a} \big)+
5 \lambda_{24} \big(s_{a+b}s_{4a}+2s_{2a}s_{3a+b} \big)
\\
&&
+10 \lambda_{33}s_{2a+b}s_{3a}+
5 \lambda_{1^24} \big(s_{4a}s_{a}s_{b}+2s_{3a+b}s_{2a} \big)
\\
&&
+10 \lambda_{123} \big(2s_{3a}s_{a+b}s_{a}+3s_{2a+b}s_{2a}s_{a}+s_{3a}s_{2a}s_{b} \big)+
15 \lambda_{2^3}s_{2a}^2s_{a+b}
\\
&&
+10 \lambda_{1^33} \big(s_{3a}s_{a}^2s_{b}+s_{2a+b}s_{a}^3 \big)+
15 \lambda_{1^22^2} \big(2s_{2a}s_{a+b}s_{a}^2+s_{2a}^2s_{b}s_{a} \big)
\\
&&
+5 \lambda_{1^42} \big(2s_{2a}s_{a}^3s_{b}+s_{a+b}s_{a}^4 \big)+
\lambda_{1^6} s_{a}^5s_{b},
\\
E S_{a}^4 S_{b}^{2}
&=&
\lambda_{6}s_{4a+2b}+
2\lambda_{15} \big(s_{a}s_{3a+2b}+s_{b}s_{4a+b} \big)+
\lambda_{24} \big(8s_{a+b}s_{3a+b}
\\
&&
+6s_{2a}s_{2a+2b}+s_{2b}s_{4a} \big)+
2\lambda_{33} \big(3s_{2a+b}^2+2s_{a+2b}s_{3a} \big)
\\
&&
+\lambda_{1^24} \big(8s_{3a+b}s_{a}s_{b}+6s_{2a+2b}s_{a}^2+s_{4a}s_{b}^2 \big)+
4\lambda_{123} \big(6s_{2a+b}s_{a+b}s_{a}
\\
&&
+3s_{a+2b}s_{2a}s_{a}+3s_{2a+b}s_{2a}s_{b}+s_{3a}s_{2b}s_{a}+2s_{3a}s_{a+b}s_{b} \big)
\\
&&
+3\lambda_{2^3} \big(4s_{a+b}^2s_{2a}+s_{2a}^2s_{2b} \big)+
4\lambda_{1^33} \big(3s_{2a+b}s_{a}^2s_{b}+s_{3a}s_{a}s_{b}^2+s_{a+2b}s_{a}^3 \big)
\\
&&
+3\lambda_{1^22^2} \big(s_{2a}^2s_{b}^2
+2s_{2a}s_{2b}s_{a}^2+4s_{a+b}^2s_{a}^2+8s_{a+b}s_{2a}s_{a}s_{b} \big)
\\
&&
+\lambda_{1^42}\big(8s_{a+b}s_{a}^3s_{b}+6s_{2a}s_{a}^2s_{b}^2+s_{2b}s_{a}^4 \big)+
\lambda_{1^6} s_{a}^4s_{b}^2,
\\
E S_{a}^3 S_{b}^{3}
&=&
\lambda_{6} s_{3a+3b}+
3\lambda_{15}\sum^2s_{a}s_{2a+3b}+
3\lambda_{24} \big(s_{a+b}s_{2a+2b}+\sum^2s_{2a}s_{a+3b} \big)
\\
&&
+\lambda_{33} \big(s_{a+2b}s_{2a+b}+s_{3a}s_{3b} \big)+
3 \lambda_{1^24} \big(3s_{2a+2b}s_{a}s_{b}+\sum^2s_{a+3b}s_{a}^2 \big)
\\
&&
+3 \lambda_{123} \big(\sum^2s_{3b}s_{2a}s_{a}+3\sum^2s_{a+2b}s_{2a}s_{b}+6\sum^2s_{a+2b}s_{a+b}s_{a} \big)
\\
&&
+3 \lambda_{2^3} \big(2s_{a+b}^3+3s_{a+b}s_{2a}s_{2b} \big)+
\lambda_{1^33} \big(9\sum^2s_{2a+b}s_{a}s_{b}^2+ \sum^2 s_{3b}s_{a}^3 \big)
\\
&&
+9 \lambda_{1^22^2} \big(s_{2b}s_{2a}s_{a}s_{b}+2s_{a+b}^2s_{a}s_{b}
+s_{a+b}\sum^2s_{2a}s_{b}^2 \big)
\\
&&
+3 \lambda_{1^42} \big(3s_{a+b}s_{a}^2s_{b}^2+\sum^2s_{2a}s_{a}s_{b}^3 \big)+
\lambda_{1^6} s_{a}^3s_{b}^3,
\end{eqnarray*}
\begin{eqnarray*}
E S_{a}^{4} S_ {b} S_{c}
&=&
\lambda_{6}s_{4a+b+c}+
\lambda_{15}\big(4s_{a}s_{3a+b+c}+\sum^2s_{4a+b}s_{c} \big)+
\lambda_{24}\big(4\sum^2s_{a+b}s_{3a+c}
\\
&&
+6s_{2a}s_{2a+b+c}+s_{b+c}s_{4a} \big)+
2 \lambda_{33} \big(2s_{a+b+c}s_{3a}+3s_{2a+b}s_{2a+c} \big)
\\
&&
+\lambda_{1^24} \big(4s_{a}\sum^2s_{3a+b}s_{c}+6s_{2a+b+c}s_{a}^2+s_{4a}
s_{b}s_{c} \big)
\\
&&
+2 \lambda_{123} \big(2s_{3a}\sum^3s_{a+b}s_{c}+6s_{a+b+c}s_{2a}s_{a}
+3s_{2a}\sum^2s_{2a+b}s_{c}
\\
&&
+6s_{a}\sum^2s_{2a+b}s_{a+c} \big)+
3 \lambda_{2^3} \big(4s_{2a}s_{a+b}s_{a+c}+s_{2a}^2s_{b+c} \big)
\\
&&
+2 \lambda_{1^33} \big(2s_{3a}s_{a}s_{b}s_{c}+3s_{a}^2\sum^2s_{2a+b}s_{c}
+2s_{a+b+c}s_{a}^3 \big)
\\
&&
+3 \lambda_{1^22^2} \big(4s_{2a}s_{a}\sum^2s_{a+b}s_{c}+s_{2a}^2s_{b}s_{c}
+2s_{2a}s_{b+c}s_{a}^2+4s_{a+c}s_{a+b}s_{a}^2 \big)
\\
&&
+\lambda_{1^42} \big(6s_{2a}s_{a}^2s_{b}s_{c}+4s_{a}^3\sum^2s_{a+b}s_{c}
+s_{b+c}s_{a}^4 \big)+
\lambda_{1^6} s_{a}^4s_{b}s_{c},
\\
E S_{a}^{3} S_{b}^{2} S_{c}
&=&
\lambda_{6}s_{3a+2b+c}+
\lambda_{15}3s_{a}s_{2a+2b+c}+2s_{b}s_{3a+b+c}+s_{c}s_{3a+2b}
\\
&&
+\lambda_{24} \big(6s_{a+b}s_{2a+b+c}+3s_{a+c}s_{2a+2b}+3s_{2a}s_{a+2b+c}+2s_{b+c}s_{3a+b}
\\
&&
+s_{2b}s_{3a+c} \big)+
\lambda_{33}\big(6s_{a+b+c}s_{2a+b}+3s_{a+2b}s_{2a+c}+s_{3a}s_{2b+c} \big)
\\
&&
+\lambda_{1^24} \big(6s_{2a+b+c}s_{a}s_{b}+3s_{2a+2b}s_{a}s_{c}+3s_{a+2b+c}s_{a}^2+
2s_{3a+b}s_{b}s_{c}
\\
&&
+s_{3a+c}s_{b}^2 \big)+
\lambda_{123} \big(
3s_{2b+c}s_{2a}s_{a}
+3s_{a+2b}s_{2a}s_{c}
+3s_{2a+c}s_{2b}s_{a}
\\
&&
+12s_{a+b+c}s_{a+b}s_{a}
+6s_{2a+b}s_{a+c}s_{b}
+6s_{2a+b}s_{a+b}s_{c}
+6s_{2a+b}s_{b+c}s_{a}
\\
&&
+6s_{a+2b}s_{a+c}s_{a}
+6s_{2a+c}s_{a+b}s_{b}
+6s_{a+b+c}s_{2a}s_{b}
+s_{3a}s_{2b}s_{c}
\\
&&
+2s_{3a}s_{b+c}s_{b} \big) +
3 \lambda_{2^3} \big( 2s_{a+b}^2s_{a+c}+2s_{a+b}s_{2a}s_{b+c}+s_{2a}s_{a+c}s_{2b} \big)
\\
&&
+\lambda_{1^33} \big(6s_{2a+b}s_{a}s_{b}s_{c}+6s_{a+b+c}s_{a}^2s_{b}
+3s_{2a+c}s_{a}s_{b}^2+3s_{a+2b}s_{a}^2s_{c}
\\
&&
+s_{2b+c}s_{a}^3+s_{3a}s_{b}^2s_{c} \big)+
\lambda_{1^22^2} \big(
3s_{a+c}s_{2a}s_{b}^2
+3s_{a+c}s_{2b}s_{a}^2
\\
&&
+3s_{2a}s_{2b}s_{a}s_{c}
+6s_{a+b}^2s_{c}s_{a}
+6s_{b+c}s_{a+b}s_{a}^2
+6s_{a+b}s_{2a}s_{c}s_{b}
\\
&&
+6s_{b+c}s_{2a}s_{a}s_{b}
+12s_{a+c}s_{a+b}s_{a+b} \big)+
\lambda_{1^42} \big(6s_{a+b}s_{a}^2s_{b}s_{c}
\\
&&
+3s_{2a}s_{a}s_{b}^2s_{c}
+3s_{a+c}s_{a}^2s_{b}^2+2s_{b+c}s_{a}^3s_{b}+s_{2b}s_{a}^3s_{c} \big)+
\lambda_{1^6}s_{a}^3s_{b}^2s_{c},
\\
E S_{a}^{2} S_{b}^{2} S_{c}^{2}
&=&
\lambda_{6} s_{2a+2b+2c}+
2 \lambda_{15} \sum^3s_{a}s_{a+2b+2c}+
\lambda_{24} \big(4\sum^3s_{a+b}s_{a+b+2c}
\\
&&
+\sum^3s_{2a}s_{2b+2c} \big)+
2 \lambda_{33} \big(2s_{a+b+c}^2+\sum^3s_{a+2b}s_{a+2c} \big)
\\
&&
+\lambda_{1^24} \big(4\sum^3s_{a+b+2c}s_{a}s_{b}+\sum^3s_{2a+2b}s_{c}^2 \big)+
2 \lambda_{123} \big(4s_{a+b+c}\sum^3s_{a+b}s_{c}
\\
&&
+2\sum^6s_{2a+b}s_{b+c}s_{c}
+\sum^6s_{a+2b}s_{2c}s_{a} \big)+
\lambda_{2^3}\big(2\sum^3s_{a+b}^2s_{2c}
\\
&&
+8s_{a+b}s_{a+c}s_{b+c}+s_{2c}s_{2a}s_{2b} \big)+
2 \lambda_{1^33} \big(4s_{a+b+c}s_{a}s_{b}s_{c}
\\
&&
+\sum^6s_{b+2c}s_{a}^2s_{b} \big)+
\lambda_{1^22^2} \big(8\sum^3s_{b+c}s_{a+b}s_{c}s_{a}+4\sum^3s_{2c}s_{a+b}s_{b}s_{a}
\\
&&
+\sum^3s_{2b}s_{2a}s_{c}^2+2\sum^3s_{a+b}^2s_{c}^2 \big)+
\lambda_{1^42} \big(4\sum^3s_{a+b}s_{a}s_{b}s_{c}^2
\\
&&
+\sum^3s_{2a}s_{b}^2s_{c}^2 \big)+
\lambda_{1^6} s_{a}^2s_{b}^2s_{c}^2,
\end{eqnarray*}
\begin{eqnarray*}
E S_{a}^{2} S_{b}^{2} S_{c} S_{d}
&=&
\lambda_{6}s_{2a+2b+c+d}+
\lambda_{15} \big(2\sum^2s_{a}s_{a+2b+c+d}+\sum^2s_{c}s_{2a+2b+d} \big)
\\
&&
+\lambda_{24}\big(
4s_{a+b}s_{a+b+c+d}
+2\sum^4s_{a+c}s_{a+2b+d}
+\sum^2s_{2a}s_{2b+c+d}
\\
&&
+s_{c+d}s_{2a+2b} \big)+
\lambda_{33} \big(4s_{a+b+c}s_{a+b+d}+2 \sum^2 s_{a+2b}s_{a+c+d}
\\
&&
+ \sum^2 s_{2a+c}s_{2b+d} \big)+
\lambda_{1^24}\big(4s_{a+b+c+d}s_{a}s_{b}
+2 \sum^4 s_{a}s_{c}s_{a+2b+d}
\\
&&
+ \sum^2 s_{a}^2s_{2b+c+d}s_{a}^2+s_{2a+2b}s_{c}s_{d} \big)
+\lambda_{123} \big(
\sum^4 s_{2a}s_{c}s_{2b+d}
\\
&&
+2\sum^4 s_{a}s_{2b+c}s_{a+d}
+4s_{a+b}\sum^2 s_{c}s_{a+b+d}
+4\sum^4s_{a}s_{b+c}s_{a+b+d}
\\
&&
+2\sum^4s_{c}s_{a+2b}s_{a+d}
+2s_{c+d}\sum^2s_{a}s_{a+2b}
+4s_{a+b}\sum^2s_{a}s_{b+c+d}
\\
&&
+2\sum^2s_{a}s_{2b}s_{a+c+d} \big)+
\lambda_{2^3} \big(2s_{a+b}^2s_{c+d}
+4\sum^2s_{a+b}s_{a+c}s_{b+d}
\\
&&
+2\sum^2 s_{2a}s_{b+c}s_{b+d}+s_{2a}s_{2b}s_{c+d} \big)
+\lambda_{1^33}\big(4 \sum^2 s_{a}s_{b}s_{c}s_{a+b+d}
\\
&&
+2\sum^2s_{a}^2s_{b}s_{b+c+d}
+2 \sum^2 s_{a}s_{c}s_{d}s_{a+2b}
+ \sum^4 s_{a}^2s_{c}s_{2b+d} \big)
\\
&&
+\lambda_{1^22^2}\big(
2 \sum^4 s_{2a}s_{c}s_{b}s_{b+d}
+2\sum^2  s_{a}^2 s_{b+c} s_{b+d}
\\
&&
+s_{c+d}\sum^2 s_{a}^2  s_{2b}
+4 \sum^4 s_{a+b}s_{b+d}s_{a}s_{c}
 +4\sum^2s_{a+d}s_{b+c}s_{a}s_{b}
\\
&&
+4s_{a}s_{b}s_{c+d}s_{a+b}
+2 s_{c}s_{d}  s_{a+b}^2
+s_{2a}s_{2b}s_{c}s_{d} \big)
\\
&&
+\lambda_{1^42} \big(4s_{a+b}s_{a}s_{b}s_{c}s_{d}
+2 \sum^4s_{a}^2s_{b}s_{c}s_{b+d}
 +s_{c+d}s_{a}^2s_{b}^2
\\
&&
+s_{c}s_{d} \sum^2 s_{2a}s_{b}^2 \big)+\lambda_{1^6}  s_{a}^2s_{b}^2s_{c}s_{d},
\end{eqnarray*}
\begin{eqnarray*}
E S_{a}^{3} S_{b} S_{c} S_{d}
&=&
\lambda_{6} s_{3a+b+c+d}+
\lambda_{15} \big(3s_{a}s_{2a+b+c+d}+\sum^3s_{b}s_{3a+c+d} \big)
\\
&&
+\lambda_{24} \big(3\sum^3s_{a+b}s_{2a+c+d}+3s_{2a}s_{a+b+c+d}+\sum^3s_{b+c}s_{3a+d} \big)
\\
&&
+\lambda_{33} \big(3 \sum^3 s_{2a+b}s_{a+c+d}
+s_{3a}s_{b+c+d} \big)+
\lambda_{1^24} \big(3 s_{a} \sum^3 s_{b}s_{2a+c+d}
\\
&&
+3s_{a}^2s_{a+b+c+d}
+ \sum^3 s_{b}s_{c}s_{3a+d} \big)+
\lambda_{123} \big(
s_{3a} \sum^3 s_{b}s_{c+d}
+3s_{a}s_{2a}s_{b+c+d}
\\
&&+3s_{2a}\sum^3s_{b}s_{a+c+d}
+6s_{a}\sum^3s_{a+b}s_{a+c+d}
+3s_{a}\sum^3s_{b+c}s_{2a+d}
\\
&&
+3\sum^6s_{b}s_{a+c}s_{2a+d} \big)+
3 \lambda_{2^3} \big( s_{2a}\sum^3s_{a+b}s_{c+d}+2s_{a+d}s_{a+b}s_{a+c} \big)
\\
&&
+\lambda_{1^33} \big(3s_{a}\sum^3s_{2a+b}s_{c}s_{d}
+3s_{a}^2\sum^3s_{a+b+c}s_{d}
+s_{b+c+d}s_{a}^3
\\
&&
+s_{3a}s_{b}s_{c}s_{d} \big)+
3 \lambda_{1^22^2} \big(
s_{2a} \sum^3 s_{b}s_{c}s_{a+d}
+2s_{a}\sum^3s_{b}s_{a+c}s_{a+d}
\\
&&
+s_{a}^2\sum^3s_{a+b}s_{c+d}
+s_{a}s_{2a}\sum^3s_{b}s_{c+d} \big)+
\lambda_{1^42} \big(3s_{2a}s_{a}s_{b}s_{c}s_{d}
\\
&&
+3s_{a}^2\sum^3s_{b}s_{c}s_{a+d}
+s_{a}^3\sum^3s_{b}s_{c+d} \big)+
\lambda_{1^6}  s_{a}^3s_{b}s_{c}s_{d},
\\
E S_{a}^{2} S_{b} S_{c} S_{d} S_{e}
&=&
\lambda_{6} s_{2a+b+c+d+e}+
\lambda_{15} \big(2s_{a}s_{a+b+c+d+e}
+\sum^{4}s_{b}s_{2a+c+d+e} \big)
\\
&&
+\lambda_{24} \big(2 \sum^4 s_{a+b}s_{a+c+d+e}+s_{2a}s_{b+c+d+e}+\sum^{6}s_{b+c}s_{2a+d+e} \big)
\\
&&
+\lambda_{33}\big(2 \sum^3 s_{a+b+c}s_{a+d+e}
+\sum^{4}s_{2a+b}s_{c+d+e} \big)
\\
&&
+\lambda_{1^24} \big(2 s_{a} \sum^4 s_{b}s_{a+c+d+e}
+s_{b+c+d+e}s_{a}^2
+\sum^{6}s_{2a+d+e}s_{b}s_{c} \big)
\\
&&
+\lambda_{123}
\big(2 \sum^{12}s_{b}s_{a+c}s_{a+d+e}
+s_{2a} \sum^4 s_{b} s_{c+d+e}
+2 s_{a} \sum^6 s_{b+c}s_{a+d+e}
\\
&&
+2 s_{a} \sum^4 s_{a+b}s_{c+d+e}
+ \sum^{12} s_{b}s_{c+d}s_{2a+e} \big)+
\lambda_{2^3} \big( 2\sum^6s_{a+b}s_{a+c}s_{d+e}
\\
&&
+s_{2a}\sum^3s_{b+c}s_{d+e} \big)+
\lambda_{1^33}\big(
2s_{a}\sum^6s_{b}s_{c}s_{a+b+c}
+s_a^2\sum^4s_{b}s_{c+d+e}
\\
&&
+\sum^4s_{2a+b}s_{c}s_{d}s_{e} \big)+
\lambda_{1^22^2} \big(
2 \sum^6 s_{b}s_{c}s_{a+d}s_{a+e}
+ s_{2a} \sum^6 s_{b}s_{c}s_{d+e}
\\
&&
+ s_{a}^2 \sum^3 s_{b+c}s_{d+e}
+ 2s_{a} \sum^{12} s_{b}s_{a+c}s_{d+e} \big)+
\lambda_{1^42} \big( 2s_{a}\sum^4s_{a+b}s_{c}s_{d}s_{e}
\\
&&
+s_{a}^2\sum^6s_{b+c}s_{d}s_{e}
+s_{2a}s_{b}s_{c}s_{d}s_{e} \big)+
\lambda_{1^6} s_{2a}s_{b}s_{c}s_{d}s_{e}.
\end{eqnarray*}


\begin{thebibliography}{999}


\bibitem{}
Carver, H. C. (1930).
Fundamentals of the theory of sampling.
Annals of Mathematical Statistics, 1, 101-121.


\bibitem{}
Carver, H. C. (1930).
Fundamentals of the theory of sampling.
Annals of Mathematical Statistics, 1, 260-274.

\bibitem{}
Dwyer, P. S. and Tracy, D. S. (1980).
Expectation and estimation of
product moments in sampling from a finite population.
Journal of the American Statistical Association, 75, 431-437.

\bibitem{}
Dwyer, P. S. (1938).
On combined expansions of products of symmetric
power sums and of sums of symmetric power products with applications to sampling (continued).
Annals of Mathematical Statistics, 9, 97-132.

\bibitem{}
Fisher, R. A. (1929).
Moments and product moments of sampling distributions.
Proceedings of the London Mathematical Society Series 2, 30, 199-238.

\bibitem{}
James, G. S. (1958).
On moments and cumulants of systems of statistics.
Sankhy\=a, 20, 1-30.

\bibitem{}
Nath, S. N. (1968).
On product moments from a finite universe.
Journal of the American Statistical Association, 63, 535-541.

\bibitem{}
Nath, S. N. (1969).
More results on product moments from a finite universe.
Journal of the American Statistical Association, 64, 864-869.

\bibitem{}
Pierce, J. A. (1940).
A study of a universe of $n$ finite populations
with application to moment-function adjustments for grouped data.
Annals of Mathematical Statistics, 11, 311-334.

\bibitem{}
Raghunandanan, K. and Srinivasan, R. (1973).
Some product moments useful in sampling theory.
Journal of the American Statistical Association, 68, 409-413.

\bibitem{}
Skellam, J. G. (1949).
The distribution of moment statistics of samples
drawn without replacement from a finite population.
Journal of the Royal Statistical Society, B, 11, 291-296.

\bibitem{}
Stuart, A. and Ord, J. K. (1987).
Kendall's Advanced Theory of Statistics, fifth edition, volume 1.
Griffin, London.

\bibitem{}
Sukhatme, P. V. (1944).
Moments and product moments of moment-statistics
for samples of the finite and infinite populations.
Sankhy\=a, 6, 363-382.


\bibitem{}
Tiit, E. (1988).
Unbiased and $n^{-k}$-biased estimations of entire rational functions of moments.
Statistical and Probabilistic Models, Report 798, Tartu State University, Tartu, Estonia, pp. 3-17.

\bibitem{}
Wishart, J. (1952).
Moment coefficients of the $k$-statistics in samples from a finite population.
Biometrika, 39, 1-13.


\end{thebibliography}
\end{document}